\documentclass[10pt]{article}
\usepackage{amsmath}
\usepackage{amsthm}
\usepackage{graphicx}
\usepackage{latexsym}
\usepackage{amssymb}
\usepackage{amsbsy}
\usepackage[]{authblk}

\DeclareSymbolFont{msbm}{U}{msb}{m}{n}
\DeclareMathSymbol{\C}{\mathalpha}{msbm}{'103}
\DeclareMathSymbol{\R}{\mathalpha}{msbm}{'122}
\DeclareMathSymbol{\Z}{\mathalpha}{msbm}{'132}
\DeclareMathSymbol{\N}{\mathalpha}{msbm}{'116}

\newtheorem{remark}{Remark}

\def\vr{\boldsymbol{\varrho}}
\def\RR{\mathbb R}

\def\be{\begin{equation}}
\def\ee{\end{equation}}
\def\bea{\begin{eqnarray}}
\def\ba{\begin{array}{l}\displaystyle}
\def\eea{\end{eqnarray}}
\def\ea{\end{array}}

\def\IR{\mathop{\mbox{\rm Iround}}\nolimits}

\begin{document}

\title{\small \hfill 
    \Large
\begin{center}
FLUID SIMULATIONS WITH LOCALIZED BOLTZMANN UPSCALING BY DIRECT
SIMULATION MONTE-CARLO
\thanks{Acknowledgements: This work was supported by the Marie Curie Actions of the European
Commission in the frame of the DEASE project (MEST-CT-2005-021122)}
\end{center}}
\author[1,2]{Pierre Degond}
\author[1,2]{Giacomo Dimarco\footnote
{Corresponding author address: Institut de Math\'{e}matiques de
Toulouse, UMR 5219 Universit\'{e} Paul Sabatier, 118, route de
Narbonne 31062 TOULOUSE Cedex, FRANCE. \\
\emph{E-mail}:pierre.degond@math.univ-toulouse.fr,giacomo.dimarco@math.univ-toulouse.fr}}

\affil[1]{Universit\'{e} de Toulouse; UPS, INSA, UT1, UTM ;

Institut de Math\'{e}matiques de Toulouse ; F-31062 Toulouse,
France.} \affil[2]{CNRS; Institut de Math\'{e}matiques de Toulouse
UMR 5219;

F-31062 Toulouse, France.} \maketitle

\begin{abstract}
In the present work, we present a novel numerical algorithm to
couple the Direct Simulation Monte Carlo method (DSMC) for the
solution of the Boltzmann equation with a finite volume like method
for the solution of the Euler equations. Recently we presented in
\cite{degond2},\cite{dimarco3},\cite{dimarco4} different
methodologies which permit to solve fluid dynamics problems with
localized regions of departure from thermodynamical equilibrium. The
methods rely on the introduction of buffer zones which realize a
smooth transition between the kinetic and the fluid regions. In this
paper we extend the idea of buffer zones and dynamic coupling to the
case of the Monte Carlo methods. To facilitate the coupling and
avoid the onset of spurious oscillations in the fluid regions which
are consequences of the coupling with a stochastic numerical scheme,
we use a new technique which permits to reduce the variance of the
particle methods \cite{dimarco1}. In addition, the use of this
method permits to obtain estimations of the breakdowns of the fluid
models less affected by fluctuations and consequently to reduce the
kinetic regions and optimize the coupling. In the last part of the
paper several numerical examples are presented to validate the
method and measure its computational performances.

\end{abstract}

 \vskip 0.4cm

\setcounter{equation}{0}
\section{Introduction}
\label{intro}

Many problems of interests in applications involve non equilibrium
gas flows as hypersonic objects simulations or
micro-electro-mechanical devices. Often, these kinds of problems are
characterized by breakdowns of fluid models, either Euler or
Navier-Stokes. Now, when the breakdown is localized both in space
and time we must deal with conjunctions of continuum and non
equilibrium regions. To face such problems, the most natural
approach is to try to combine numerical schemes for continuum models
with microscopic kinetic models which guarantee a more accurate
description of the physics when far from thermodynamical equilibrium
conditions are reached. In fact, it is not always necessary to solve
the Boltzmann equation, which is computationally more expensive by
several orders of magnitude than a continuum description, in the
entire domain. It is, instead, sufficient to use the microscopic
description in regions where the departure from equilibrium is
strong and manifests itself as boundary layers or shocks.

The construction of such hybrid kinetic-fluid methods involve three
main problems. The first one is how to accurately identify the
different regions. We refer, for instance, to the works of
Wijesinghe and Hadjiconstantinou~\cite{he}, Levermore, Morokoff, and
Nadiga~\cite{Levermore}, and Wang and Boyd~\cite{Wang}, in which
various breakdown criteria are analyzed and then proposed. The
second main problem is how to, in practice, realize the coupling. In
other words, how to match the two models at the interfaces. Several
different ideas are described in the works of Bourgat, LeTallec,
Perthame, and Qiu \cite{BLPQ}, Bourgat, LeTallec and Tidriri
\cite{BLT}, LeTallec and Mallinger \cite{Letallec}, Aktas and Aluru
\cite{AA}, Roveda, Goldstein and Varghese \cite{RGV1,RGV2}, Sun,
Boyd and Candler \cite{Boyd2,WSB}, Wadsworth and Erwin \cite{Wad}.
Domain decomposition approaches are also popular in others fields
such as, for instance in epitaxial growth \cite{SSE} or for problems
involving different scaling such as the diffusive scalings
\cite{degond1,Klar}. The decomposition between equilibrium and non
equilibrium states can be realized also in the velocity space
instead of physical space \cite{degond}, \cite{dimarco5}. It is
important to stress that most of the mentioned methods use a static
interface between kinetic and fluid regions. However, this approach
appears as somehow inadequate and inefficient for many realistic
problems. To overcome this difficulty, automatic domain
decomposition methods have also been developed in the recent past.
See for example Kolobov et al.~\cite{Kobolov}, or
Tiwari~\cite{tiwari_JCP,tiwari_JCP1} and Degond, Dimarco and
Mieussens~\cite{dimarco3,dimarco4}. The third problem is the choice
of the appropriate numerical methods for solving the macroscopic and
the microscopic model. A first possibility are the
deterministic-deterministic methods in which both the macroscopic
and the microscopic models are solved by means of finite volume
techniques such as for instance in \cite{degond2,dimarco3,dimarco4}.
These methods permit to obtain very accurate results but generally
they are too expensive. For this reason simplified models are then
used to describe the kinetic part of the problem. In practice, in
most of the cases, the collision operator is replaced by the BGK
relaxation operator which permits to reduce the computational cost.
A second approach consists in using particle-particle methods
\cite{TK,Boyd2}. These methods avoid the complexity associated with
using two independent methods for different regions and provide a
simple solution to the interface problem. On the other hand, the
solutions contain large statistical errors which are the typical
drawback of Monte Carlo methods. Recently, some low diffusion
particle-particle numerical methods have been proposed by Boyd and
coauthors in \cite{Boyd}. A third approach is the
deterministic-–particle approach
\cite{Schwartz,Wad,Wij,HashHassan,WSB}. Here, the fluid model is
discretized by using a finite volume method, while the kinetic model
is solved by the DSMC method. The popularity of DSMC is due to
several advantages compared to other simulation methods for kinetic
equations, including advantages related to efficiency, memory usage
and implementation of additional physical phenomena. The
disadvantages of DSMC, however, are the increase of the
computational cost close to thermodynamical equilibrium and the
large statistical error which scales with the inverse square root of
the number of particles. This last aspect is the cause of the big
difficulties met in the coupling, both for defining the kinetic
regions and to avoid the introduction of oscillations in the fluid
part.

In the present paper, we present a novel numerical algorithm to
couple a DSMC method for the solution of the Boltzmann model and a
finite volume like method for the Euler equations. The method relies
on two main aspects. The first one is the use of a new low variance
Direct Simulation Monte Carlo method to compute the solution of the
kinetic part. Monte Carlo methods, as already observed, permit more
realistic descriptions of the physical problems. On the other hand,
the solutions computed with these techniques are affected by large
statistical noise. In order to reduce this effect, we extend the
variance reduction technique (called moment guided method) which has
been recently developed in \cite{dimarco2} to the coupling case. The
moment guided method permits to reduce the statistical error through
a matching with a set of suitable macroscopic moment equations. The
basic idea, on which the method relies, consists in guiding the
particle positions and velocities through moment equations so that
the concurrent solution of the moment and kinetic models furnishes
the same macroscopic quantities. The method is based on the
observation that the deterministic solution of the moment equation
leads to a more accurate solution, in term of statistical
fluctuations, than the DSMC method. We experimentally showed in
\cite{dimarco1} that this is indeed the case. Now, a crucial point,
in the implementation of the domain decomposition strategies, is the
identification of the zones in which the microscopic description is
necessary. In other words, kinetic regions have to be as thin as
possible, with the constraint that the effective solution is
correctly captured. Thus, the direct consequence of using the moment
guided method to solve the Boltzmann equation is that the research
of the interface location becomes a less difficult task to
accomplish. Moreover, the use of this technique permits to reduce
significantly the propagation of spurious oscillations in the fluid
regions.

The second main aspect of the method consists of the introduction of
a buffer zone which realizes a smooth transition between the kinetic
and fluid regions. It differs from the method of
\cite{degond2,degond3,dimarco3,dimarco4} by the way the solution of
the Boltzmann equation is decomposed into a kinetic and a fluid
component. The new decomposition we adopt is better adapted to a
DSMC method while the earlier one was more adapted to grid based
solutions of the Boltzmann equation. As in \cite{dimarco1}, we
suggest a methodology to allow for the time evolution of the kinetic
and fluid regions which permits to follow discontinuities in time
and space and solve the microscopic model only where necessary.
Thanks to this technique, it is possible to achieve considerable
computational speedup compared with steady interface coupling
strategies, without otherwise, losing accuracy in the solution.
Finally, we observe that the use of a smooth transition between the
micro and the macro models further reduces the propagation of
fluctuations in the fluid regions and thus it permits to obtain more
accurate results.

The outline of the article is the following. In section
\ref{sec_Boltzmann}, we introduce the Boltzmann equation, its
properties and the fluid limit. In section \ref{sec moment_guided}
we describe the moment guided method while in section
\ref{sec_coupling} we present the coupling strategy. Section
\ref{sec_num_approx} is devoted to the discretization of the
Boltzmann equation and of the finite volume scheme for the Euler
equations. In section \ref{sec moving}, we describe the breakdown
criterion. Several numerical tests are presented in section
\ref{sec_tests} to illustrate the properties of our method and to
demonstrate its efficiency in comparison with DSMC schemes. A short
conclusion is given in section \ref{sec_conclu}.

\setcounter{equation}{0}
\section{The Boltzmann equation and its fluid limit}
\label{sec_Boltzmann} We consider the following model
\cite{cercignani}  \be
\partial_t f + v\cdot\nabla_{x}f = Q(f,f),
\label{eq:B} \ee with the initial data
\begin{equation*}
  f(x,v,t=0)=f_{0},
\end{equation*}
where $f=f(x,v,t)$ is the density of particles that have velocity $v
\in \R^3 $ and position $x \in \Omega \subset \R^3$ at time $ t> 0$.
The collision operator $Q$ locally acts in space and time and takes
into account interactions between particles. It is written: \be
Q(f,f)=\int_{\R^3}\int_{S^2} B\left(|q|,\frac{q\cdot n}{|q|}\right)
[f(v')f(v'_*)-f(v)f(v_*)]dn dv_*.\ee In the above expression $S^2$
is the unit sphere in $\R^3$, $q=v-v_*$ the relative velocity, $n\in
S^2$ the unit vector in the direction of the scattered velocity. The
collision kernel $B(|q|,q\cdot n/|q|)$ characterizes the detail of
the collision and is defined as {\be
B(|q|,\cos\theta)=|q|\sigma(|q|,\theta), \ (0 \leq \theta \leq \pi)
\ee} where $\cos \theta=q\cdot n/|q|$ and $\sigma(|q|,\theta)$ is
the collision cross section which corresponds to the scattering
angle $\theta$. Finally the post collisional velocity are computed
by {\be v'=\frac{1}{2}(v+v_*+|q|n), \ v'_*=\frac{1}{2}(v+v_*-|q|n).
\ee} The collisional operator locally satisfies the conservation of
mass, momentum and energy: \be \langle m Q(f,f)\rangle=0\ee for
every $f$, where we denote weighted integrals of $f$ over the
velocity space by \be\langle \phi f\rangle=\int_{\R^3}
\phi(v)f(v)dv,\ee with $\phi(v)$ any function of $v$, and
$m(v)=(1,v,|v|^2)$ are the collisional invariants. The
multiplication of (\ref{eq:B}) by $m(v)$ and the integration in
velocity space leads to the following system of local conservation
laws \be
\partial_t \langle m f\rangle +\nabla_{x}\cdot\langle v m f\rangle = 0. \label{consl}\ee
Now, the well known Boltzmann's H-theorem implies that any
equilibrium distribution function, i.e. any function $f$ for which
$Q(f, f)=0$, has the form of a locally Maxwellian distribution: \be
 E[\boldsymbol{\varrho]}(v)=\frac{\varrho}{(2\pi \theta)^{3/2}}\exp\left(\frac{-|u-v|^{2}}{2\theta}\right) ,
\label{eq:M} \ee where $\boldsymbol{\varrho}=(\varrho,\varrho
u,\varrho e)$, with $\varrho$ and $u$ the density and mean velocity
while $\theta=RT$, with $T$ the temperature of the gas and $R$ the
gas constant. The macroscopic values $\varrho$, $u$ and $T$ are
related to $f$ by: \be \varrho=\int_{\R^3} fdv, \qquad \varrho
u=\int_{\R^3} vfdv, \qquad \theta=\frac{1}{3\varrho}
\int_{\R^3}|v-u|^{2}fdv, \label{eq:Mo} \ee while the fluid energy
$e$ is defined as \be e=\frac{1}{2\varrho} \int_{\R^3}|v|^{2}fdv =
\frac{1}{2}|u|^2 + \frac{3}{2}  \theta. \label{eq:E} \ee

When the mean free path between particles is very small compared to
the size of the computational domain, the space and time variables
can be rescaled to \be x'=\varepsilon x, \ t'=\varepsilon
t\label{eq:scaling}\ee where $\varepsilon$ is the ratio between the
microscopic and the macroscopic scale (the so-called Knudsen
number). Using these new variables in~(\ref{eq:B}), we get \be
\partial_{t'} f^\varepsilon + v\cdot\nabla_{x'}f^\varepsilon =
\frac{1}{\varepsilon}Q(f^{\varepsilon},f^{\varepsilon}).\ee If the
Knudsen number $\varepsilon$ tends to zero, formally the
distribution function converges towards the local Maxwellian
equilibrium $E^\varepsilon[\boldsymbol{\varrho}]$. Inserting this
relation into the conservation laws~(\ref{consl}) gives the set of
compressible Euler equations for the macroscopic quantities
$\boldsymbol{\varrho}$: \be
\partial_{t'} \boldsymbol{\varrho}
+\nabla_{x'}\cdot F(\boldsymbol{\varrho}) = 0,\ee where
$F(\boldsymbol{\varrho})=\langle v m
E^{\varepsilon}[\boldsymbol{\varrho}]\rangle$.

In the sequel we will omit the primes to simplify notations wherever
they will not cause in confusion.

\section{The moment guided method}
\label{sec moment_guided}

In this section we introduce the method which permits to reduce
statistical fluctuations in the DSMC scheme. We refer to
\cite{dimarco1, dimarco2} for details on the method. The starting
point of the moment guided method consists in decomposing the
distribution function $f$ as a part in equilibrium and a part in non
equilibrium. We call this decomposition micro-macro decomposition
\be f=E[\vr]+g.\ee

The function $g$ represents the non-equilibrium part of the
distribution function. From the definition above, it follows that
$g$ is in general non positive. Moreover since $f$ and $E[\vr]$ have
the same moments we have \be \langle m(v)g\rangle=0.\ee Finally, it
is possible to show that $E[\vr]$ and $g$ satisfy the following
coupled system of equations
\begin{eqnarray}
& &\partial_t \boldsymbol{\varrho}+\nabla_x \cdot F
(\boldsymbol{\varrho})+\nabla_x\cdot\langle vm(v)g\rangle=0,\label{momentg}\\
& &\partial_t f + v \cdot \nabla_x f =\frac{1}{\varepsilon} Q(f,f)
\label{momentg1}
\end{eqnarray}  We skip the proof of the above statement and
we refer to \cite{degond3} for details on the micro-macro
decomposition and its properties.

The idea of the moment guided strategy is to solve the kinetic
equation (\ref{momentg1}) with the Monte Carlo method, and
concurrently the fluid equation with a finite volume scheme. The
correction term which is necessary to close the macroscopic
equations $\nabla_x\cdot\langle vm(v)g\rangle$ is evaluated using
particle moments. We observe that the two equations
(\ref{momentg}-\ref{momentg1}), except for numerical errors, give
the same results in terms of macroscopic quantities. We assume that
the set of moments obtained from the fluid system represents a
better statistical estimate of the true moments of the solution,
since the resolution of the moment equations does not involve any
stochastic process. We experimentally showed that this is the case
in \cite{dimarco1}.

Consider a time discretization of the problem
(\ref{momentg}-\ref{momentg1}), then the method is summarized in the
following steps: at each time step $t^n$
\begin{enumerate}
 \item Solve the kinetic equation (\ref{momentg1}) with any type of DSMC method and obtain a first
 set of moments.
 \item Solve the fluid equation (\ref{momentg}) with any type of finite volume scheme where
 particles are used to evaluate $\nabla_x\cdot\langle vm(v)g\rangle$ and obtain a second set of moments.
 \item Match the moments of the kinetic solution with the fluid solution through a transformation of the
 samples values.
 \item Restart the computation to the next time step.
 \end{enumerate}

In the following we will consider the case in which the perturbation
$g$ is very small in a part of the domain while it is large in the
complementary part. In the small perturbation regions it will be
suitable to avoid the solution of the kinetic equation
(\ref{momentg1}) and solve only the Euler equations which are able
to furnish an accurate description of the problem under
consideration. However, the direct passage from the micro to the
macro model does not prevent the fluctuations, which has been
reduced with the moment guided method but which are not completely
disappeared, to propagate to the fluid regions. In order to avoid
this loss of accuracy, we introduce in the next section a model
which realizes a smooth transition between the different regions.

\section{The coupling method}
\label{sec_coupling}  Let $\Omega_1$, $\Omega_2$, and $\Omega_3$ be
three disjointed sets such that $\Omega_1\cup \Omega_2 \cup
\Omega_3=\R^3$. The first set $\Omega_1$ is supposed to be a domain
in which the flow is far from the equilibrium (the "kinetic zone",
$g$ large), while the flow is supposed to be close to equilibrium in
$\Omega_2$ (the "fluid zone", $g\simeq 0$). The third region is the
buffer zone and also in $\Omega_3$ we suppose the gas to be close to
thermodynamical equilibrium. We define a function $h(x,t)$ such that
\begin{equation}
h(x,t)=\left\lbrace
\begin{array}{lll}
\displaystyle 1, & \mbox{for} & x \in \Omega_1, \\
0, & \mbox{for} & x \in \Omega_2, \\
0 \leq h(x,t) \leq 1, & \mbox{for} & x \in \Omega_3. \\
\end{array}
\right.
\end{equation}
Note that the time dependence of $h$ means that we account for
dynamically changing the shape of the transition function. The
topology and geometry of the different zones is directly encoded in
$h$ and may change dynamically as well.

Multiplying (\ref{momentg1}) by $h$ and $1-h$, we get
\begin{eqnarray}
& &\partial_t \boldsymbol{\varrho}+\nabla_x \cdot F
(\boldsymbol{\varrho})+\nabla_x\cdot\langle hvm(v)g\rangle+\nabla_x\cdot\langle (1-h)vm(v)g\rangle= 0, \label{hrho} \\
& & \partial_t (hf) + v \cdot \nabla_x (hf)= \frac{h}{\varepsilon} Q(f,f) + \left(\frac{\partial_t h}{h} + v \cdot \frac{\nabla_x h}{h} \right) (hf), \label{hf} \\
& &\partial_t ((1-h)f) + v \cdot \nabla_x ((1-h)f)=
\frac{1-h}{\varepsilon} Q(f,f) - \left(\frac{\partial_t h}{h}+ v
\cdot \frac{\nabla_x h}{h} \right) (hf). \label{(1-h)f}
\end{eqnarray}
Denoting $f_K = hf$, $f_F = (1-h) f$, $g_K=hg$ and $g_F=(1-h)g$ the
ensemble $(f_K,f_F,\boldsymbol{\varrho},g_K,g_F)$ satisfies the
system
\begin{eqnarray}
& &\partial_t \boldsymbol{\varrho}+\nabla_x \cdot F
(\boldsymbol{\varrho})+\nabla_x\cdot\langle
vm(v)g_K\rangle+\nabla_x\cdot\langle vm(v)g_F\rangle=0, \label{rho_K} \\
& & \partial_t f_K + v \cdot \nabla_x f_K=
\frac{h}{\varepsilon}Q(f_K+f_F,f_K+f_F) + \left(\frac{\partial_t h+
v \cdot \nabla_x h
}{h}\right) f_K, \label{fK} \\
& & \partial_t f_F + v \cdot \nabla_x f_F =
\frac{1-h}{\varepsilon}Q(f_K+f_F,f_K+f_F) - \left(\frac{\partial_t
h+ v \cdot \nabla_x h }{h}\right) f_K , \label{fF}
\end{eqnarray}
We observe that the system (\ref{rho_K})-(\ref{fF}) is equivalent to
equation (\ref{eq:B}) in the sense that if $f$ is a solution of
(\ref{eq:B}) with initial datum $f(t=0) = f_0, \
g(t=0)=g_0=f(t=0)-E[\vr(t=0)]$ then $f_K = hf$, $f_F = (1-h) f$,
$g_K=hg$, $g_F=(1-h)g$ and $\boldsymbol{\varrho}=\langle
m(v)f\rangle$ are solutions of system (\ref{rho_K})-(\ref{fF}) with
initial data $f_K(t=0) = hf_0$, $f_F(t=0) = (1-h)f_0$, $g_K(t=0) =
hg_0$, $g_F(t=0) = (1-h)g_0$ $\boldsymbol{\varrho}(t=0)=\langle
m(v)f_0\rangle$. Reciprocally, if $f_K$, $f_F$, $g_K$, $g_F$ and
$\boldsymbol{\varrho}$ are solutions of (\ref{rho_K})-(\ref{fF})
with initial data $f_K(t=0) = hf_0$, $f_F(t=0) = (1-h)f_0$,
$g_K=hg_0$, $g_F=(1-h)g_0$ and $\boldsymbol{\varrho}(t=0)=\langle
m(v)f_0\rangle$ then $f=f_K + f_F, \ g=g_K+g_F=f-E[\vr]$ is the
solution of (\ref{eq:B}) with initial data $f(t=0) = f_0$ and
$g(t=0)=g_0$.

Now assume that the flow is very close to equilibrium in
$\Omega_2\cup\Omega_3$. This means that we assume the distribution
function to be Maxwellian in this set: $f_F =
E[\boldsymbol{\varrho}_F]$ and thus $g_F=0$. Then, equation
(\ref{fF}) can be eliminated and we get:
\begin{eqnarray}
& &\partial_t \boldsymbol{\varrho}+\nabla_x \cdot F
(\boldsymbol{\varrho})+\nabla_x\cdot\langle
vm(v)g_K\rangle=0, \label{rho_K1} \\
& & \partial_t f_K + v \cdot \nabla_x f_K=
\frac{h}{\varepsilon}Q(f_K+E[\vr_F],f_K+E[\vr_F]) +
\left(\frac{\partial_t h+ v \cdot \nabla_x h }{h}\right) f_K.
\label{fK_f1}
\end{eqnarray}
 System~(\ref{rho_K1}--\ref{fK_f1})
represents our kinetic-fluid model. In this model the local
distribution function $f$ is approximated by $E[\vr]+g_K$ where
$g_K=f_K-E[\vr_K]$. The distribution $E[\vr_K]$ is defined as the
Maxwellian whose moments are given by $\langle m(v)f_K\rangle$,
while $E[\vr_F]$ is the Maxwellian whose moments are given by
$\vr-\langle m(v)f_K\rangle$. Observe that $f_K$ is zero in the
fluid zone $\Omega_2$ as well as $g_K$. This means that, in these
regions, we simply solve the compressible Euler equations. On the
other hand, in $\Omega_1$ the function $\langle m(v)f_K\rangle=\vr$
and thus the system (\ref{rho_K1}--\ref{fK_f1}) reduces to the
system (\ref{momentg}--\ref{momentg1}) which is equivalent to the
Boltzmann equation (\ref{eq:B}) where we employ the moment guided
method to reduce the variance of the DSMC method.

In the next section we will discuss the numerical scheme which
discretizes the system (\ref{rho_K1}--\ref{fK_f1}) as well as the
boundary conditions between the different regions.

\section{The numerical schemes}
\label{sec_num_approx}

\subsection{The DSMC method for the coupled Boltzmann equation}

The classical approach used to solve the Boltzmann equation by Monte
Carlo methods is the time splitting approach \cite{Babovsky,Nanbu80}
between the free transport term \be
\partial_t f_K + v\cdot\nabla_{x}f_K
=0,\label{split1}\ee and the collision term \be
\partial_t f_K=\frac{h}{\varepsilon}Q(f_K+E[\vr_F],f_K+E[\vr_F]).\label{split2}\ee
We start to discuss the discretization of the transport term which,
in our model, is replaced by the free transport term plus the term
which involves the time and space derivative of the transition
function $h$ ($h$-term).
\subsubsection{Free transport step} We
introduce a space discretization of mesh size $\Delta x$ and a time
discretization of time step $\Delta t$. The discretization of the
domain is not needed for the transport step which is solved exactly
for particles but it is necessary to solve the collision part of the
problem which acts locally in space and thus it is necessary to
solve the full problem. In Monte Carlo simulations the distribution
function $f_K$ is discretized by a finite set of particles
\begin{eqnarray}
& & f_K = \frac{m_p}{N}\sum_{i=1}^N  \alpha_i(t) \, \delta(x-X_i(t))
\delta(v-V_i(t)), \label{particle}
\end{eqnarray}
where $X_i(t)$ represents the particle position in the three spatial
directions, $V_i(t)$ the particle velocities in the velocity space,
$m_p$ a constant and $\alpha_i(t)$ the weight to associate to each
particle. During the transport step (\ref{split1}), the particles
move to their next positions according to \be X_i(t+\Delta
t)=X_i(t)+V_i(t)\Delta t\label{transport}\ee where (\ref{particle})
with (\ref{transport}) and $\alpha_i=1$ represents an exact solution
of equation (\ref{split1}). On the other hand, an exact solution of
the modified transport equation \be
\partial_t f_K + v\cdot\nabla_{x}f_K
=\left(\frac{\partial_t h+ v \cdot \nabla_x h }{h}\right)
f_K\label{split1_bis}\ee is a distribution of the form
(\ref{particle}) if $X_i(t)$ and $\alpha_i(t)$ satisfy the system
\begin{eqnarray}
& & \dot X_i = V_i \, , \label{pos_ev} \\
& & \alpha_i(t) =h(X_i,t)  . \label{alpha_ev}
\end{eqnarray}
This can be shown by inserting (\ref{particle}) inside
(\ref{split1_bis}). The previous solution of equation
(\ref{split1_bis}) holds for all choices of $V_i$. The transition
function $h$ is discretized on the mesh of size $\Delta x$ and thus
particles which belongs to the same cell have the same weight. These
relations mean that the weights vanish in the buffer zone
proportionally to the local value of $h$. Consequently, at the
boundary of the buffer zones on the fluid side, the particles are
weightless and can be removed. Simultaneously, the velocity is not
affected by the value of the cut-off function (as for instance in
the decomposition used in \cite{degond2} and \cite{dimarco3}), which
means that there is no accumulation of weightless particles at the
boundaries.

Thus, the Monte Carlo solution of the transport step consists in
moving the particles in space according to equation
(\ref{transport}) as in the classical transport case.
The only difference is that the mass of each particle changes in
time and space according to the transition function value $h$ and it
takes values between $0$ when $h(x,t)=0$ (which means in cells in
which the gas is in thermodynamical equilibrium) and $m_p$ when
$h(x,t)=1$ (which means in the cells in which the gas is far from
equilibrium). The constant $m_p$ is defined at the beginning of the
computation in the following way \be
m_p=\frac{1}{N}\int_{\Omega}\int_{\R^{3}} f_K(t=0) dv dx.\ee Thus,
the transport step is composed of two substeps: transport of
particles and mass rescaling. Finally, to guarantee preservation of
macroscopic quantities we need boundary conditions. This means we
need to inject weightless particles at the boundaries of the buffer
zones on the fluid sides. The gas, being in thermodynamical
equilibrium in these regions, we just sample the number of requested
particles from a Maxwellian distribution whose associated
macroscopic quantities are the ones given by the solution of the
fluid model.
\subsubsection{The moment matching} \label{sec:mm} We
discuss now the matching method between the kinetic equation and the
macroscopic equations. Imagine that the value of the macroscopic
moments $\vr^{n+1}$ obtained from the solution of system
(\ref{rho_K1}) is known at time $n+1$. Now, in the kinetic regions
we want to match these moments with the moments
$\widetilde{\vr}_K^{n+1}$ obtained solving the transport part of the
kinetic equation (\ref{split1_bis}). On the other hand, in buffer
regions, the mass of particles decreases linearly with the
transition function $h$. This means that the matching has to be done
between $h\vr^{n+1}$ and $\widetilde{\vr}_K^{n+1}$, where when $h=1$
we have the matching between the macroscopic equations and the
Boltzmann equation and when $h=0$ we do not perform any matching.
Observe that, the collision step conserves the moments while it
changes the shape of the distribution in velocity space. Thus
$\widetilde{\vr}_K^{n+1}$ represents the moments solution obtained
by solving the kinetic equation (\ref{fK_f1}) at time $n+1$
independently of the type and number of collisions.

We first discuss the matching of momentum and energy which can be
realized through the following linear transformation: let consider
the set of velocities $V_1,\ldots,V_{N_{\mathcal{I}_j}}$ inside the
cell $\mathcal{I}_j$ at time $n+1$, where $N_{\mathcal{I}_j}$ is the
number of particles inside the cell $j$. In the Monte Carlo setting
the first two moments are given by \be
\mu_1=\frac{1}{N_{\mathcal{I}_j}}\sum_{i\in
\mathcal{I}_{j}}^{N_{\mathcal{I}_j}} V_i \qquad
\mu_2=\frac{1}{N_{\mathcal{I}_j}}\sum_{i\in
\mathcal{I}_{j}}^{N_{\mathcal{I}_j}} \frac{1}{2}|V_i|^2, \ee where
$\mu_1$ is a vector representing the mean velocity in the three
spatial directions. Suppose now that better estimates $\sigma_1$ and
$\sigma_2$ of the same moments are available. They are, for
instance, obtained by resolving the moment equations (\ref{rho_K1}).
We apply the transformation described in \cite{Cf} and get
velocities $V_i^*$ given by \be V_i^*=(V_i-\mu_1)/c+\sigma_1\quad
c=\sqrt{\frac{\mu_2-\mu_1^2}{\sigma_2-\sigma_1^2}},\quad
i=1,\ldots,\mathcal{I}_j \label{eq:4}\ee to get
\[
\frac{1}{N_{\mathcal{I}_j}}\sum_{i\in
\mathcal{I}_{j}}^{N_{\mathcal{I}_j}} V_i^*=\sigma_1,\qquad
\frac{1}{N_{\mathcal{I}_j}}\sum_{i\in
\mathcal{I}_{j}}^{N_{\mathcal{I}_j}} \frac{1}{2}|V_i^*|^2=\sigma_2.
\]
Let us now discuss the matching of the zero-th order moment: the
mass density. In this case, we first observe that the above
renormalization is not possible if we want to keep the mass of the
particles constant and uniform inside the cells. We denote by
$\mu_0$ an estimate of the zero-th order moment and by $\sigma_0$
its better evaluation obtained solving the moment equations
(\ref{rho_K1}). Among the possible techniques that can be used to
restore a prescribed density, we choose to replicate or discard
particles inside the cells. Other possibilities are discussed in
\cite{dimarco2}.

In the discarding procedure case we eliminate the following number
of
particles:\be\widetilde{N_p}=\IR\left(\frac{\mu_0-\sigma_0}{m_p}\right),\label{eq:5}\ee
where $\IR$ represents a stochastic rounding. This means that a
mismatch $e$ such that $e<\pm m_p$ is unavoidable when the mass of
the particle is constant and uniform. Observe that, we already
change, at each time step, the mass of the particles according to
the transition function $h$. However, if we change the particles
mass to match the density, this implies that, due to the transport
step, we will end with particles of different weights in the same
cells. This situation will cause difficulties for the collision step
procedure that we want to avoid. In the case of change of mass due
to the transition function, we do not have this problem, because
particles which belong to the same cell have the same mass, $h$
being constant inside each cell.

In the opposite case, in which the mass of the particles inside a
cell is lower than the mass prescribed by the fluid equations
$\mu_0<\sigma_0$, the situation is less simple. Here, since the
distribution function is not known analytically, it is not possible
to sample new particles without introducing correlations between
samples. In this case we replicate \be
\widetilde{N_p}=\IR\left(\frac{\sigma_0-\mu_0}{m_p}\right)\label{eq:6}\ee
particles randomly with uniform probability. We then perform a
stochastic rounding of $(\sigma_0-\mu_0)/{m_p}$ by lower and upper
values. Note that this replication is done allowing repetitions of
the choice of the particle to replicate. After the generation step,
samples are relocated uniformly inside each spatial cell.

\subsubsection{Collision step}
We consider now the solution of the collision step (\ref{split2}).
The effect of this step is to change the shape of the velocity
distribution and to project it towards the equilibrium distribution
leaving unchanged mass momentum and energy. The collision operator
acts locally in space which means that we solve it independently in
each spatial cell where the particles are assumed to have all the
same weight ($h$ is constant in each cell). The initial data of this
step is given by the solution of the transport step after the moment
matching procedure. We assume that the collision term (\ref{split2})
can be written in the form \be
   \partial_t f_K = \frac{h}{\varepsilon} [P(f_K+E[\vr_F],f_K+E[\vr_F])-\mu (f_K+E[\vr_F])],
\label{eq:Kac_like} \ee where $\mu> 0$ is a constant which typically
corresponds to an estimation of the largest spectrum of the loss
part of the collision operator. The operator $P(f,f)$ is a non
negative bilinear operator such that \be \frac1{\mu}\int_{\R^{3}}
P(f,f)(v)m(v)\,dv=\int_{\R^3} f(v)m(v)\,dv. \ee For example for
Maxwellian molecules \be P(f,f)= Q^+(f,f)(v) = \int_{\RR^3}
\int_{S^2} b_0(\cos\theta) f(v^{\prime})f(v_{\ast}^{\prime}) \,
d\omega\,dv_{\ast}, \label{eq:Qpiu} \ee and $\mu=\varrho$. In the
numerical tests presented in this paper, we considered Maxwellian
molecules. In any case, the use of different and more general
kernels does not involve any change in the coupling method. Now,
since the assumption in our model is that the distribution
$E[\vr_F]$ is Maxwellian and constant along the collision step, we
can rewrite (\ref{eq:Kac_like}) in the following way: \be
   \partial_t \left(f_K+E[\vr_F]\right) = \frac{h}{\varepsilon} [P(f_K+E[\vr_F],f_K+E[\vr_F])-\mu (f_K+E[\vr_F])],
\label{eq:Kac_like1} \ee Let now $f_{K,j}^{n*}(v)$ be an
approximation of $f_K(v,n\Delta t,j x_j)$ after transport and
matching. The forward Euler scheme for equation (\ref{eq:Kac_like1})
writes \cite{Nanbu80} \be f_{K,j}^{n+1}+E[\vr_{F,j}^{n+1}] =
\left(1-\frac{h\mu\Delta
t}{\epsilon}\right)\left(f_{K,j}^{n*}+E[\vr_{F,j}^{n*}]\right)
+\frac{h\mu\Delta t}{\epsilon}
\frac{P(f_{K,j}^{n*}+E[\vr_{F,j}^{n*}],f_{K,j}^{n*}+E[\vr_{F,j}^{n*}])}{\mu}.
\label{eq:forward_Euler} \ee Observe that
$E[\vr_{F,j}^{n+1}]=E[\vr_{F,j}^{n*}]$. This means that at the Monte
Carlo level, the above formula can be interpreted in the following
way: in order to obtain the statistical solution of the distribution
function at time step $n+1$, $f^{n+1}_{K,j}$, we need to
\begin{itemize}
  \item sample with probability $(1-{\mu\Delta t}/{\epsilon})$
  from the distribution
$f^{n}_{K,j}$. This means that with probability $(1-{\mu\Delta
t}/{\epsilon})$ the velocity of the particle does not change.
  \item Sample with probability ${h\mu\Delta
t}/{\epsilon}$ from the distribution
$P(f_{K,j}^{n*}+E[\vr_{F,j}^{n*}],f_{K,j}^{n*}+E[\vr_{F,j}^{n*}])/\mu$.
This means that with probability ${h\mu\Delta t}/{\epsilon}$ the
velocity of the sample change as a consequence of a collision with a
particle either coming from $f_{K,j}^{n*}$ or $E[\vr_{F,j}^{n*}]$.
\end{itemize}
To sample from the distribution
$P(f_{K,j}^{n*}+E[\vr_{F,j}^{n*}],f_{K,j}^{n*}+E[\vr_{F,j}^{n*}])/\mu$,
first we need to construct the probability distribution
$f_{K,j}^{n*}+E[\vr_{F,j}^{n*}]$. For constructing this distribution
we just need to sample a number of particle corresponding to the
density $\varrho_{F,j}^{n*}$ from a Maxwellian distribution with
moments $\vr^{n*}_{F,j}$. Then, we apply the operator
$P(\cdot,\cdot)$ to the results. Note that, we request $\Delta t
\leq \epsilon/(h\mu)$ for the probabilistic interpretation to be
valid.

We finally observe that in $\Omega_1$ equation
(\ref{eq:forward_Euler}) reduces to the classical Nanbu scheme for
Maxwellian molecules because $E[\vr_F]$ is identically equal to zero
in this region as well as the transport step reduces to the
classical transport. Conversely in $\Omega_3$ we consider both
models (the macroscopic and the kinetic one) and thus we solve the
coupled collision step described by equation
(\ref{eq:forward_Euler}).

\subsection{The numerical scheme for the moment equations}

Here we discuss the discretization of the moment equations. In the
construction of the numerical scheme we take advantage of the
knowledge of the Euler part of the moment equations \be
\underbrace{\partial_t \vr +\nabla_x \cdot F(\vr)}_{\hbox{Euler
equations}}+\nabla_x \cdot\langle vm(v)g_K \rangle=0.
\label{rho_K2}\ee Thus, the method is based on solving the set of
compressible Euler equations and then considering the discretization
of the kinetic flux $\nabla_x\cdot\langle vm(v)g_K\rangle$. To this
aim, we use a classical discretization in space and time which
means: \be \frac{\vr^{n+1}_{j}-\vr^{n}_{j}}{\Delta t}
+\frac{\psi_{j+1/2}(\vr^n)-\psi_{j-1/2}(\vr^n)}{\Delta
x}+\frac{\Psi_{j+1/2}(<vm(v)g_K^n>)-\Psi_{j-1/2}(<vm(v)g_K^n>)}{\Delta
x}=0.\label{eq:discmac}\ee The numerical flux $\psi$ is an
approximation of the flux $F(\vr^n)$ obtained by the second order
MUSCL extension of a Lax-Friedrichs like scheme. For simplicity, we
indicate in the same way the numerical flux in one or in more
spatial dimensions: \be \label{eq:flux}
\psi_{j+1/2}(\vr^n)=\frac{1}{2}(F(\vr^n_{j})+F(\vr^n_{j+1}))-\frac{1}{2}
A(\vr^n_{j+1}-\vr^n_{j})+\frac{1}{4}(\sigma^{n,+}_j-\sigma^{n,-}_{j+1}),
\ee in the above relation we set \be
\sigma^{n,\pm}_j=\left(F(\vr^n_{j+1})\pm
A\vr^n_{j+1}-F(\vr^n_{j})\mp
A\vr^n_{j}\right)\varphi_{\varepsilon}(\chi^{n,\pm}_j)\ee where
$\varphi_\varepsilon$ is a modified slope limiter, $A$ is the
largest eigenvalue of the Euler system and \be
\chi^{n,\pm}_j=\frac{F(\vr^n_{j})\pm A\vr^n_{j}-F(\vr^n_{j-1})\mp
A\vr^n_{j-1}}{F(\vr^n_{j+1})\pm A\vr^n_{j+1}-F(\vr^n_{j})\mp
A\vr^n_{F,j}}\ee where the above vectors ratios are defined
componentwise. Classical slope limiters, based on the total
variation arguments, determine which regions of the domain can be
solved by a second order method and which regions need a first order
method to avoid the onset of numerical oscillations. Following the
same principle, we define a modified limiter which takes into
account also the departure from the thermodynamical equilibrium. We
introduce this modification because the consequence of the departure
from the thermodynamical equilibrium is the onset of particles in
the domain which are needed to represents the perturbation $g_K$. It
follows that, some statistical error is introduced in the solution
and the use of high order spatial discretization will keep the level
of the oscillations high. Thus, we switch from second to first order
in this case and we use the following function to perform this
switch: \be
\varphi_\varepsilon(\chi_j)=\varphi_{L}(\chi_j)(1-h(j))\ee where
$\phi_L(\chi)$ is the van Leer limiter \be
\phi_L(\chi)=\frac{|\chi|+\chi}{1+\chi}.\ee The MUSCL second order
scheme is then used in fluid regions if $\phi_L=1$, while the first
order scheme is used otherwise.

Concerning the discretization of the non equilibrium term
$\nabla_x\cdot <vm(v)g_K>$, the same space first order discrete
derivative is used as for the hydrodynamic flux $F(\vr)$:\be
\Psi_{j+1/2}(<vm(v)g_K^n>)=\frac{1}{2}\left(F(\langle
vm(v)g^n_{K,j}\rangle)+F(\langle vm(v)g^n_{K,j+1}\rangle)\right).\ee
The non equilibrium term $\langle vm(v)g_K\rangle=\langle
vm(v)(f_K-E[\vr_K])\rangle$ is computed by taking the difference
between the moments of the particle solution and those of the
Maxwellian equilibrium. Thus, the flux $\psi_{j+1/2}(\vr^n)$ can be
either first or second order while $\Psi_{j+1/2}(<vm(v)g_K^n>)$ is
always first order.

Observe that we do not need boundary condition for the moment
equations at the interfaces because they are solved in all the
domain. The only difference, between the fluid and the kinetic
regions, is that the perturbation term disappears in the equilibrium
zones. Then, thanks to the smooth transition between the two models
and the variance reduction technique employed to reduce the variance
of the DSMC method, fluctuations do not propagate to the fluid zones
as shown in the numerical test section. Finally, the CFL condition
is chosen such that the time step is always the minimum between the
relaxation parameter $\varepsilon$, the ratio between the mesh size
and the largest particle velocity ($v_{\max}$) and the ratio between
the mesh size and the largest eigenvalue of the fluid equations
($A_{\max}$): \be \Delta t=\min\left(\frac{\Delta
x}{v_{\max}},\frac{\Delta x}{A_{\max}},\varepsilon\right)\ee


\section{The breakdown criterion} \label{sec moving}
For defining the interface location, we will try to combine the
local Knudsen number with the effects that the numerical scheme has
on the solution. We observe that, the smaller uncertainly in the
macroscopic quantity values which has been obtained with the moment
guided method, permits to better define the breakdown of the fluid
model and consequently to optimize the coupling. The Knudsen number
is defined as the ratio of the mean free path of the particles
$\lambda_{path}$ to a reference length $L$: \be
\varepsilon'=\lambda_{path}/L, \label{eq:knudsen}\ee where the mean
free path is defined by
$$
\lambda_{path}=\frac{kT}{\sqrt{2}\pi p \sigma_{c}^{2}} ,
$$
with $k$ the Boltzmann constant equal to $1.380062\times 10^{-23}
JK^{-1}$, $p$ the pressure and $\pi \sigma_{c}^2$ the collision
cross section of the molecules. Now, in order to take into account
the elementary fact that, even in extremely rarefied situations, the
flow can be in thermodynamic equilibrium, as in Bird \cite{bird},
the reference length is defined as \be
L=\min\left(\frac{\varrho}{\partial\varrho/\partial x},\frac{\varrho
u}{\partial\varrho u/\partial x},\frac{\varrho e}{\partial \varrho
e/\partial x}\right).\ee  According to~\cite{Levermore}
and~\cite{Kobolov}, the fluid model is accurate enough if the local
Knudsen number is lower than the threshold value $0.05$. It is
argued that, in this way, the error between a macroscopic and a
microscopic model is less than 5\% \cite{Wang}. This parameter has
been extensively used in many works and is now considered in the
rarefied gas dynamic community as an acceptable indicator.

In the numerical tests section we distinguish between the physical
Knudsen number $\varepsilon'$ defined in (\ref{eq:knudsen}) and the
relaxation parameter $\varepsilon$ which appears in the model and in
the numerical schemes defined in the previous sections. In practice
in our tests $\varepsilon$ is fixed at the beginning of the
simulation and it does not change as a function of the macroscopic
quantities. In other words, $\varepsilon$ is just a rescaling
parameter of the Boltzmann equation, and all the quantities will be
dimensionless. This choice permits to control the relaxation towards
the equilibrium and thus to focus on the behavior the numerical
scheme.

The criterion we choose is based on measuring the ratio between the
relaxation time and the time step. We recall that, the time step is
chosen as the minimum between the relaxation time and the transport
time (both for the fluid and for the kinetic equations). Thus, if
the time step is the same as the relaxation time, it means that all
the particles collide at each time step. In this case we make the
hypothesis that we are in a fluid regime. On the other hand, if at
each time step, only few particles collide, we can be both in
equilibrium or not depending on the values of the derivatives of the
macroscopic quantities. Thus, the breakdown criterion is defined as
\be \beta=\left(1-\frac{\mu\Delta t}{\varepsilon}\right)\frac{\Delta
x}{L}.\label{break}\ee The first part of the above formula
corresponds to the first weight of the forward Euler scheme for
solving the Boltzmann collision term introduced before. The second
part is just the dimensionless reference length. Now, the smaller
the relaxation parameter is, the smaller the influence of the
derivative of macroscopic quantities on the evaluation of the
equilibrium of the gas is. On the other hand, the larger the
difference between the macroscopic quantities in adjacent cells is,
less likely the equilibrium is guarantee. The parameter $\beta$ is
substantially equivalent to the Knudsen number but it takes into
account the mesh discretization and the choice of the time step.

\begin{remark}~
In a previous work \cite{dimarco4}, we proposed a criterion to
locate the interface between the fluid and the kinetic regions based
on the analysis of the micro-macro decomposition. Unfortunately,
this criterion as well as other similar criteria which measure some
non equilibrium quantities (see for instance \cite{Levermore}) needs
a precise evaluation of the macroscopic quantities. However, when
DSMC methods are used to compute the solution especially in unsteady
situations it is very difficult to evaluate precisely these
quantities. This is moreover true when relatively few particles are
used to compute the kinetic solution. On the other hand, when the
number of particles employed to compute the solution is large enough
these indicators give good estimations of the departure from
equilibrium. It is important to remark that, the indicator $\beta$
needs less precise evaluations of the macroscopic quantities and
thus, its use is more suitable for cases in which only few particles
are used to compute the solution.
\end{remark}~

\section{Numerical tests}
\label{sec_tests}
\subsection{General setting}
\label{subsec_gen_set}

Here we present several numerical results to highlight the
performances of the method. All the tests we present represent
unsteady physical problems. By using these kinds of tests we want to
emphasize one of the main characteristics of our method: we are able
to catch departures from the thermodynamical equilibrium which move
in time. Moreover, we are able to determine the necessity of kinetic
models even using Monte Carlo methods and with relatively few
particles per cell. Of course, in the case of steady test cases, the
possibility to average quantities permits to realize simulations
using only few particles per cell (typically 50 or less particles
per cell). On the other hand, when the solutions are unsteady even
with 300-400 particles in average per cell, many details of the
solution are lost. In practice the number of particles needed
depends on the problem and on the requested resolution needed.
Nevertheless, we will show that even with a relatively low number of
particles, we are able both to catch departure from equilibrium and
to describe the solutions enough accurately. Finally, if more
precise solutions are requested, we point out that increasing the
number of particles leads to more precise interface locations and
less noisy evaluations of macroscopic quantities than classical DSMC
methods and classical coupling methods. On this subject, we recall
that the convergence of the moment guided Monte Carlo methods is
faster than the convergence of DSMC for increasing number of
particles \cite{dimarco1}.

In the first test, we consider a space dependent relaxation
frequency $\varepsilon^{-1}$. This test mimics the case in which
different regions with different collision rates have to be
considered in the same simulation. As we will see, when these
problems are present, constructing a domain decomposition method is
relatively simple, the interface location being well defined.

Next, as in \cite{dimarco4}, we consider an unsteady shock test
problem. In this case the relaxation time is taken constant in the
domain. We remark that in this situation, a standard static domain
decomposition fails. In fact the shock moves in time. When the
relaxation time is large, i.e. in rarefied regimes, it is necessary
to use a kinetic solver in the full domain. On the other hand, with
our algorithm, we are able to catch the shock region and to describe
only this part of the domain with the kinetic model.

Finally, as in \cite{dimarco4}, we use our scheme to compute the
solution of the Sod test. Again, the relaxation frequency is taken
constant in the domain. Here, the structure of the solution is more
complicated but, as described below, the method is still able to
deal with such a situation and to realize an efficient coupling.

In all tests the time step is given by the minimum of the maximum
time steps allowed by the kinetic and fluid schemes. The speed-up we
obtain is only due to the reduction of the sizes of the kinetic and
buffer regions inside the domain. This reduction is achieved through
a correct prediction of the evolution of the transition function. It
is measured through the difference between the number of particles
used for the full DSMC simulation and the number of particles used
in the DSMC/Fluid coupling considering the same number of initial
particles. We point out that all the procedure is automatic and
determined by the step by step algorithm presented in the previous
section. For all tests, we set $\beta_{thr}=2.5 \ 10^{-2}$ while the
buffer regions thickness is fixed constant for every test.


\subsection{Two relaxation frequencies test}
\label{subsec_relaxation}

In this paragraph we present a two regimes fluid-kinetic test. The
domain ranges from $x=0$ to $x=1$ and the number of cells is $200$.
For a full resolution with the DSMC method the number of particles
employed is $80000$ which corresponds to an average of $400$ per
cell. In the left half of the domain $\varepsilon$ is taken equal to
$10^{-4}$, while in the right half we take $\varepsilon=10^{-2}$.
Neumann boundary condition are chosen both for the kinetic and fluid
models. At the beginning of the simulation the gas is considered in
thermodynamical equilibrium, in other word we set $h=0$ everywhere.

\begin{figure}
\begin{center}
\includegraphics[scale=0.36]{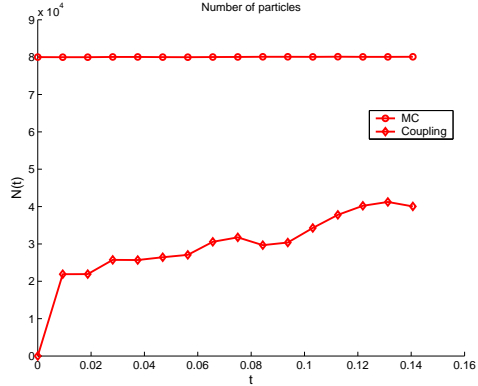}
\caption{Two relaxation frequencies test: number of particles in
time for the Monte Carlo scheme and the DSMC-Fluid Coupling.}
\label{DT03}
\end{center}
\end{figure}

\begin{figure}
\begin{center}
\includegraphics[scale=0.39]{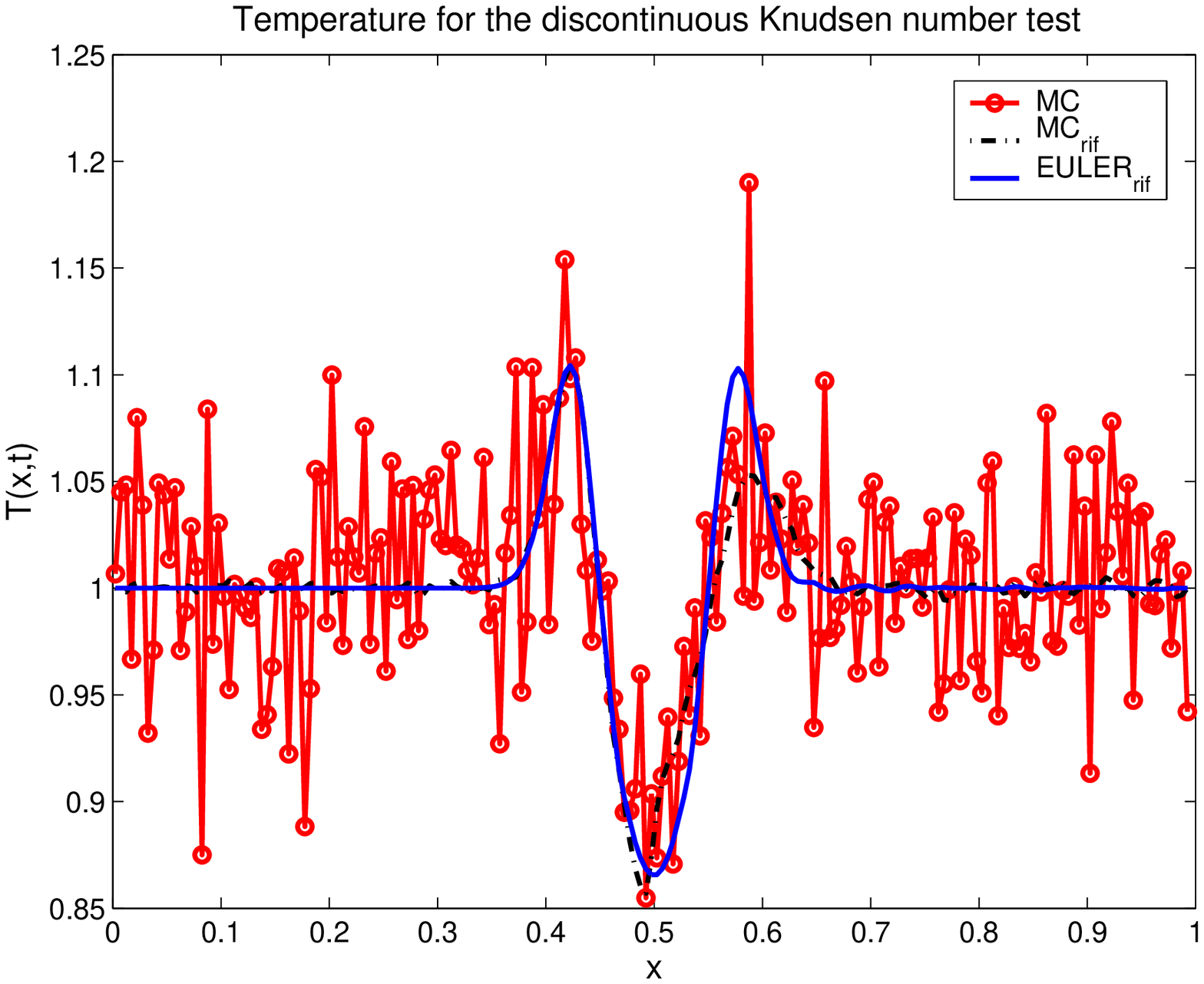}
\includegraphics[scale=0.39]{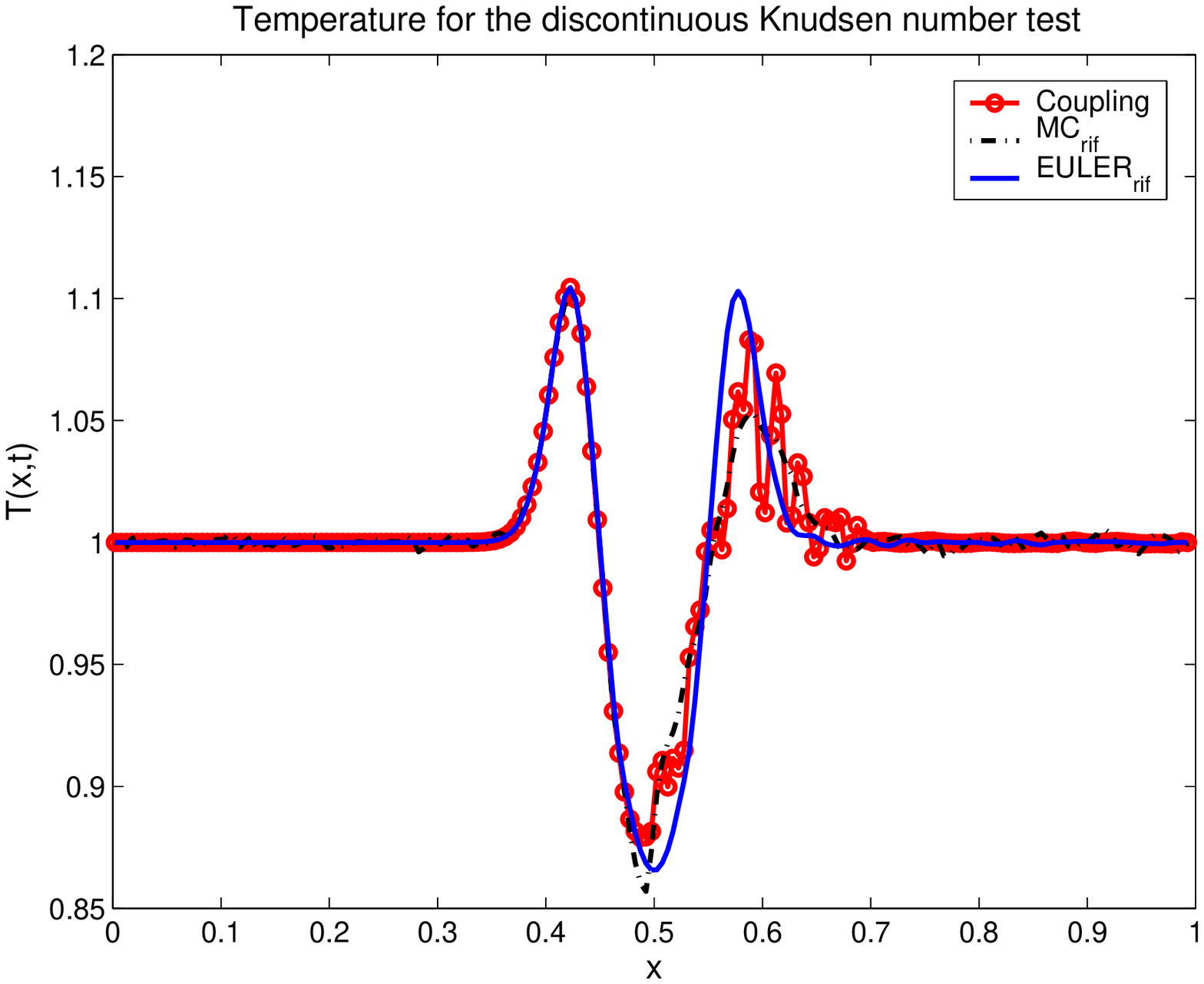}\\
\includegraphics[scale=0.39]{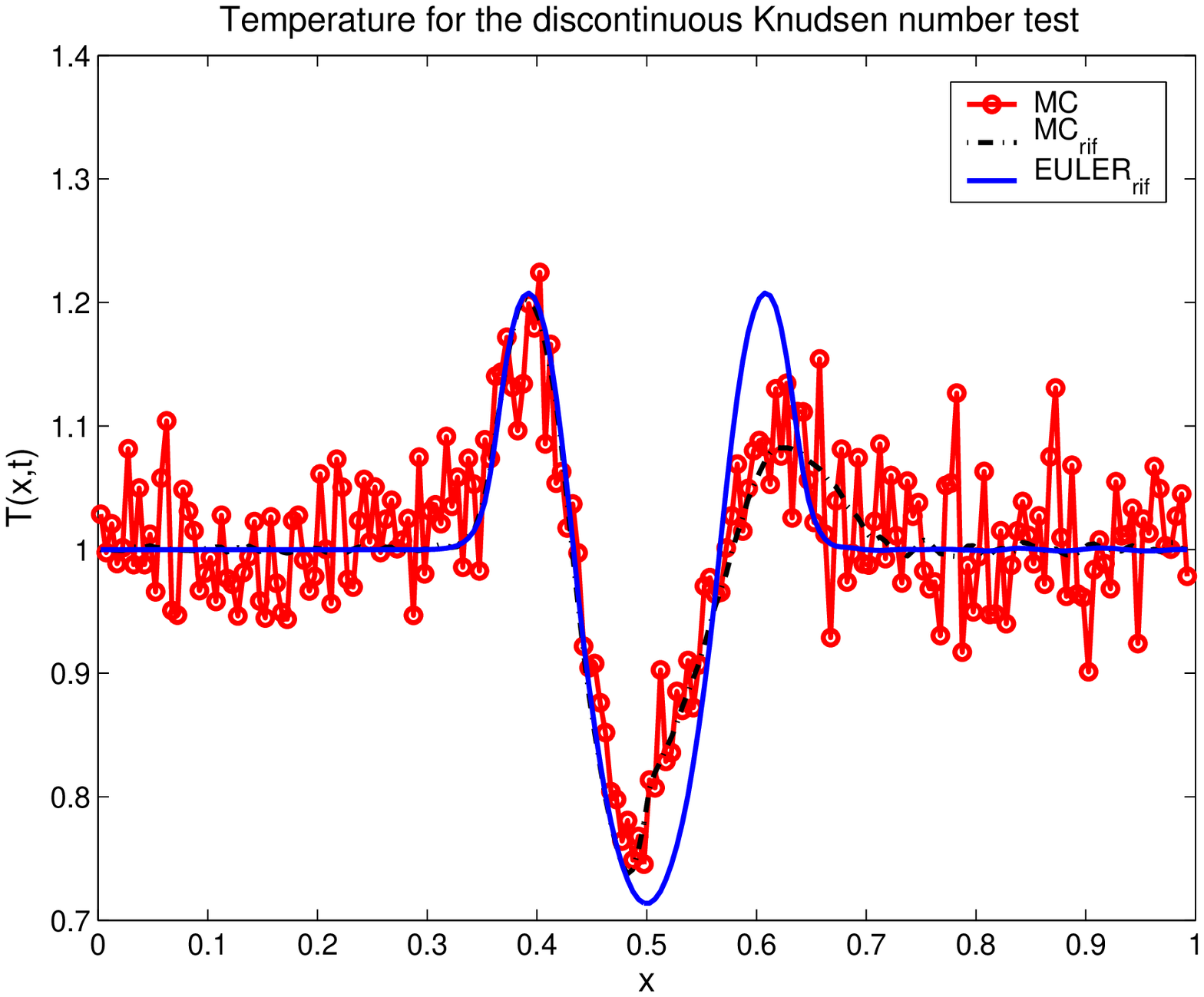}
\includegraphics[scale=0.39]{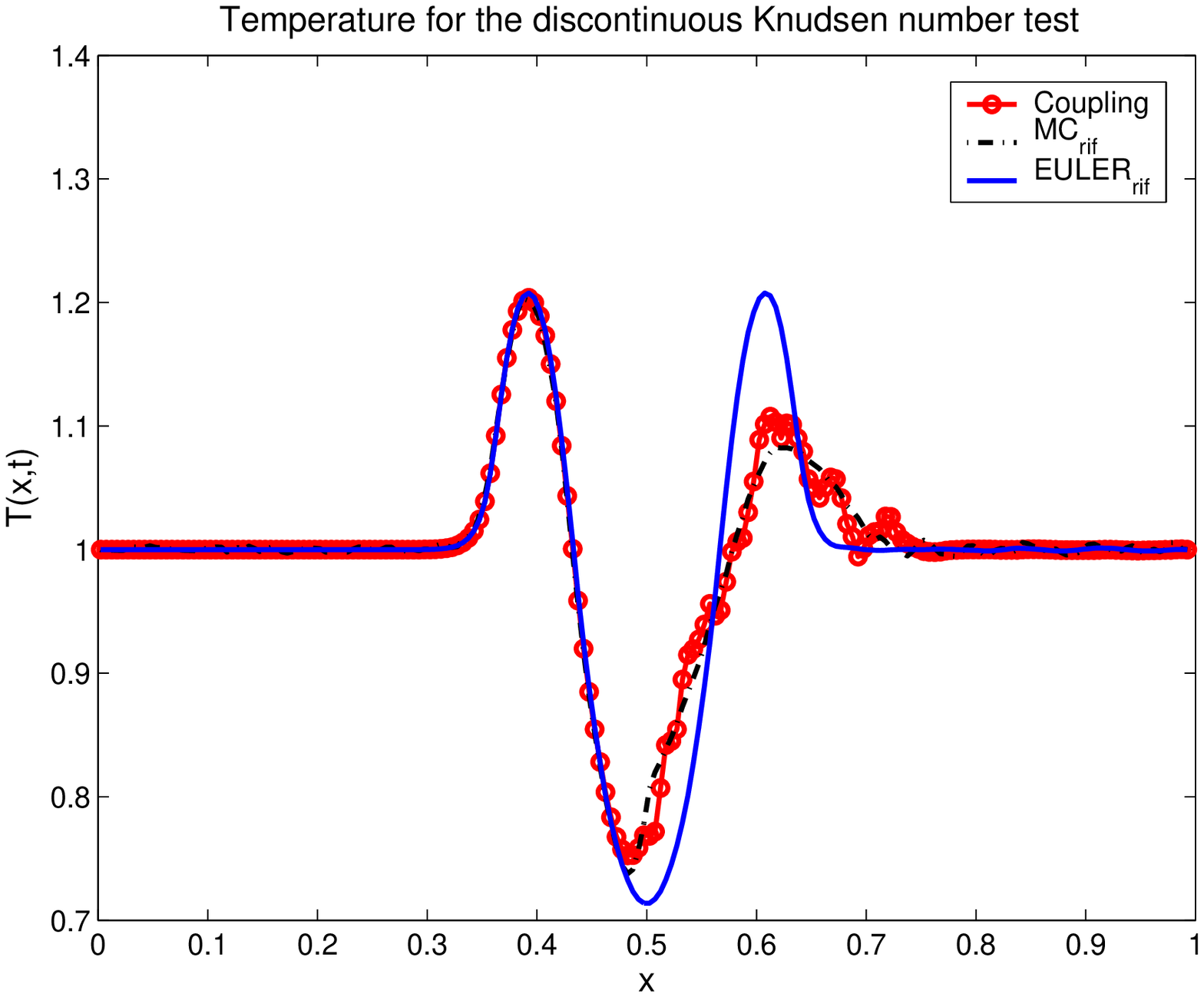}\\
\includegraphics[scale=0.39]{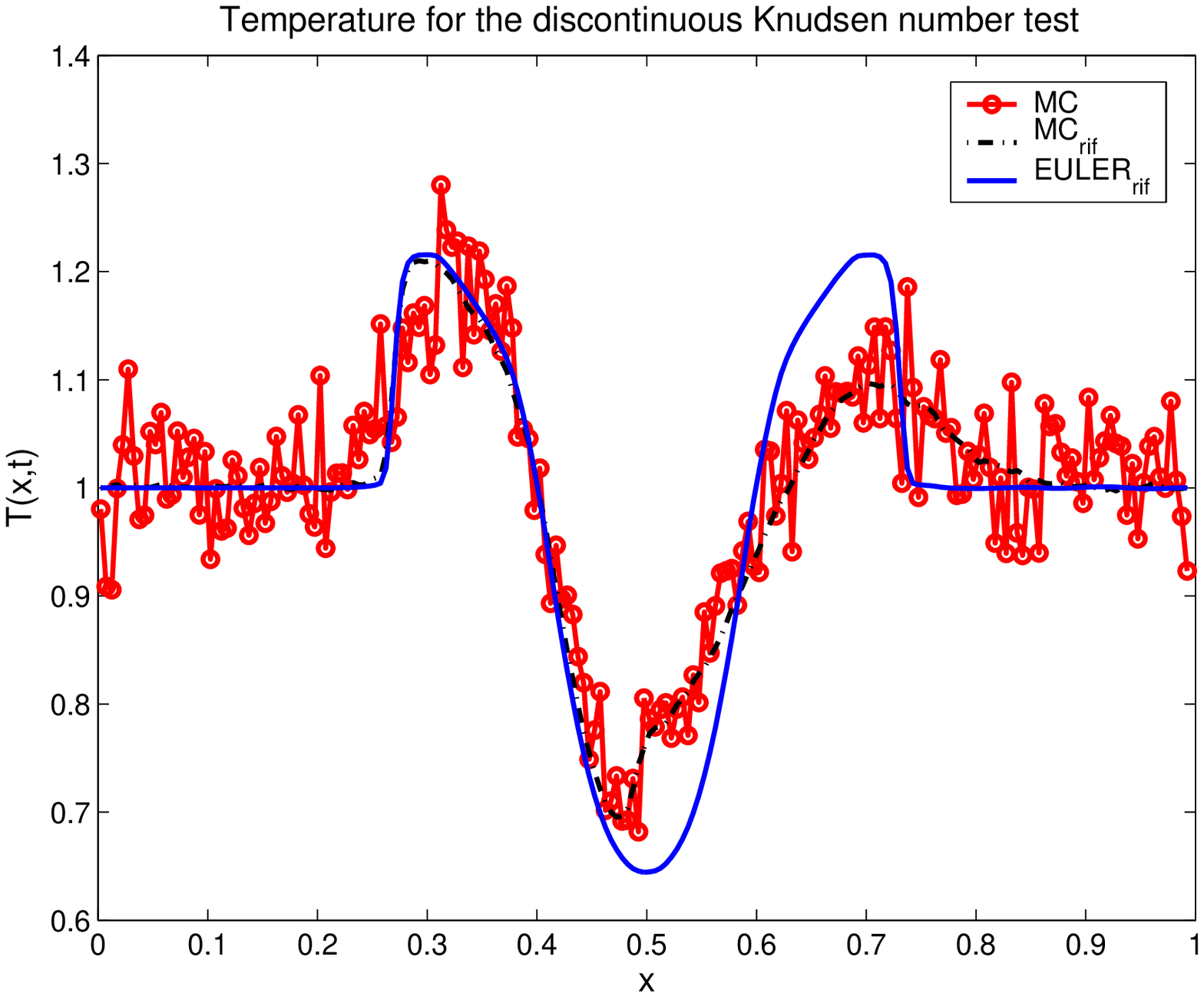}
\includegraphics[scale=0.39]{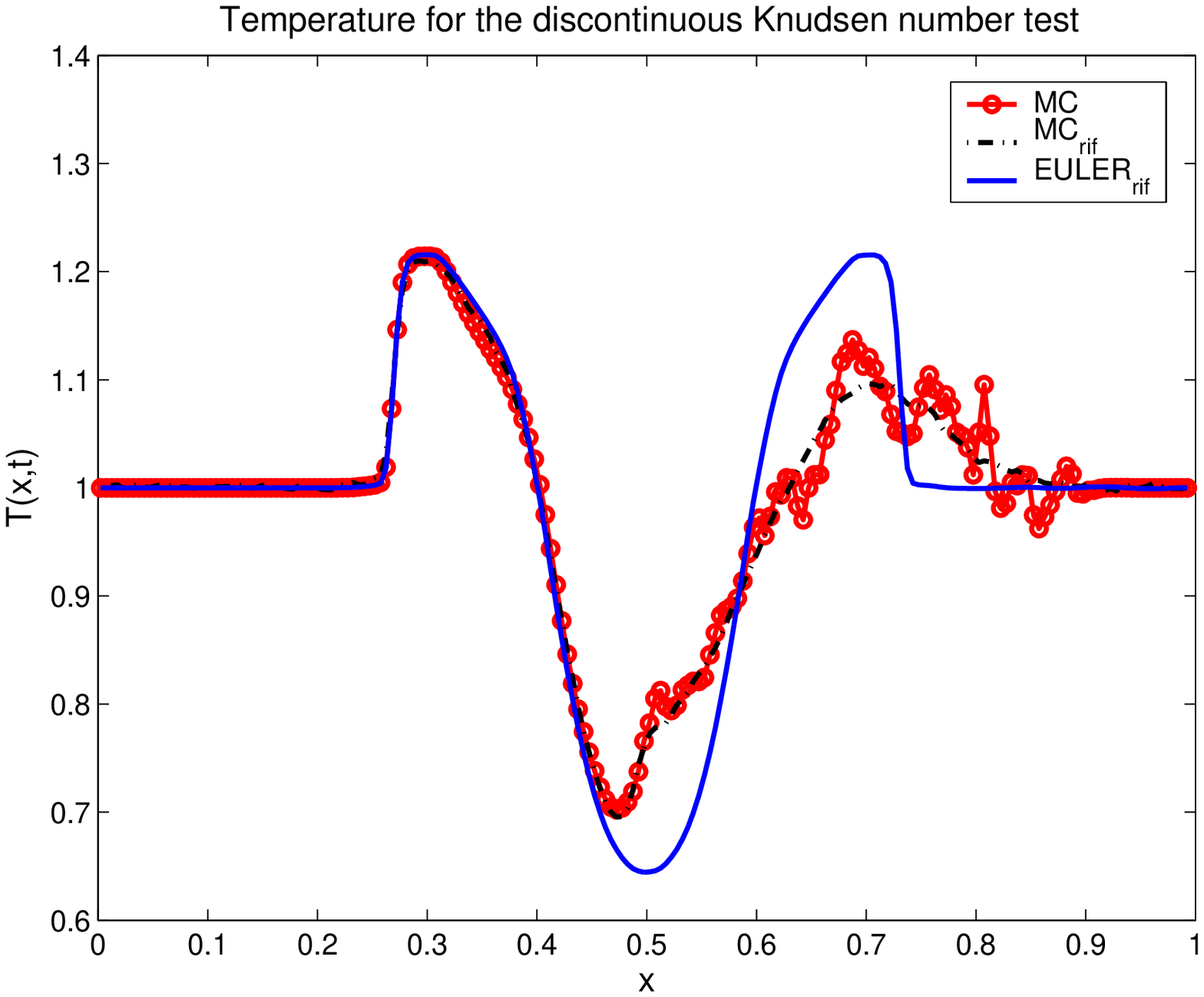}
\caption{Two relaxation frequencies test: Solution at $t=0.05$
(top), $t=0.10$ (middle) and $t=0.15$ (bottom) for the temperature.
MC method (left), Coupling DSMC-Fluid method (right). Knudsen number
$\varepsilon=10^{-4}$ in the left half of the domain and Knudsen
number $\varepsilon=10^{-2}$ in the right half of the domain.
Reference solution (dotted line), Euler solution (continuous line),
DSMC-Fluid or DSMC (circles plus continuous line).} \label{DT02}
\end{center}
\end{figure}

The initial data are: \be \nonumber
\varrho=1+0.1\sqrt{\frac{(x-0.5)^{2}}{0.02}} \ \ u=0 \ \ T=1.\ee In
figure \ref{DT03} the number of particles in time is reported for
the DSMC method and the DSMC/Fluid coupling method. The interface
position is given by the breakdown parameter $\beta$ and it appears
fixed in time and space. This fixed location reproduces the boundary
between the low and the large values of the relaxation frequency
$\varepsilon^{-1}$. A smooth transition region (ten cells thick)
between the two models is also fixed in time and space. Figure
\ref{DT02} reports the profile of the temperature for increasing
times from top to bottom, on the left side the DSMC and on the right
side the DSMC/Fluid coupling method. In each figure we plot the
solution of the Euler equation and a reference solution obtained
with the DSMC scheme where the number of particles is such that the
statistical error is very small.

We observe that, in the case in which the equilibrium and non
equilibrium zones are well identified by the number of collisions
that the particles encounter each time step, the coupling works very
well. For the coupling method case we observe on the left of the
domain an high order finite volume scheme which is able to catch
well discontinuities. On the other hand, on the right side, thanks
the moment matching method, the DSMC method exhibits reduced
fluctuations compared to the classical DSMC method. Moreover, thanks
to the smooth transition, the fluctuations do not propagate in the
fluid regions. We finally recall that, the reduction of the
statistical error due to the moment matching strongly depends on the
frequency $\varepsilon^{-1}$. In particular for large relaxation
frequencies the reduction of the fluctuations are much more obvious.
We will see, in the test cases presented in the next paragraph, that
this reduction can be larger. This is due to the fact that we will
describe with the kinetic model also regions with large relaxation.
Finally, we recall that the reduction of the computational cost is
proportional to the total number of particles used in the
simulations. In the next test cases we will consider problems with
high rarefaction states in all the domain. In these cases catching
the departure from the thermodynamical equilibrium will be less
easy.

\subsection{Unsteady shock tests}
\label{subsec_unst_shock}

\begin{figure}
\begin{center}
\includegraphics[scale=0.36]{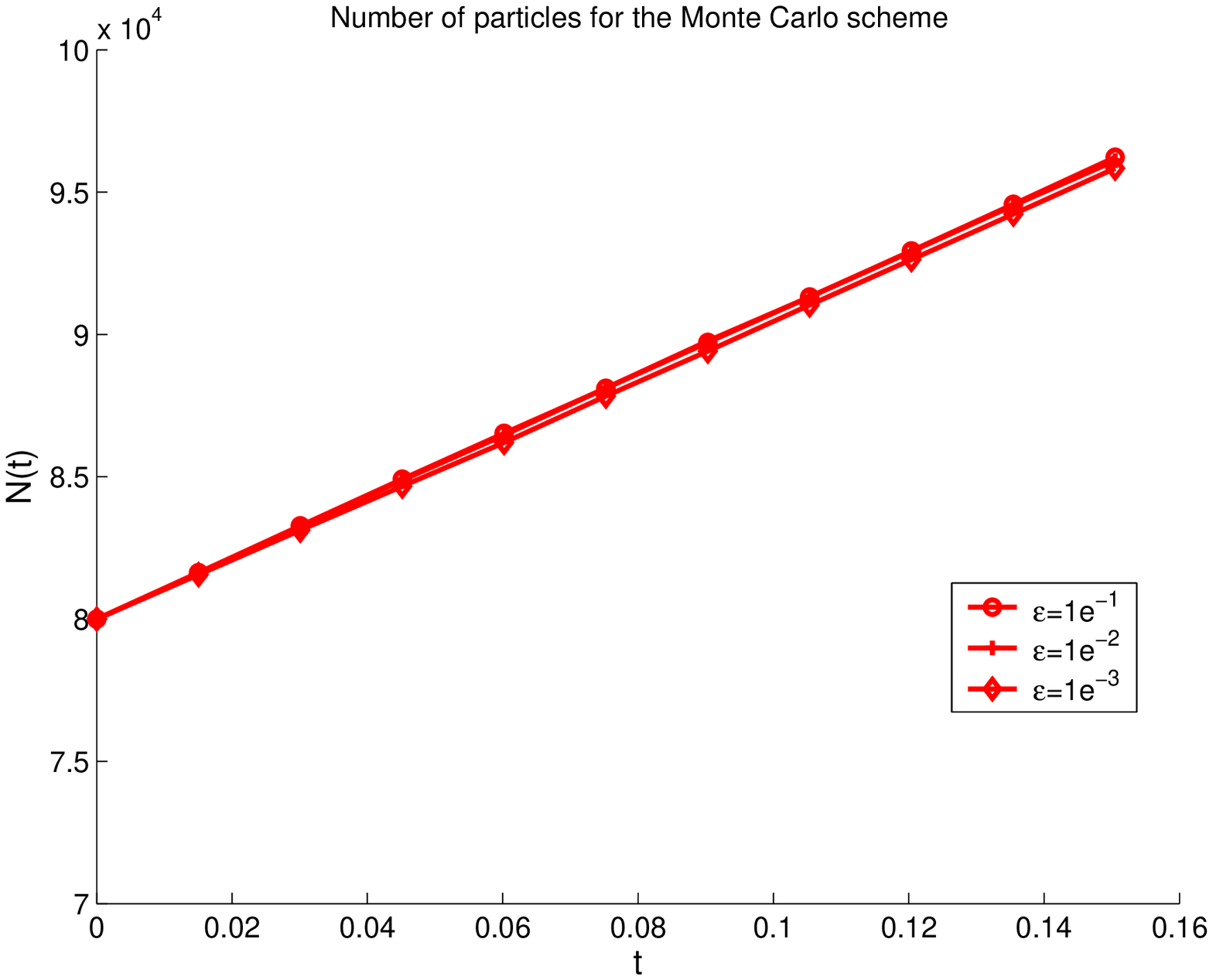}
\includegraphics[scale=0.36]{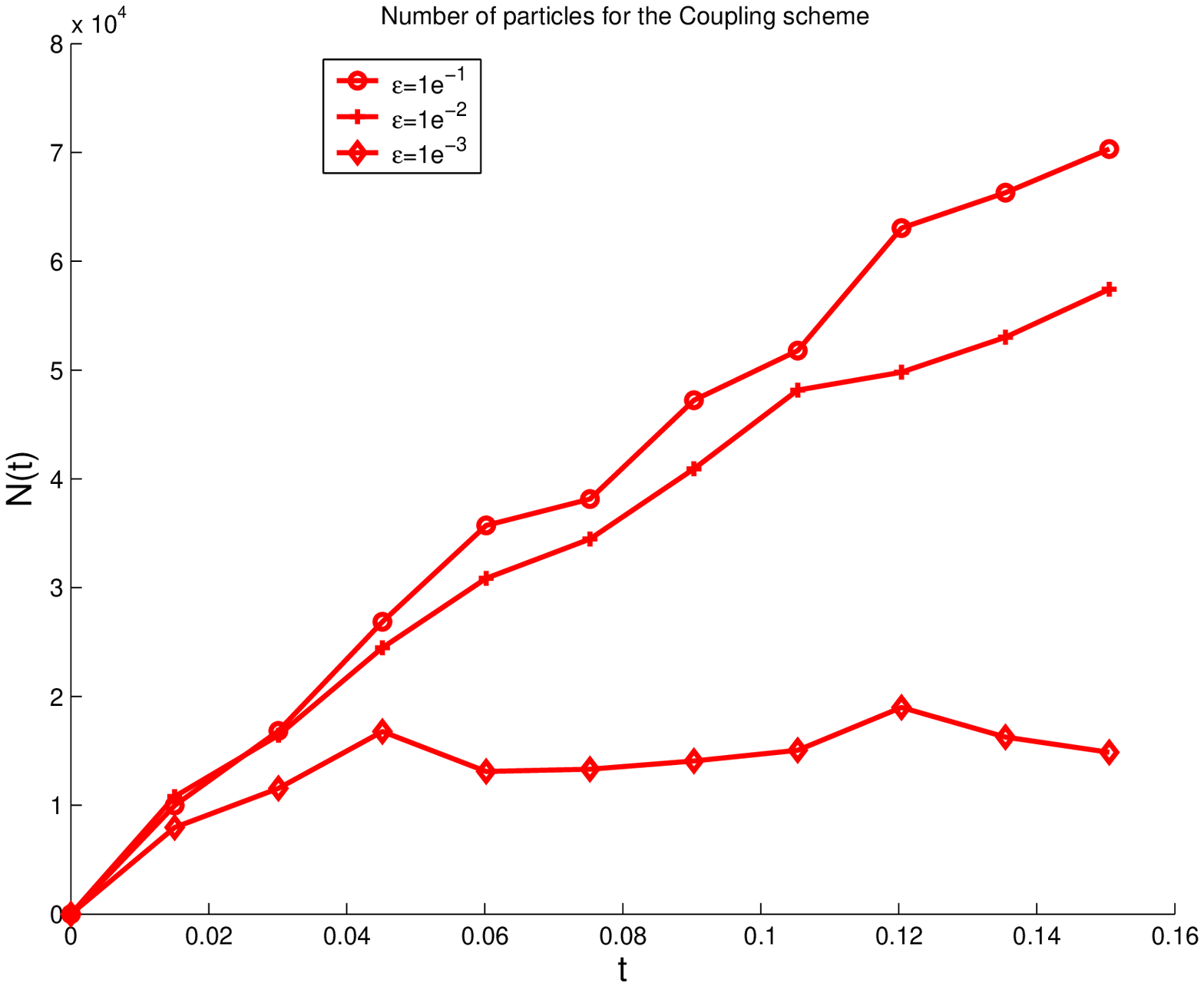}
\caption{Unsteady Shock Test: Number of particle in time for the
Monte Carlo scheme (left) and the DSMC-Fluid Coupling (right) for
different values of the Knudsen Number (from $\varepsilon=10^{-1}$
to $\varepsilon=10^{-3}$).} \label{ST24}
\end{center}
\end{figure}

\begin{figure}
\begin{center}
\includegraphics[scale=0.27]{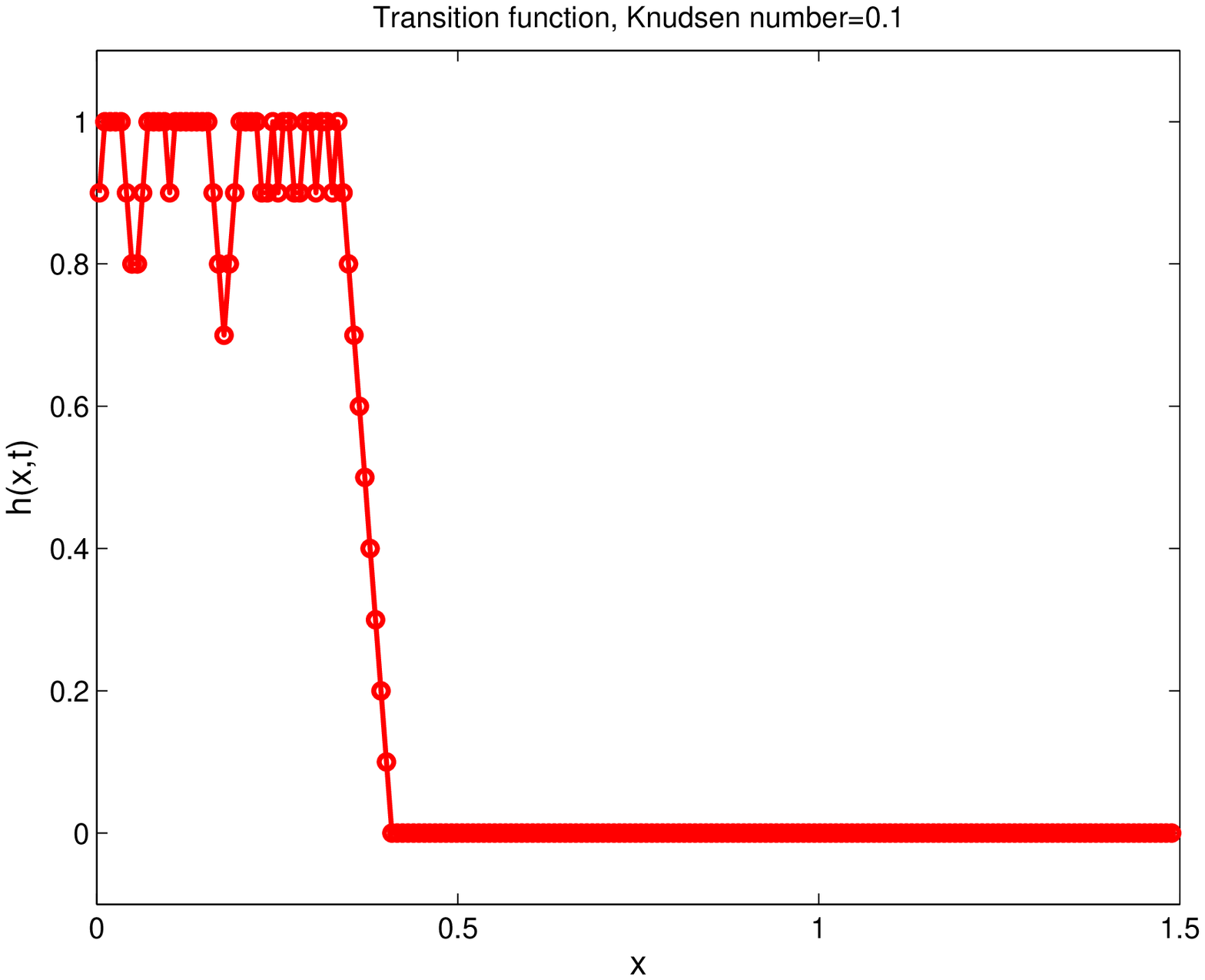}
\includegraphics[scale=0.27]{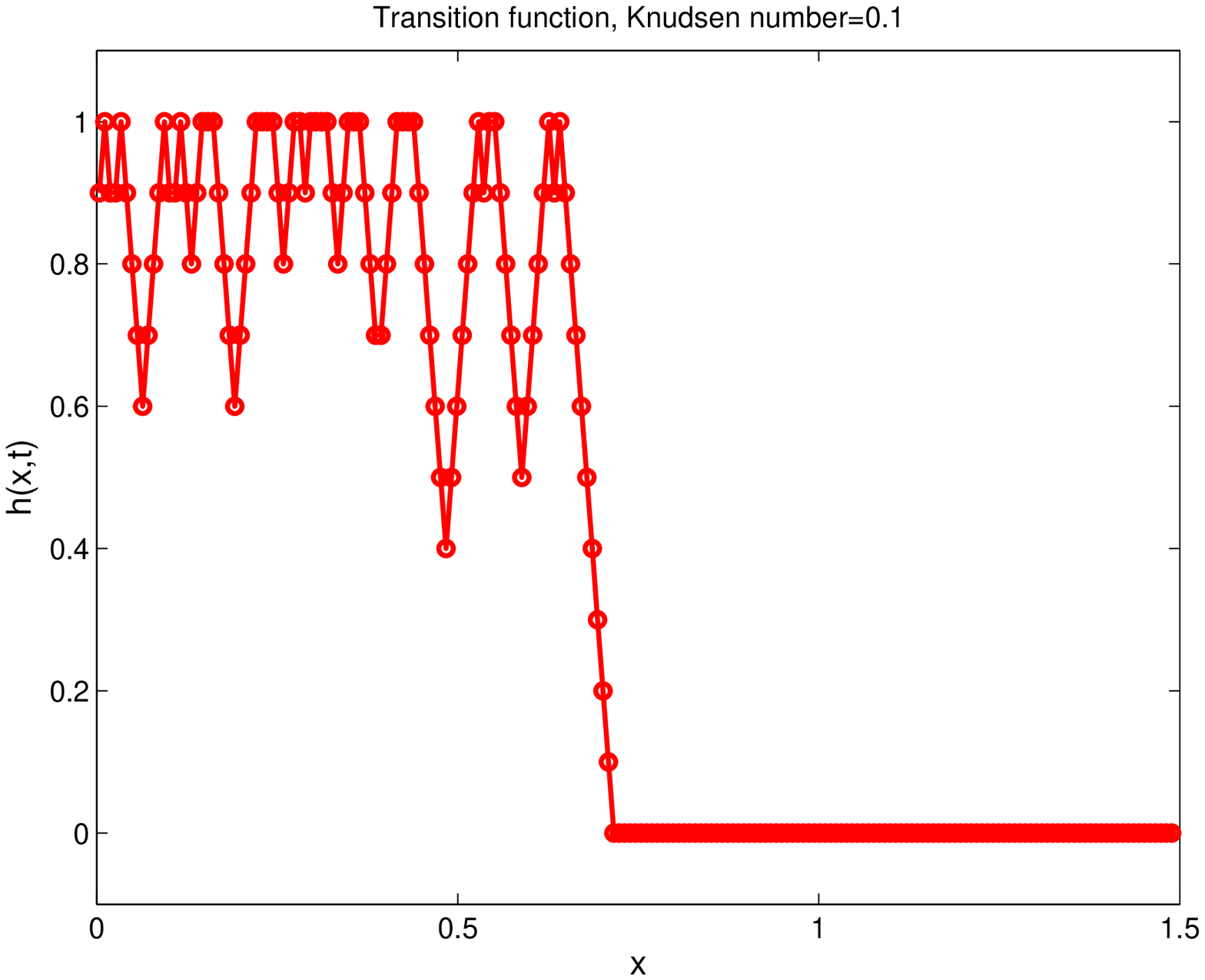}
\includegraphics[scale=0.27]{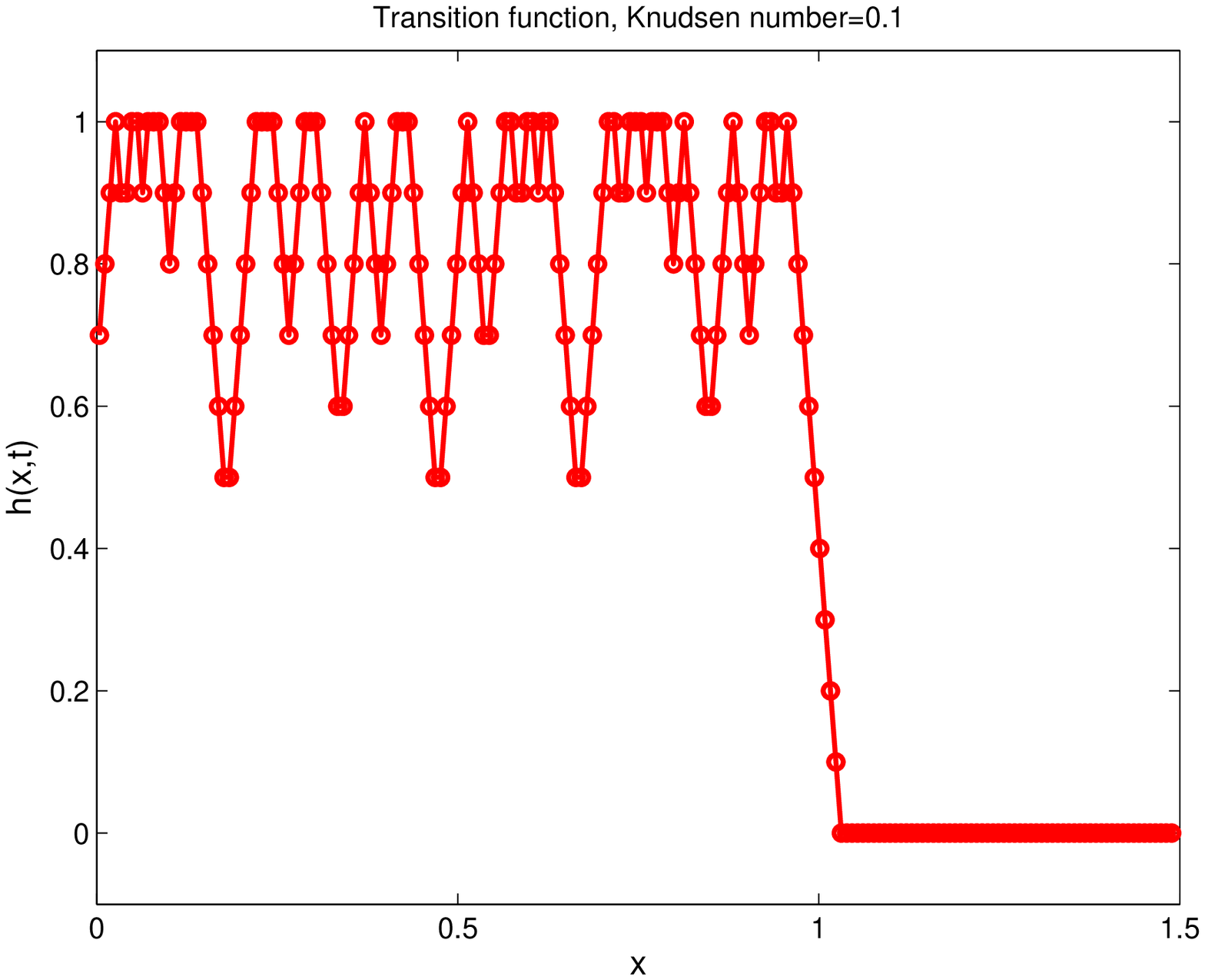}
\caption{Unsteady Shock Test: Transition function at $t=0.05$
(left), $t=0.10$ (right) and $t=0.15$ (bottom). Knudsen number
$\varepsilon=10^{-1}$.} \label{ST03}
\end{center}
\end{figure}

We consider an unsteady shock that propagates from left to right in
the computational domain $x=0$, $x=1.5$ discretized with $200$ cells
in space. The number of particles is such that $400$ particles
correspond to $\varrho=1$. The shock is produced pushing a gas
against a wall which is located on the left boundary. In our test
the particles are specularly reflected and the wall adopts the
temperature of the gas instantaneously.
The gas is supposed in thermodynamic equilibrium at the beginning of
the simulation. The computation is stopped at the final time
$t=0.15$. The transition function $h$ is initialized as $h=0$ (fluid
region) everywhere.

\begin{figure}
\begin{center}
\includegraphics[scale=0.27]{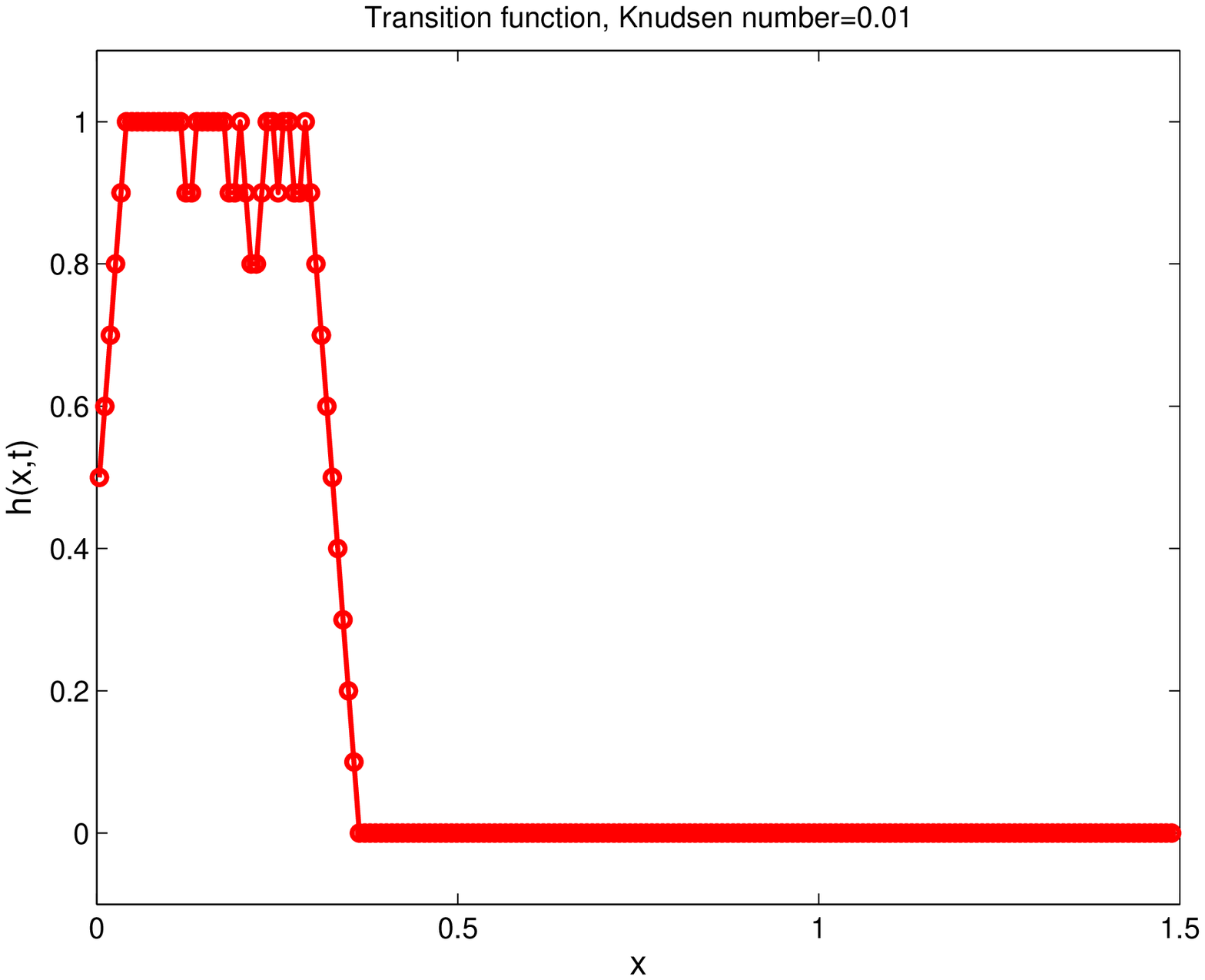}
\includegraphics[scale=0.27]{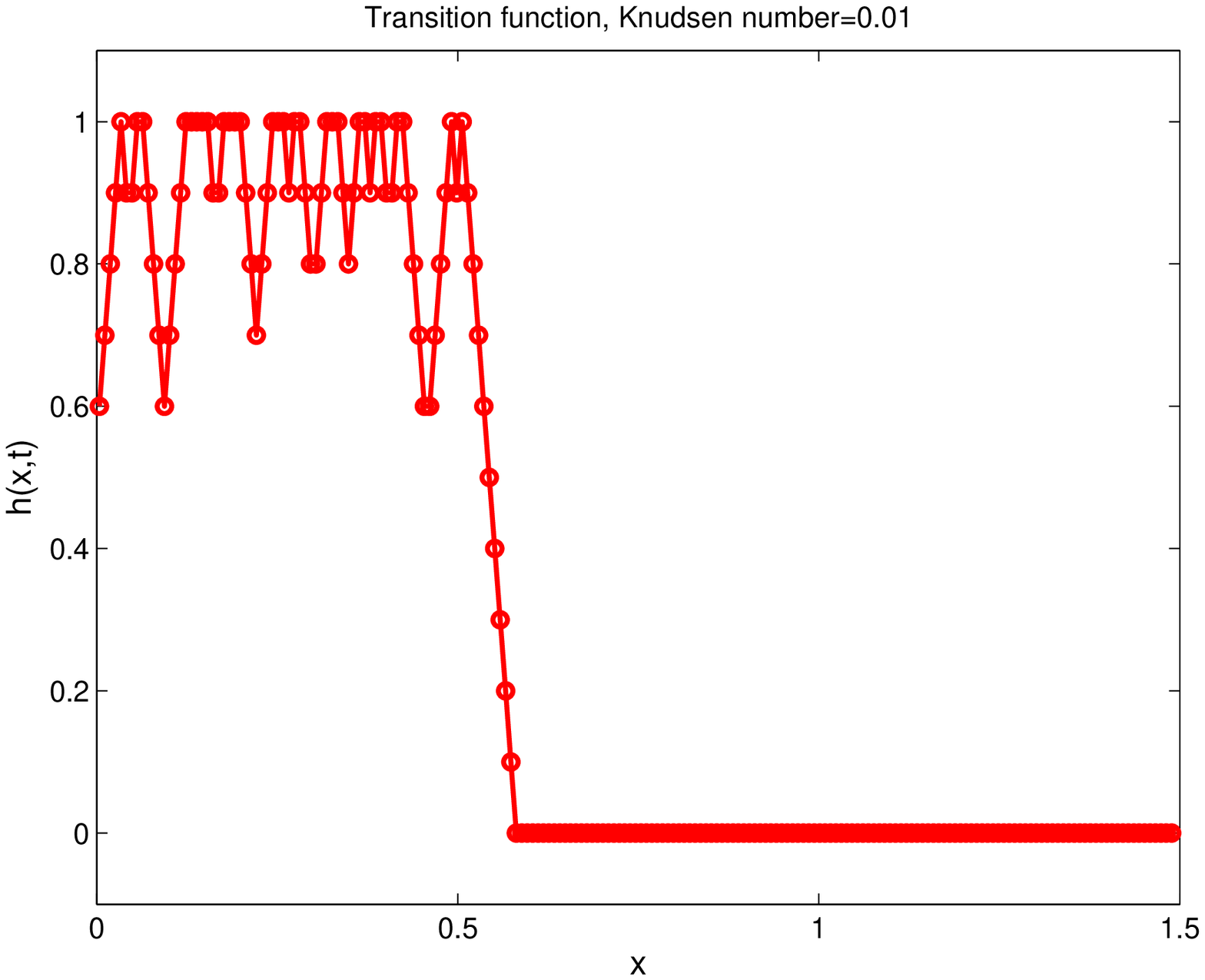}
\includegraphics[scale=0.27]{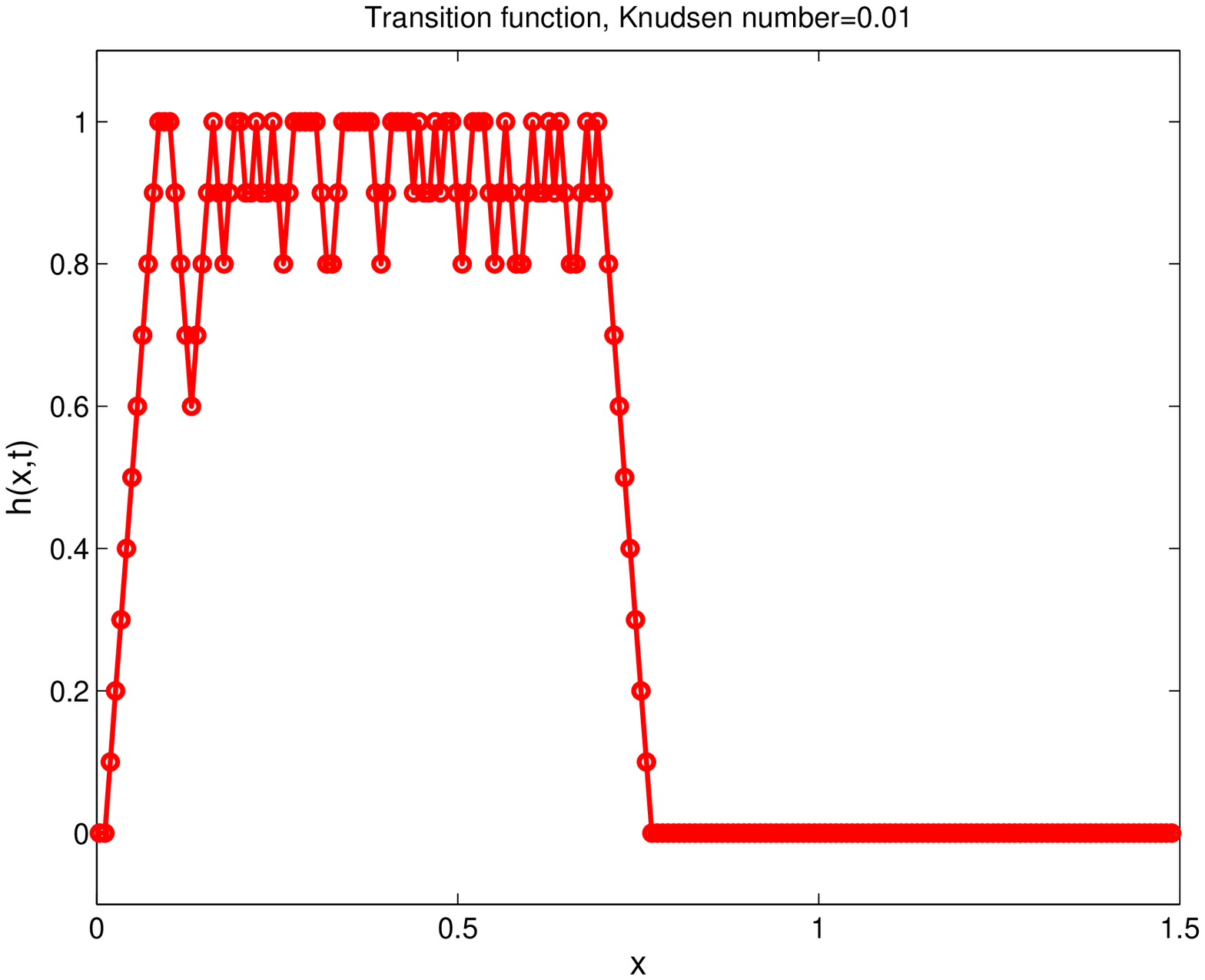}
\caption{Unsteady Shock Test: Transition function at $t=0.05$
(left), $t=0.10$ (right) and $t=0.15$ (bottom). Knudsen number
$\varepsilon=10^{-2}$.} \label{ST13}
\end{center}
\end{figure}

\begin{figure}
\begin{center}
\includegraphics[scale=0.27]{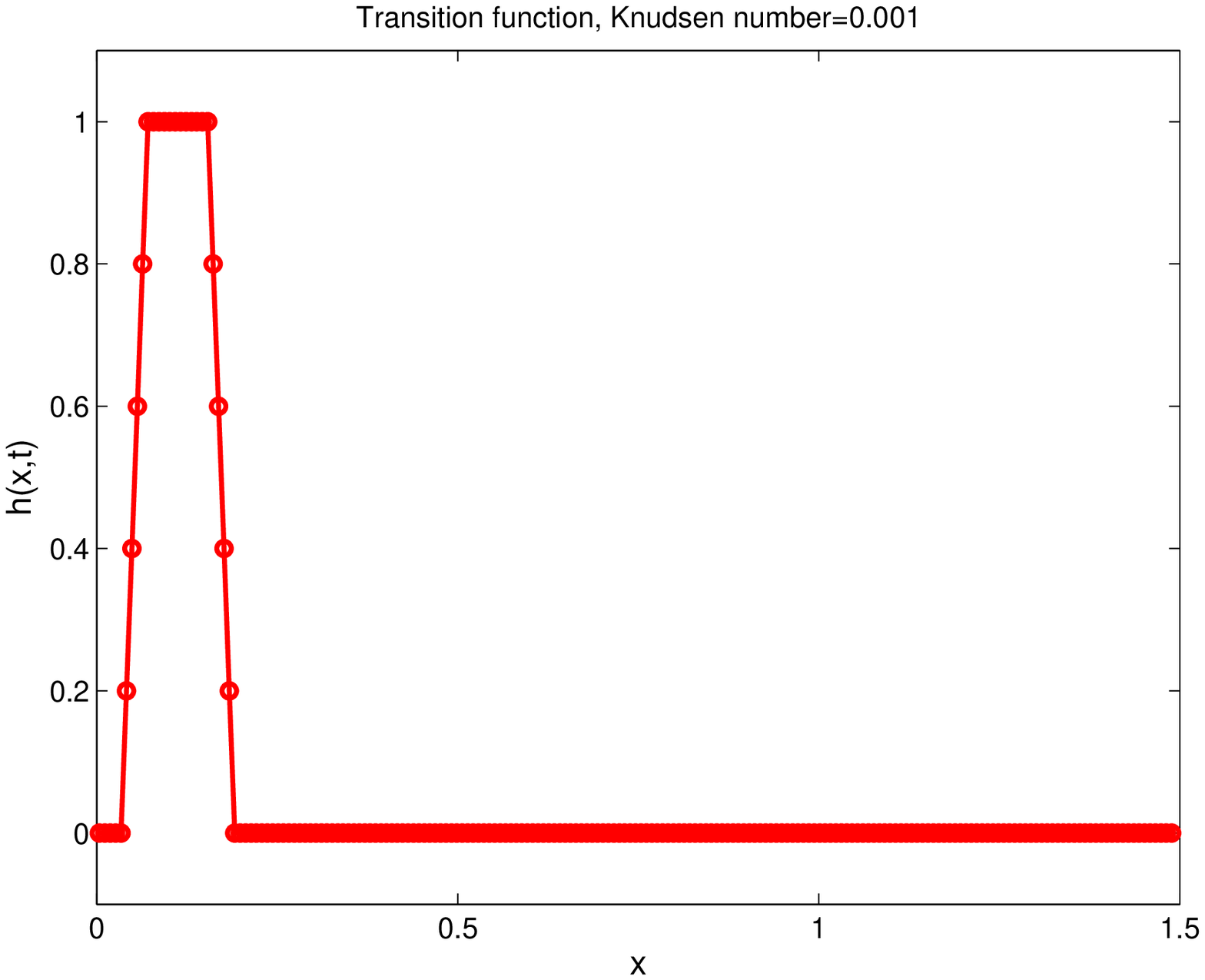}
\includegraphics[scale=0.27]{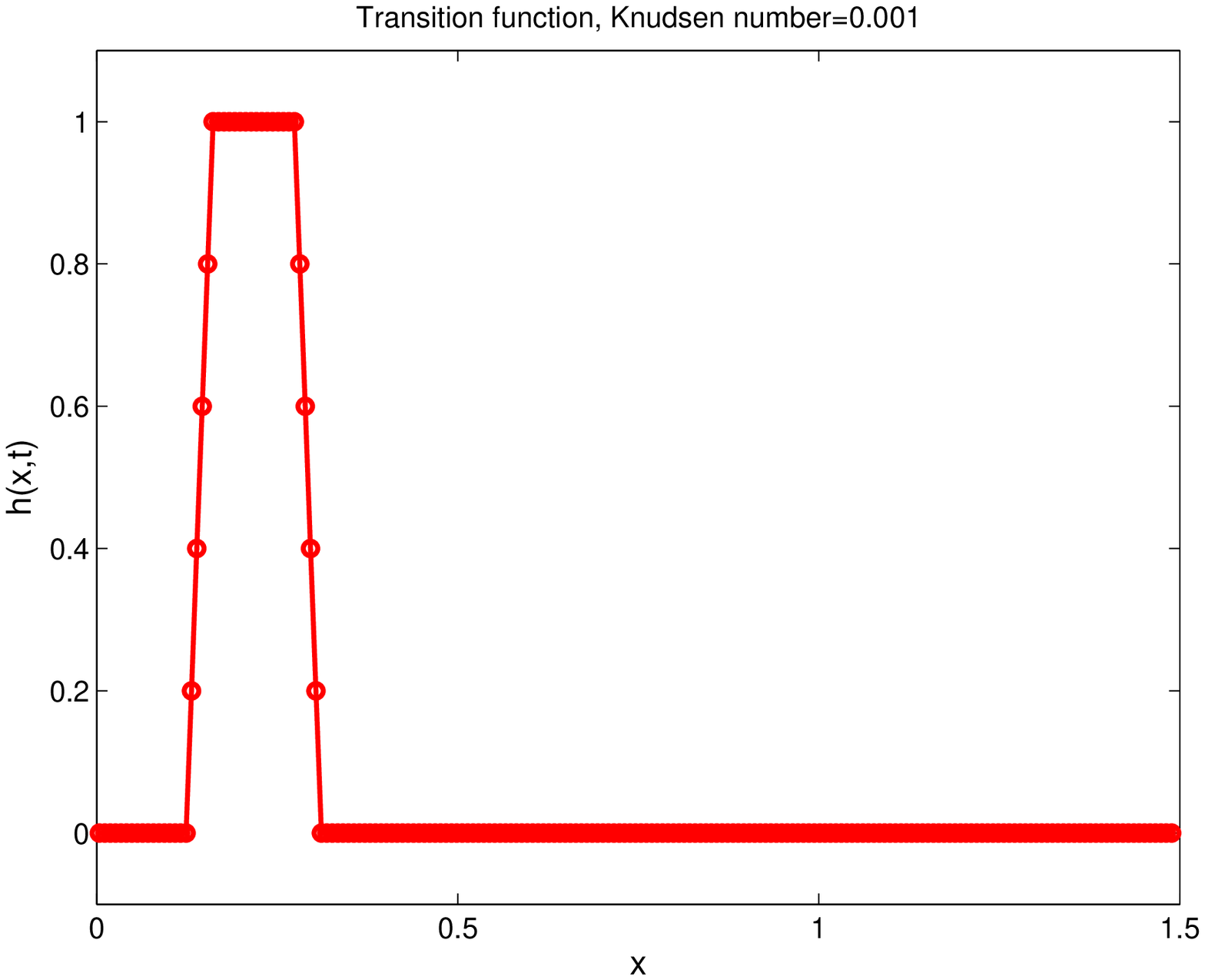}
\includegraphics[scale=0.27]{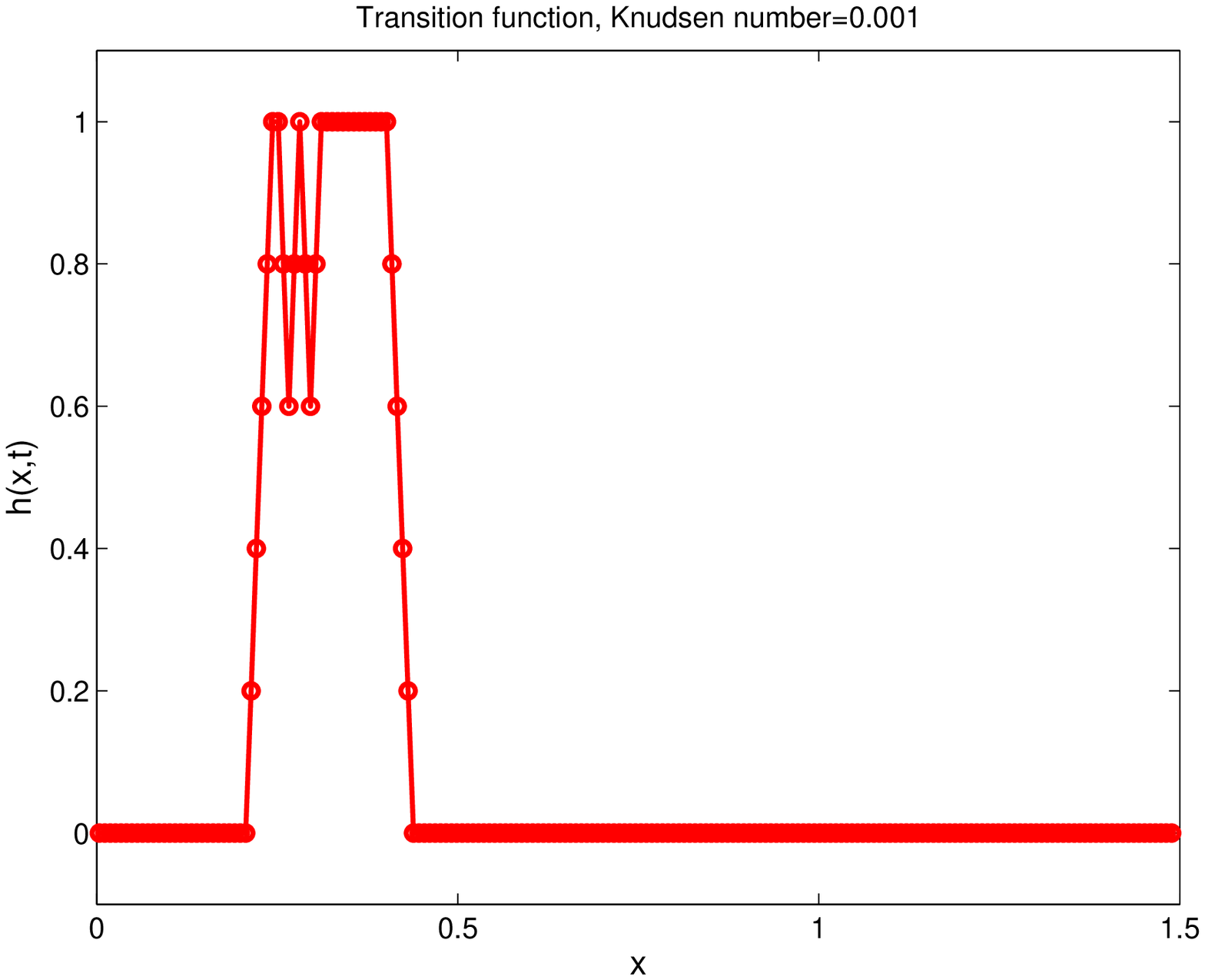}
\caption{Unsteady Shock Test: Transition function at $t=0.05$
(left), $t=0.10$ (right) and $t=0.15$ (bottom). Knudsen number
$\varepsilon=10^{-3}$.} \label{ST23}
\end{center}
\end{figure}

We repeat the test with the same initial data $\varrho=1$, $u=-2$
and $T=4$, but changing the collision frequency. The relaxation
parameter is $\varepsilon=10^{-1}$ in the first case,
$\varepsilon=10^{-2}$ in the second case and $\varepsilon=10^{-3}$
in the last case. The thickness of the transition regions is fixed
for every test depending on the frequency of collisions: ten cells
for $\varepsilon=10^{-1}$ and $\varepsilon=10^{-2}$ while five cells
for $\varepsilon=10^{-3}$.

In figure \ref{ST24} we have reported the number of particles in
time for the DSMC method (left) and the coupling method (right) for
the different values of the frequency. As expected this number is
independent from the choice of $\varepsilon$ for the DSMC method
while in the case of the coupling the quantity of particles used
varies strongly. This number is a measure of the computational cost
of the method, the coupling procedure being computationally not
expensive as well as the moment guided method and the solution of
the fluid equations. In figure \ref{ST03} the transition function is
depicted for three different times in the case of
$\varepsilon=10^{-1}$. The same function is depicted in figure
\ref{ST13} for $\varepsilon=10^{-2}$ and in figure \ref{ST23} for
$\varepsilon=10^{-3}$. These figures show the dynamics of the
kinetic and fluid regions in time and the capability of the method
to follow the shock which moves in the domain. We remark also that
the thickness of the kinetic region decreases when the collision
frequency increases as expected.

\begin{figure}[h!]
\begin{center}
\includegraphics[scale=0.39]{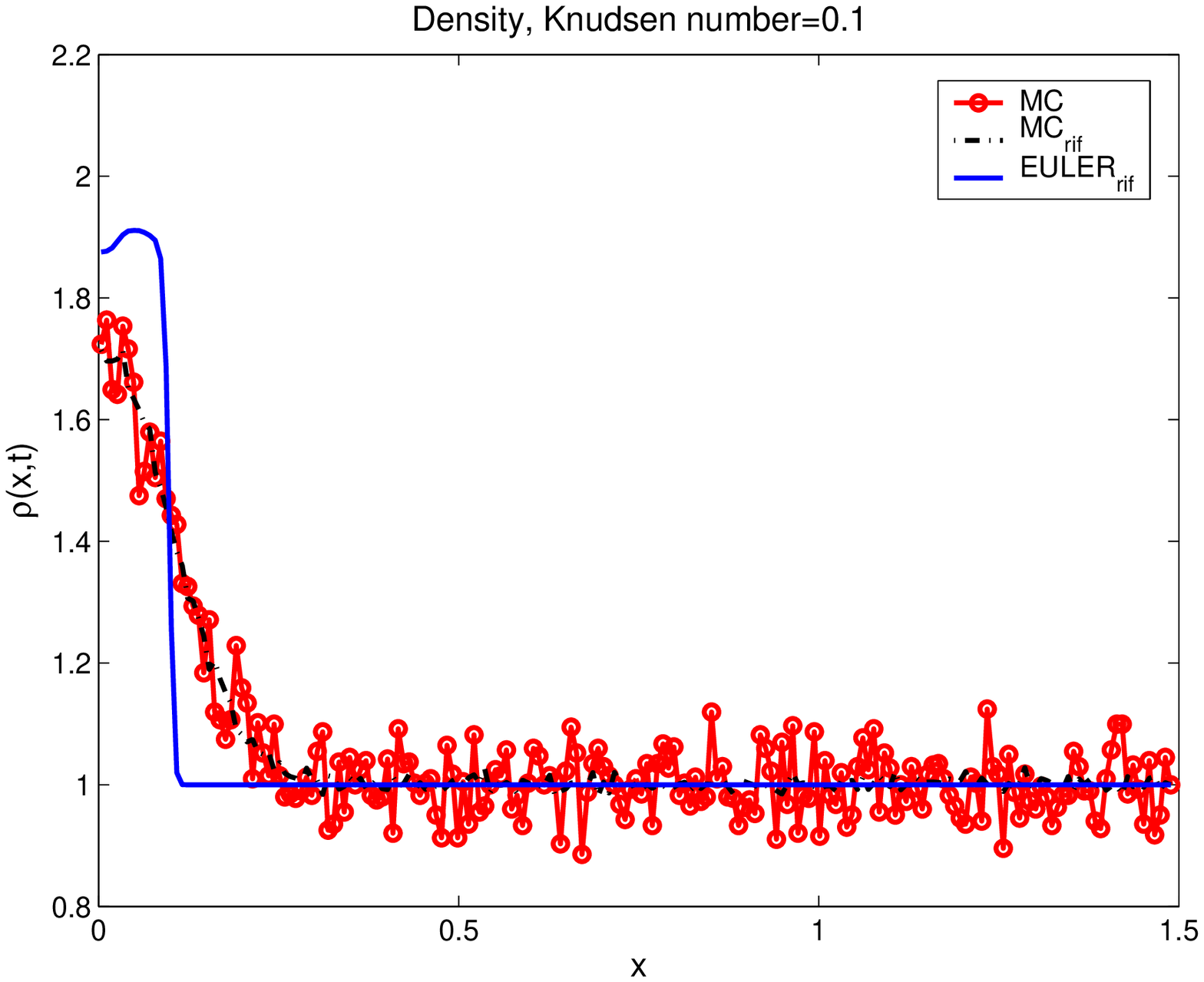}
\includegraphics[scale=0.39]{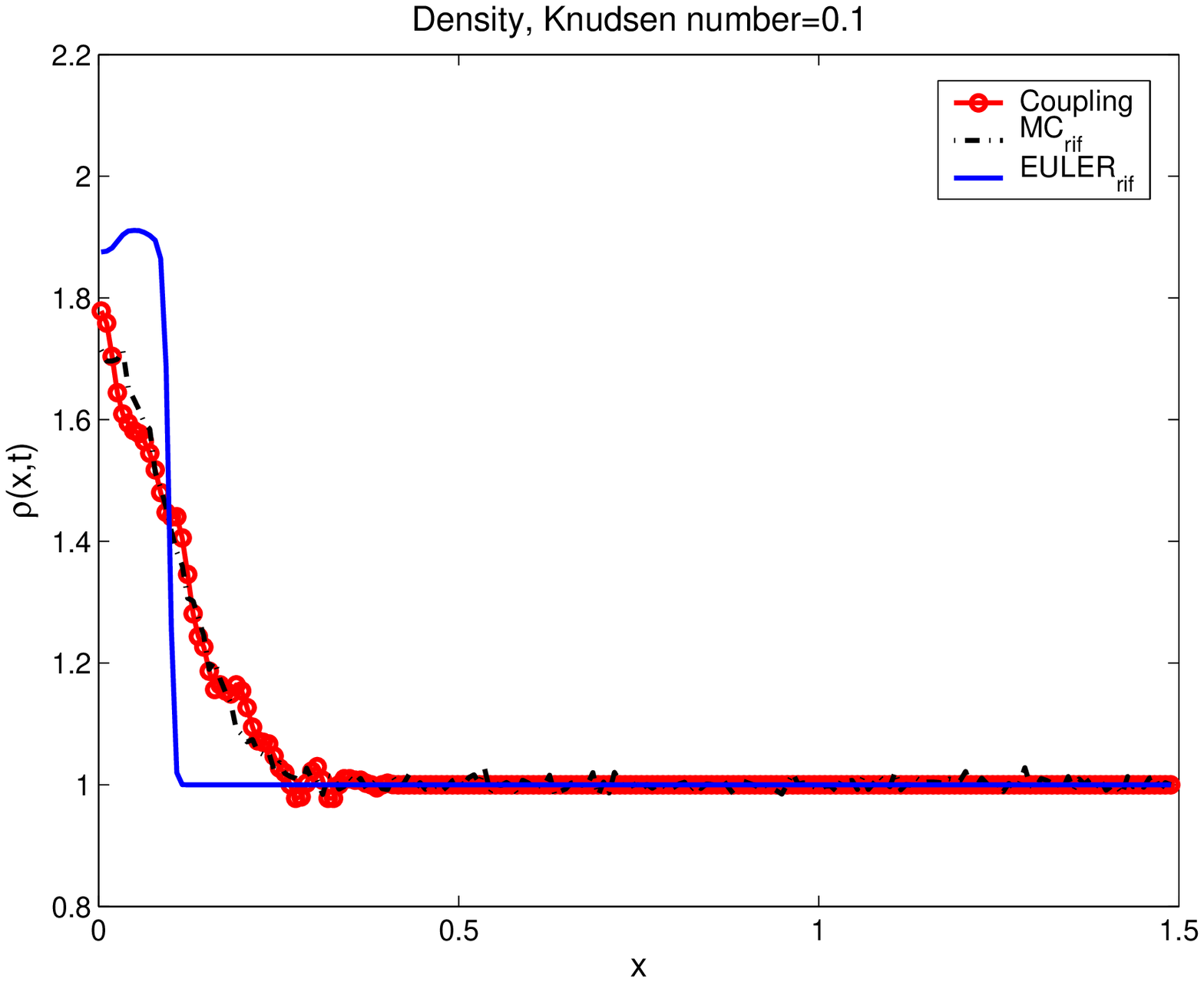}\\
\includegraphics[scale=0.39]{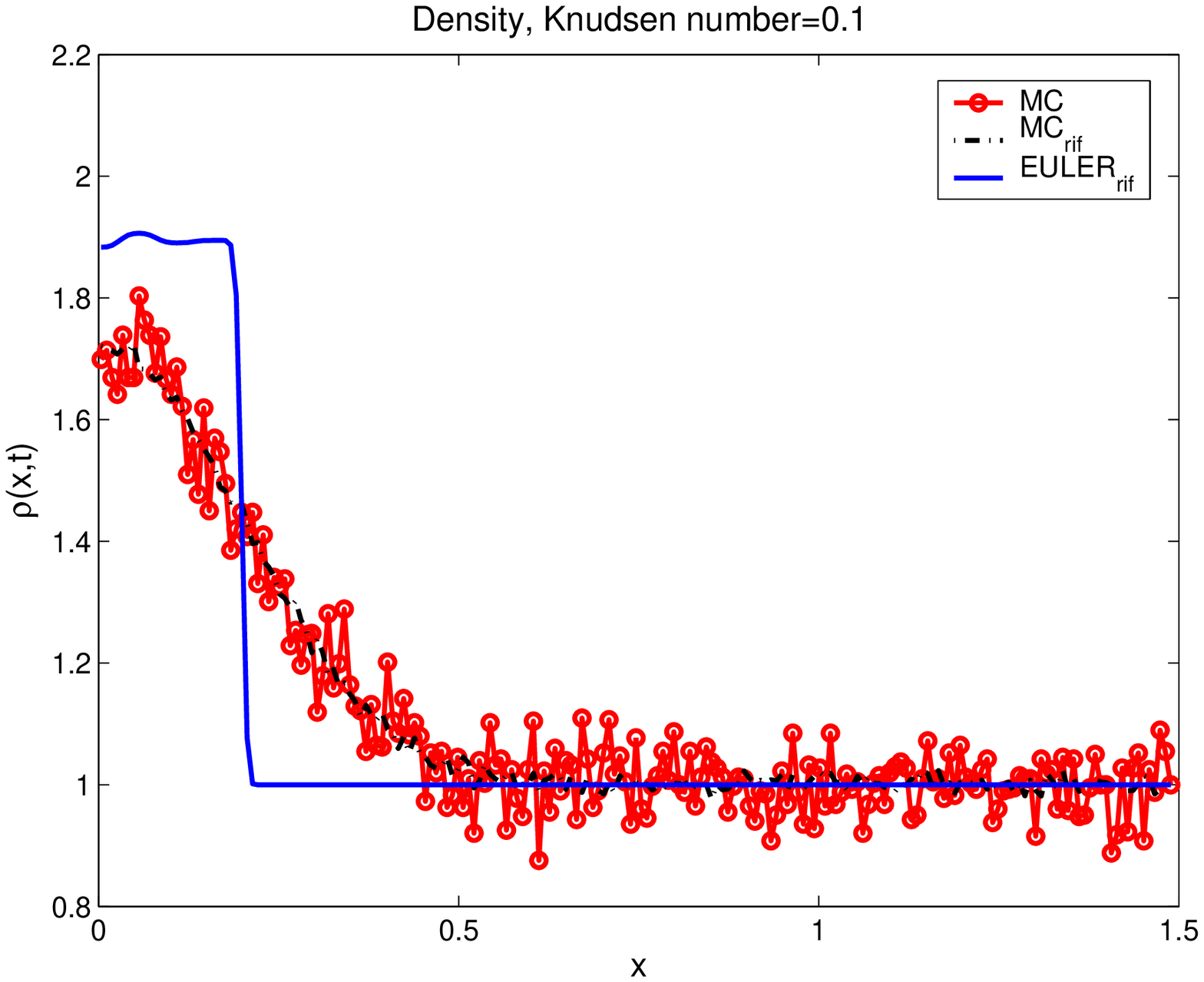}
\includegraphics[scale=0.39]{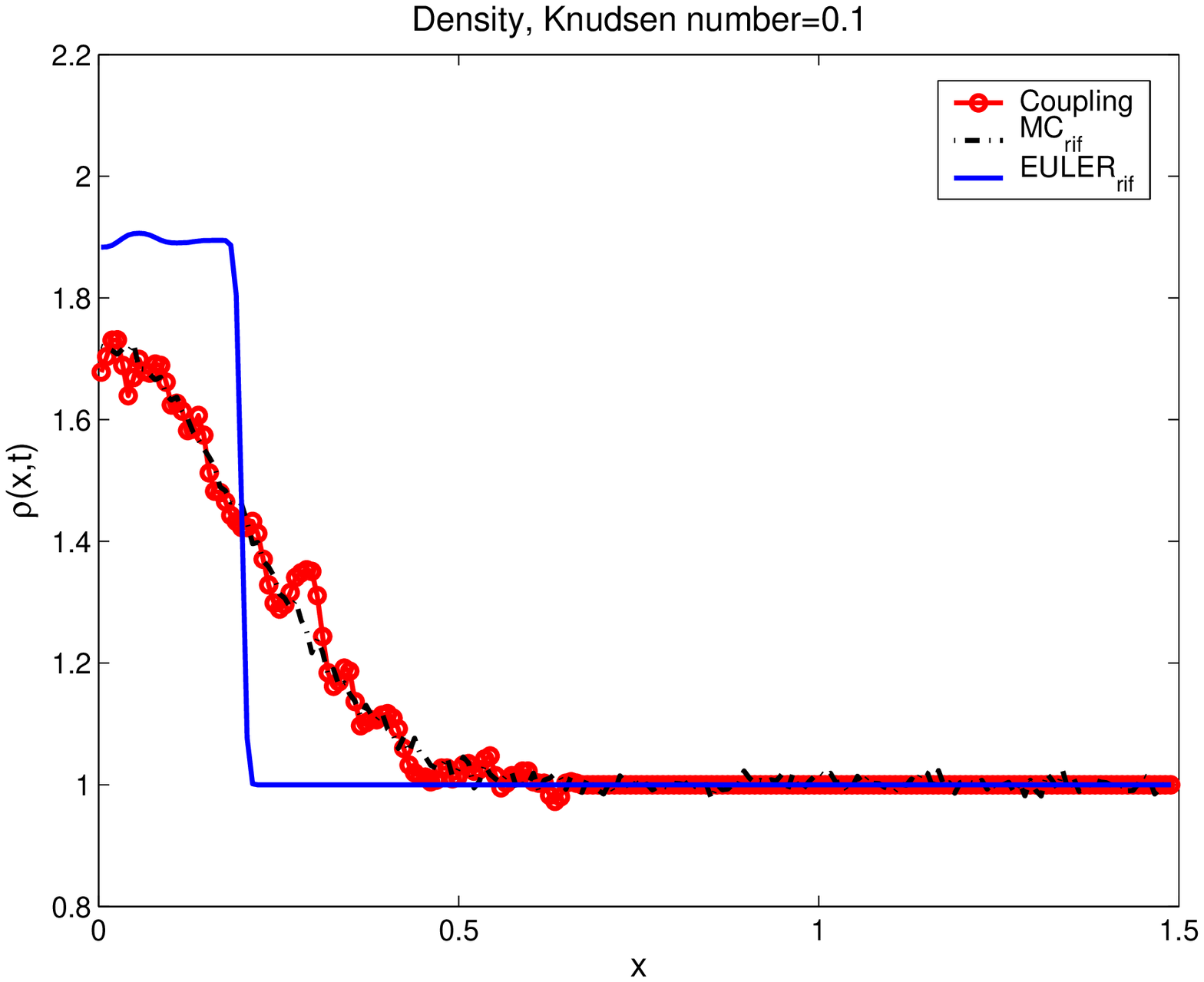}\\
\includegraphics[scale=0.39]{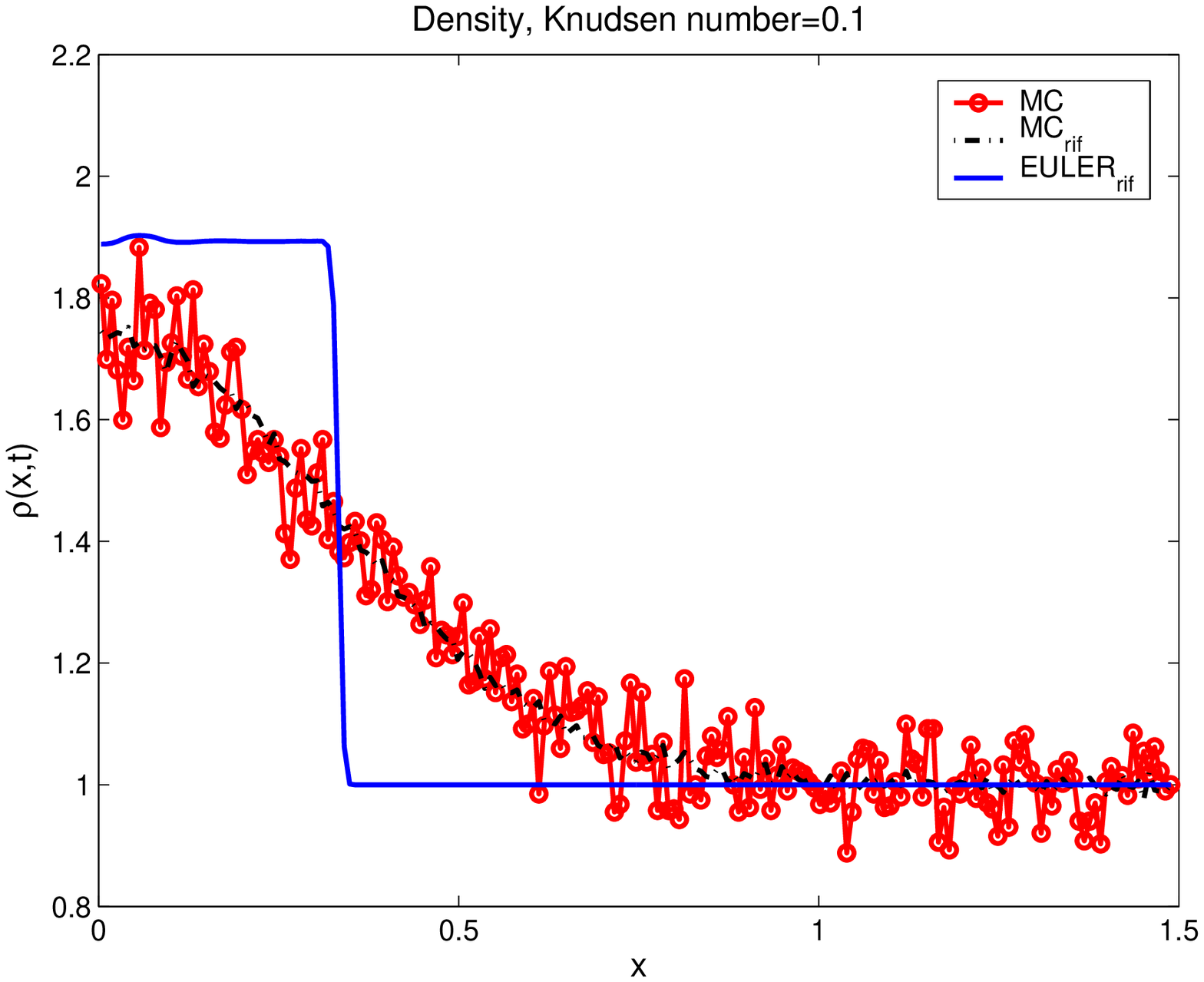}
\includegraphics[scale=0.39]{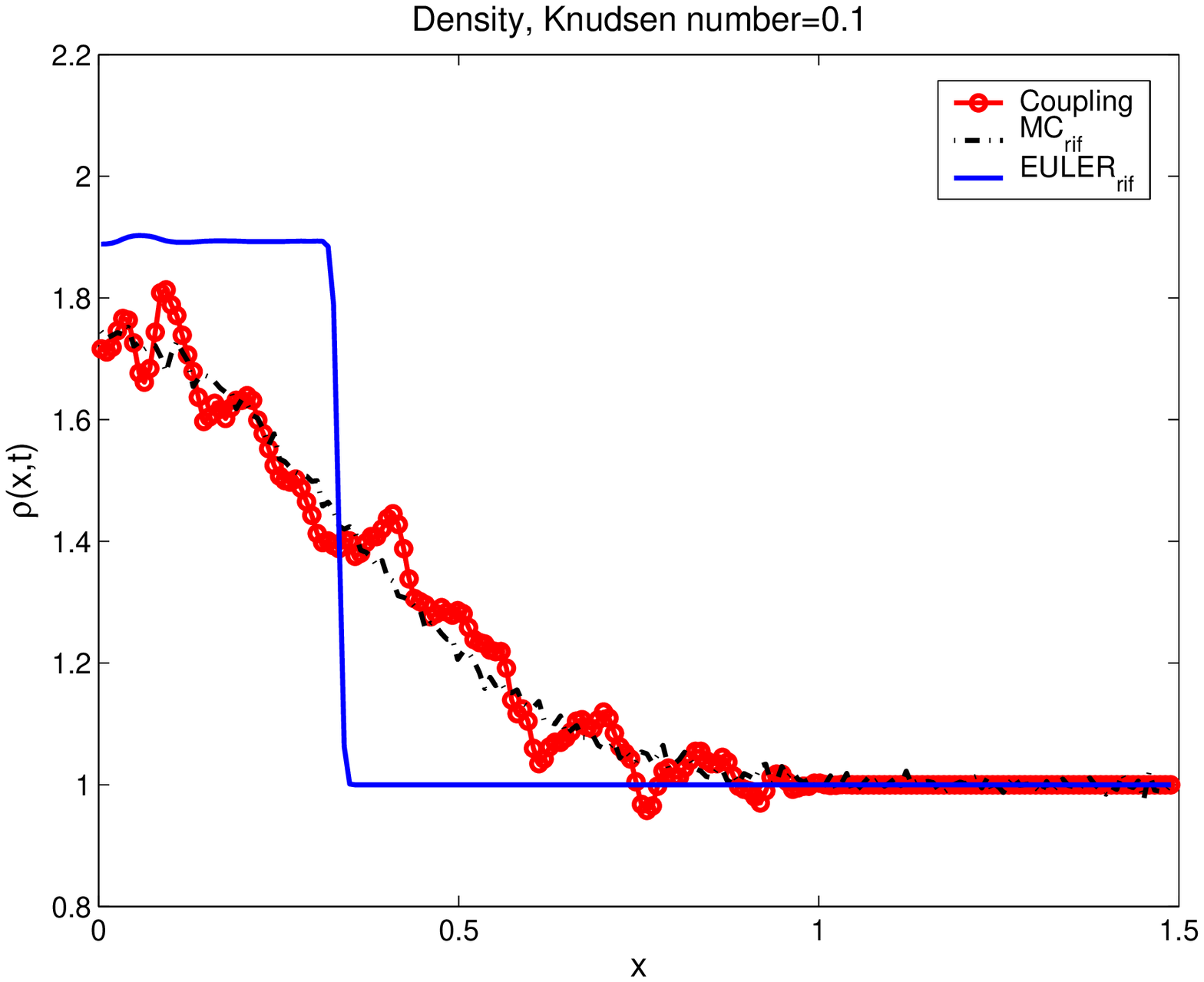}
\caption{Unsteady Shock Test: Solution at $t=0.05$ (top), $t=0.10$
(middle) and $t=0.15$ (bottom) for the density. MC method (left),
Coupling DSMC-Fluid method (right). Knudsen number
$\varepsilon=10^{-1}$. Reference solution (dotted line), Euler solution (continuous line),
DSMC-Fluid or DSMC (circles plus continuous line).}\label{ST00}
\end{center}
\end{figure}

In figure \ref{ST00} we have reported the density on the left for
the DSMC method and the on the right for the DSMC/Fluid coupling. In
this figure $\varepsilon=10^{-1}$. From top to bottom, time
increases from $t=0.05$ (top) to $t=0.15 $ (bottom) with $t=0.1$ in
the middle. In each of the plots we reported the solution computed
with our coupling method or with the DSMC method. The solution of
the compressible Euler equation and a reference solution are also
reported. The reference solution has been computed with a DMSC
method where the number of particles is taken very high to reduce
the statistical fluctuations

\begin{figure}
\begin{center}
\includegraphics[scale=0.39]{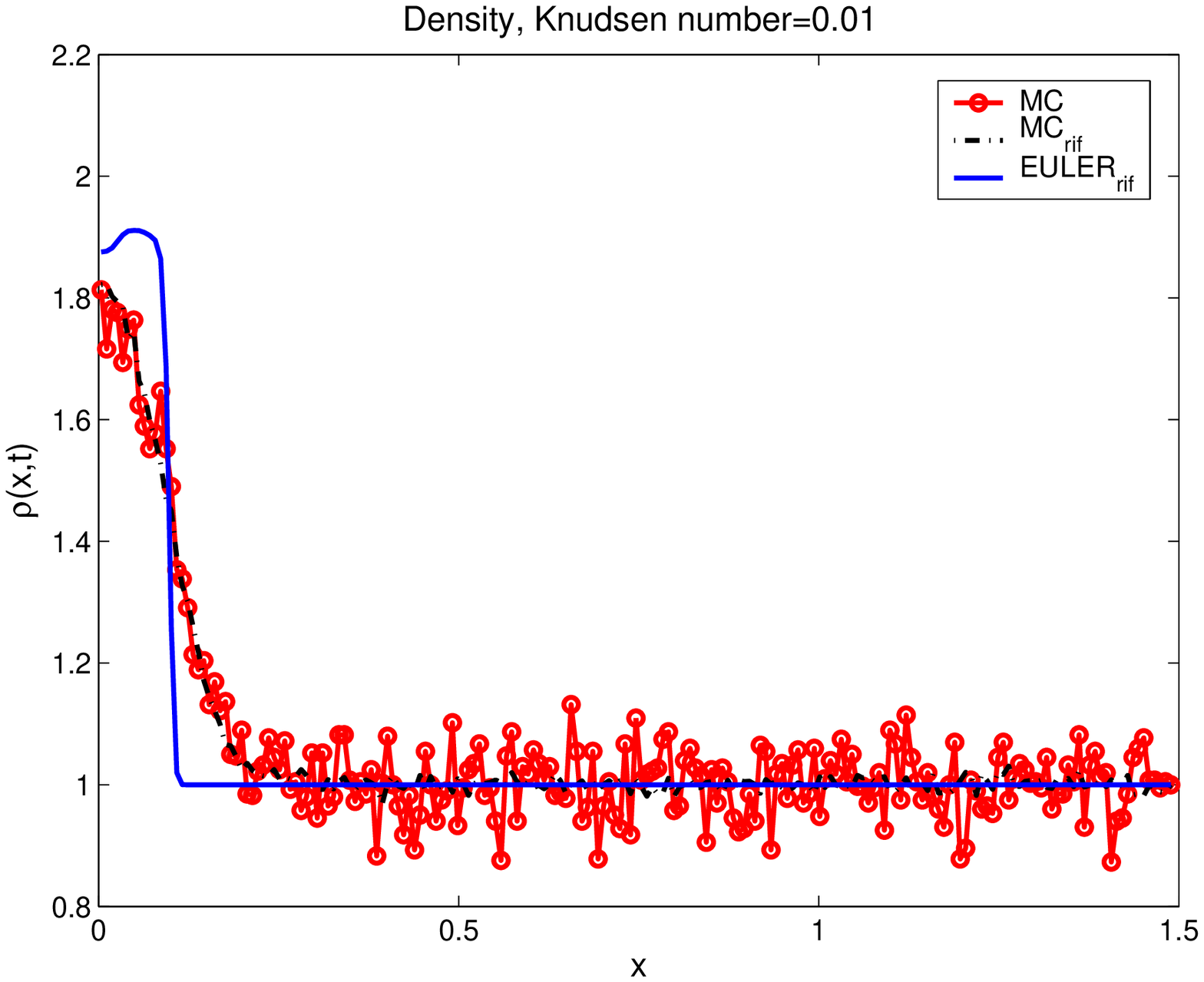}
\includegraphics[scale=0.39]{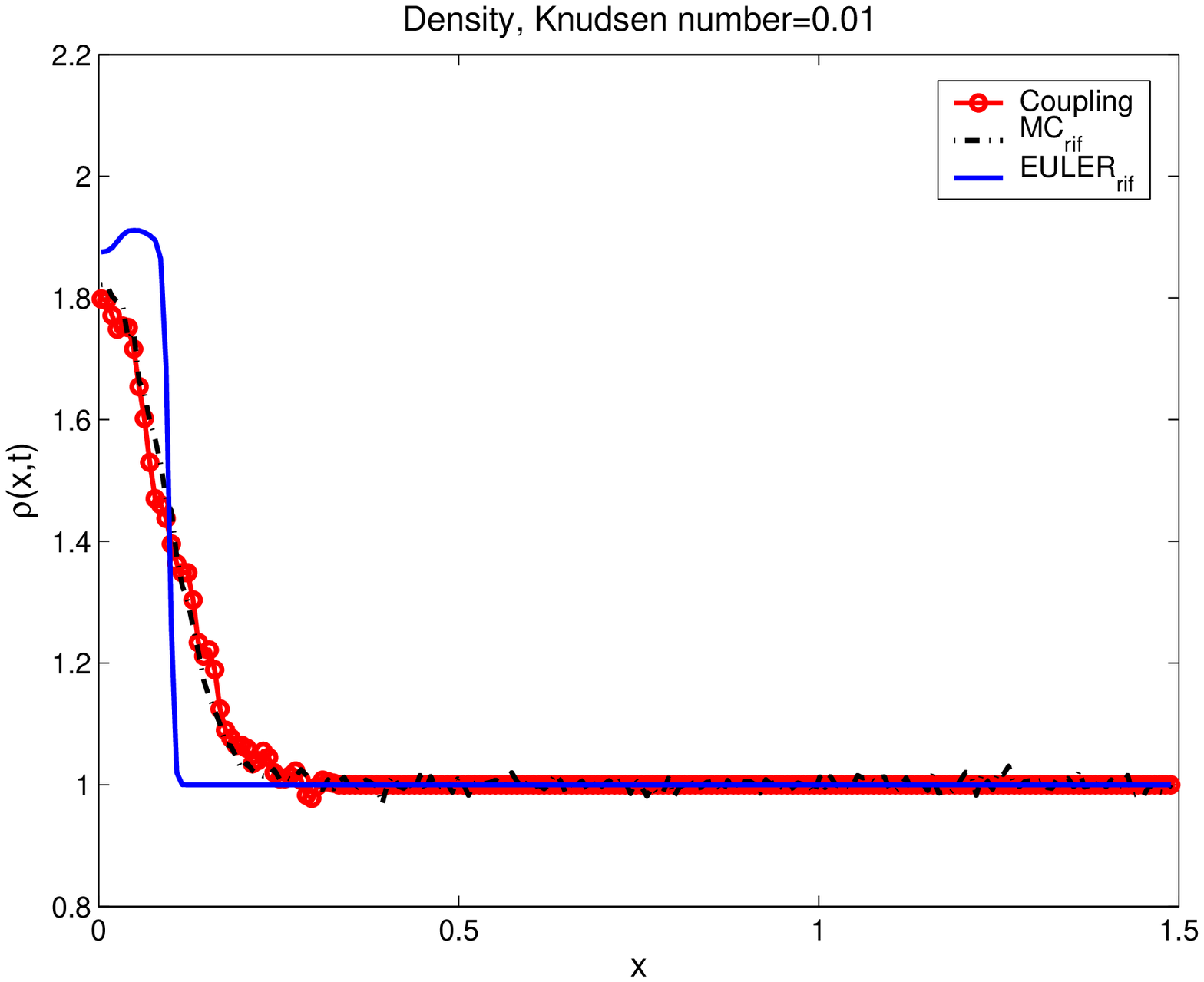}\\
\includegraphics[scale=0.39]{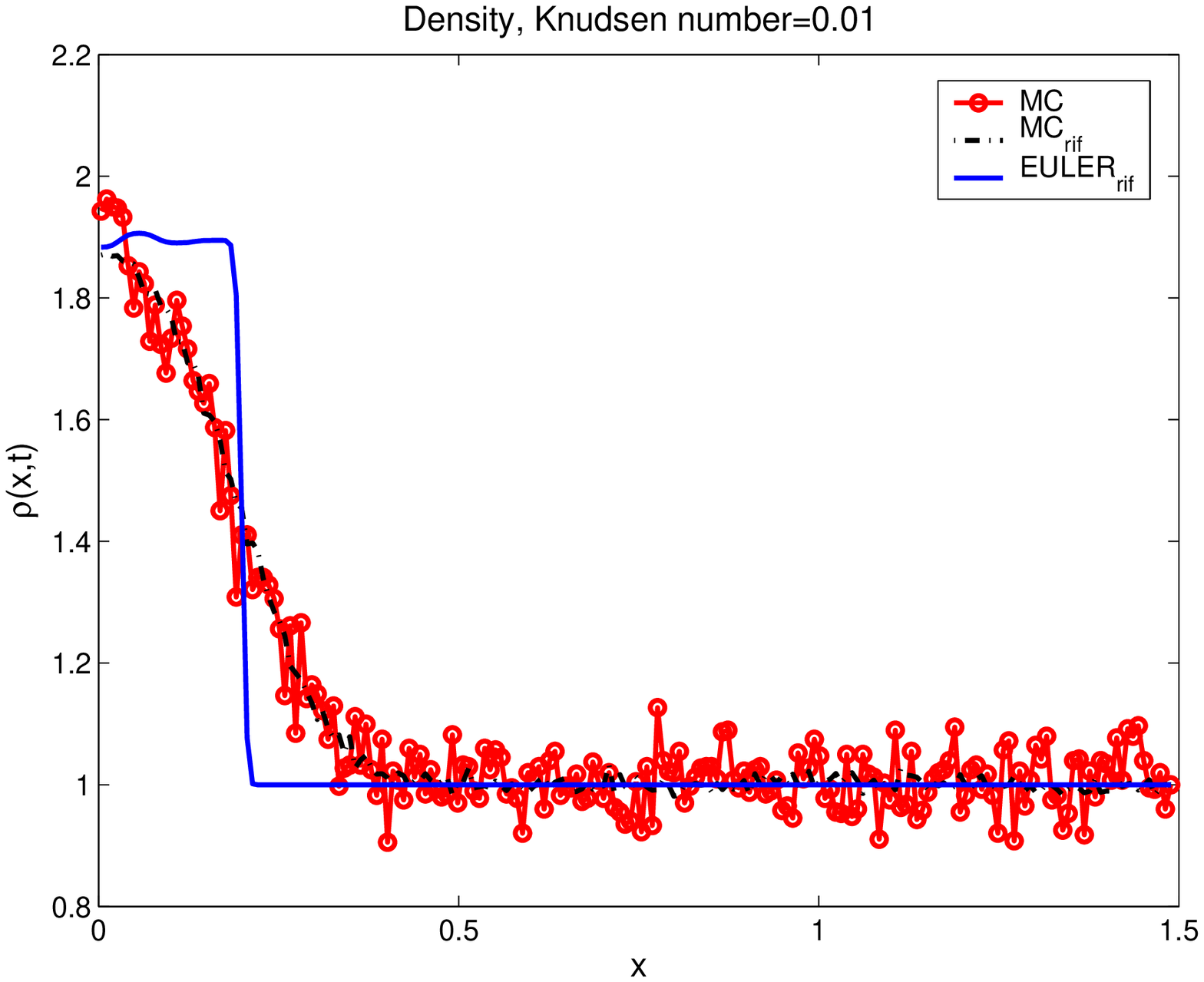}
\includegraphics[scale=0.39]{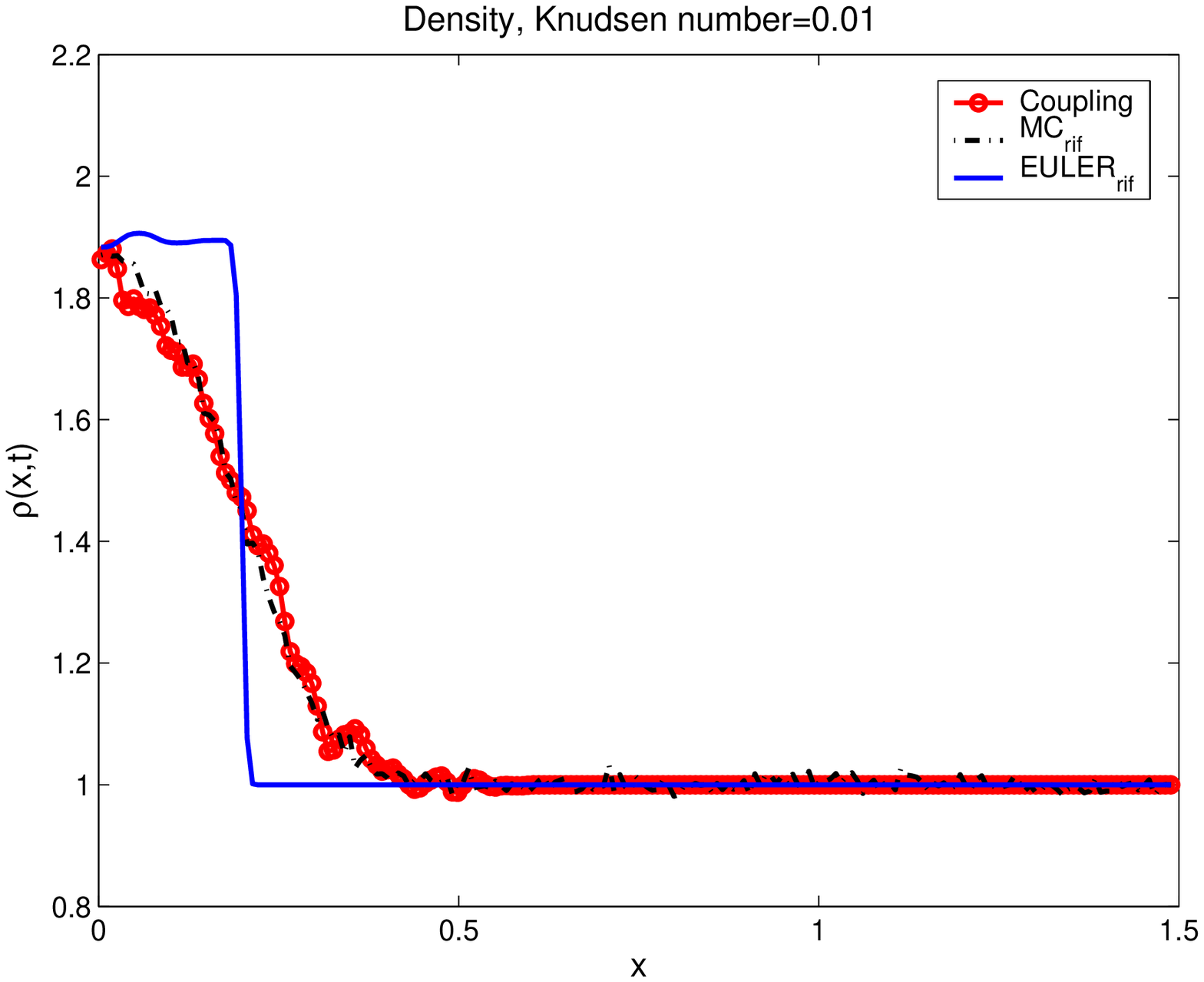}\\
\includegraphics[scale=0.39]{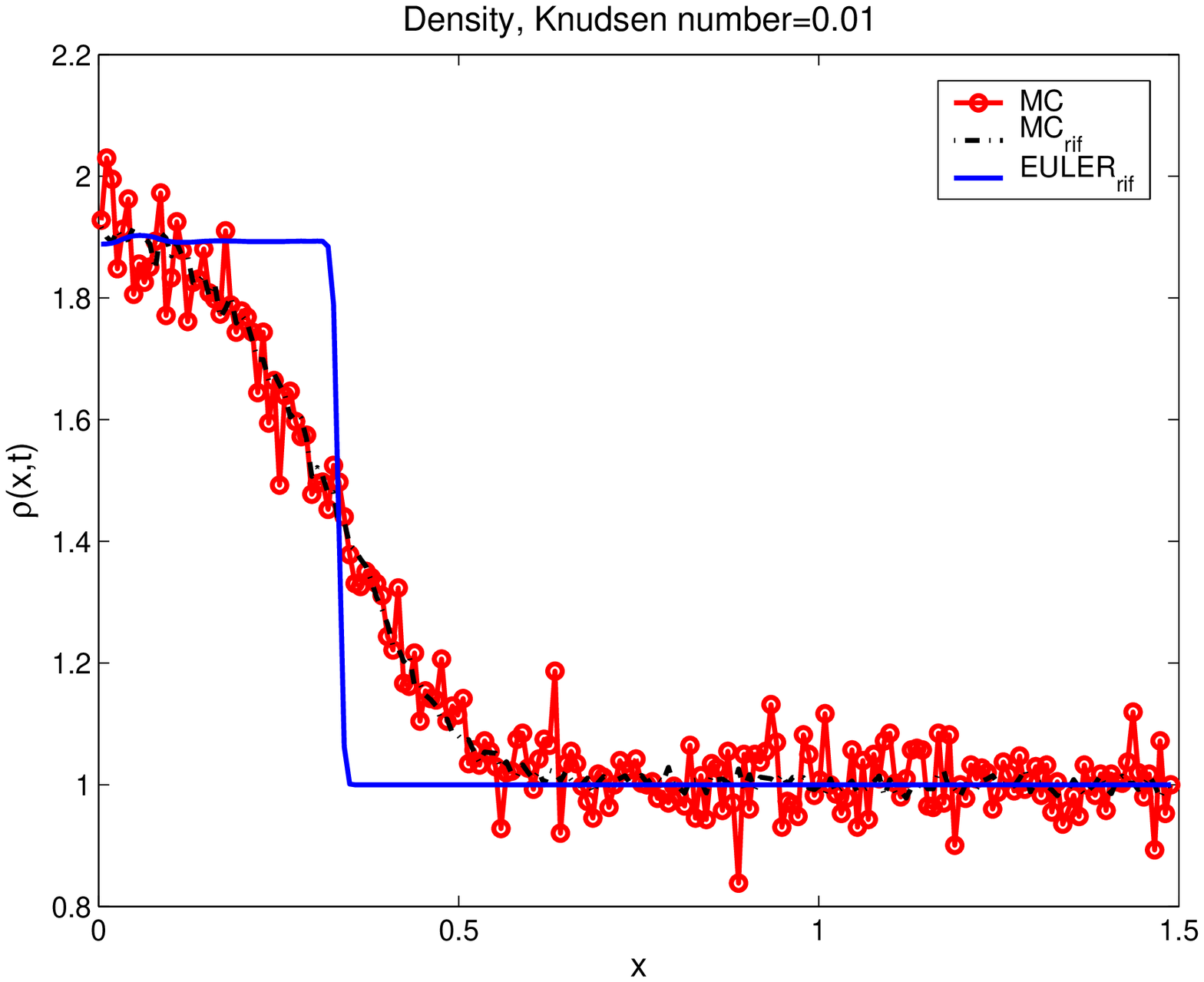}
\includegraphics[scale=0.39]{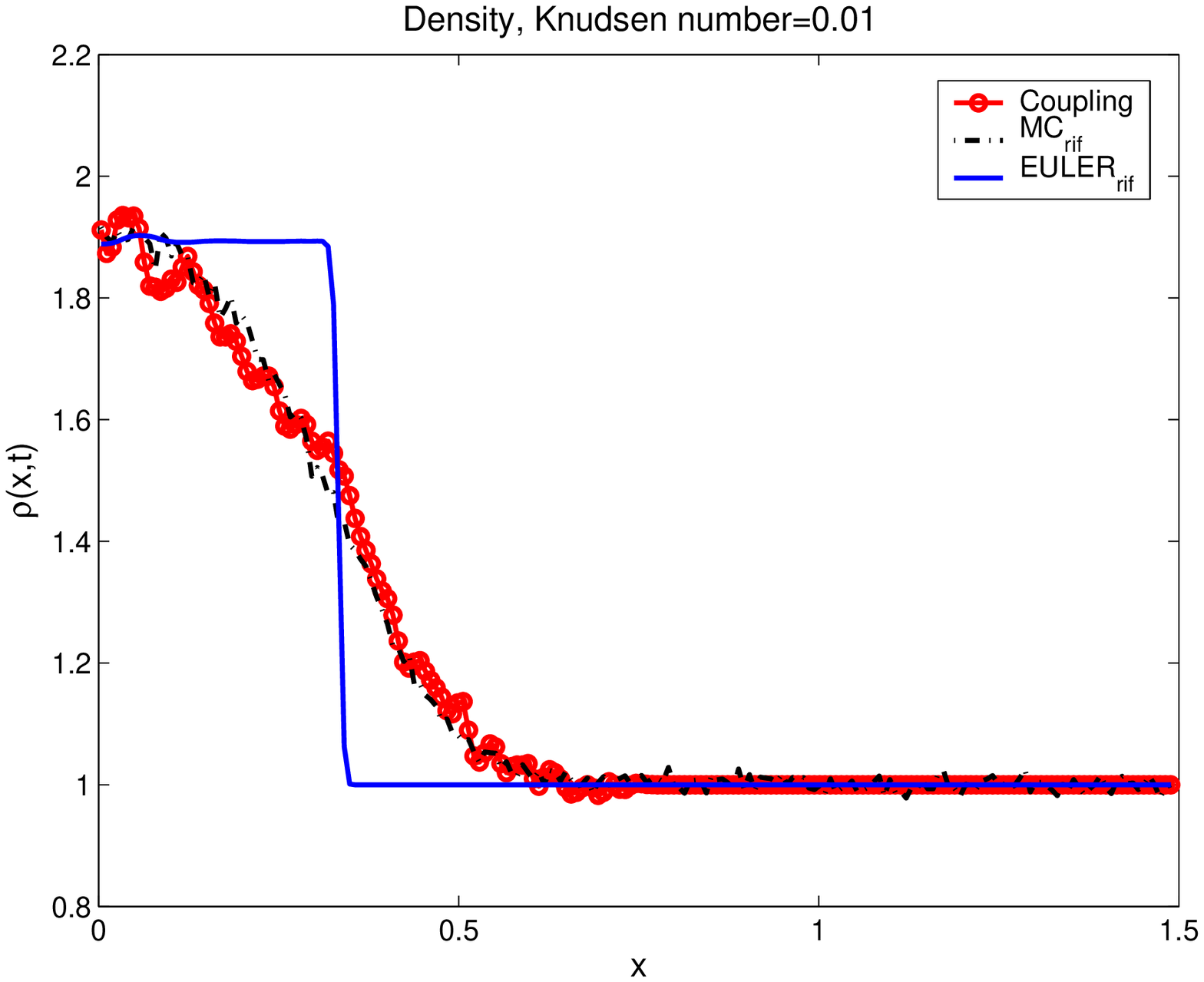}
\caption{Unsteady Shock Test: Solution at $t=0.05$ (top), $t=0.10$
(middle) and $t=0.15$ (bottom) for the density. MC method (left),
Coupling DSMC-Fluid method (right). Knudsen number
$\varepsilon=10^{-2}$. Reference solution (dotted line), Euler solution (continuous line),
DSMC-Fluid or DSMC (circles plus continuous line).}\label{ST10}
\end{center}
\end{figure}

\begin{figure}
\begin{center}
\includegraphics[scale=0.39]{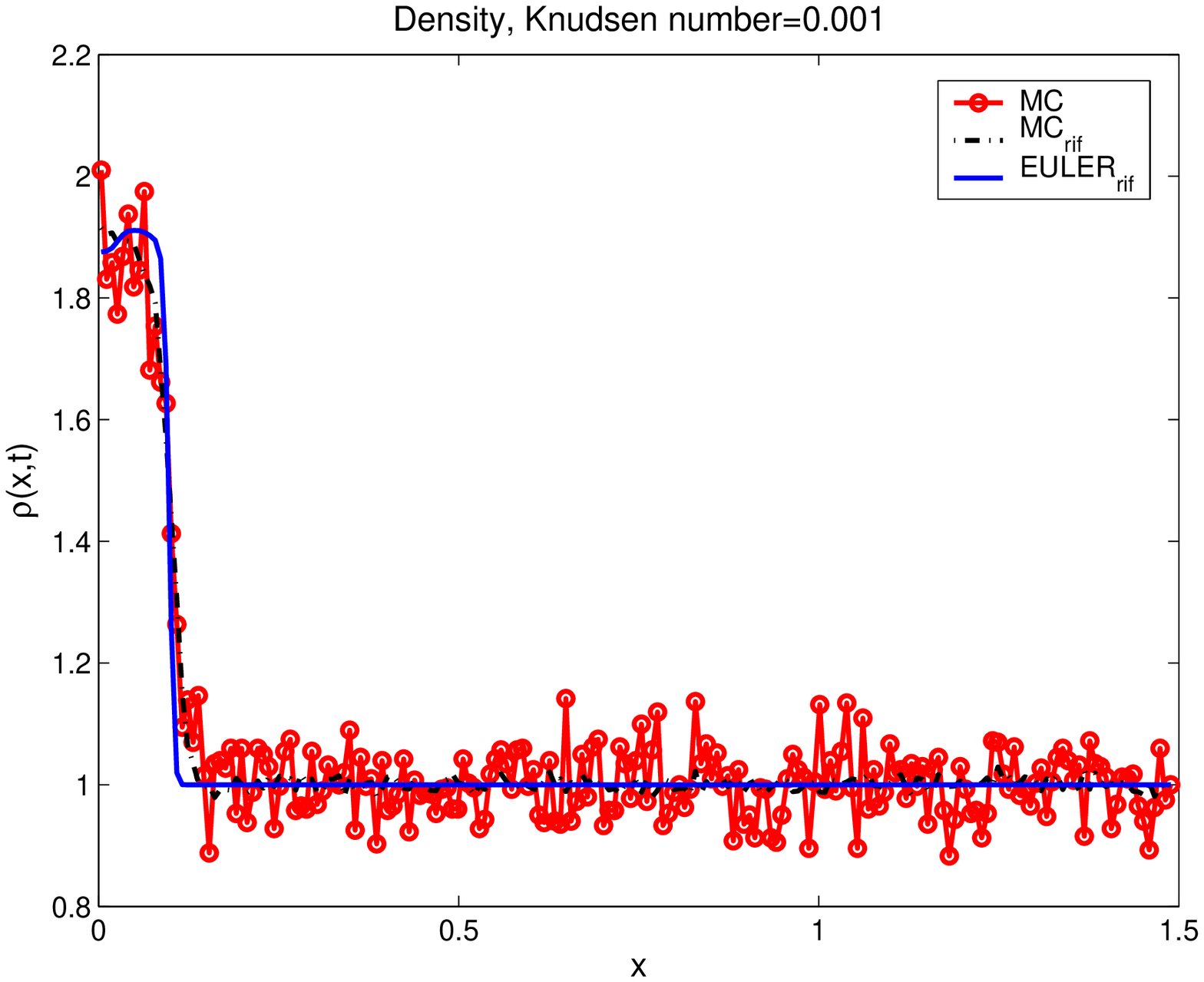}
\includegraphics[scale=0.39]{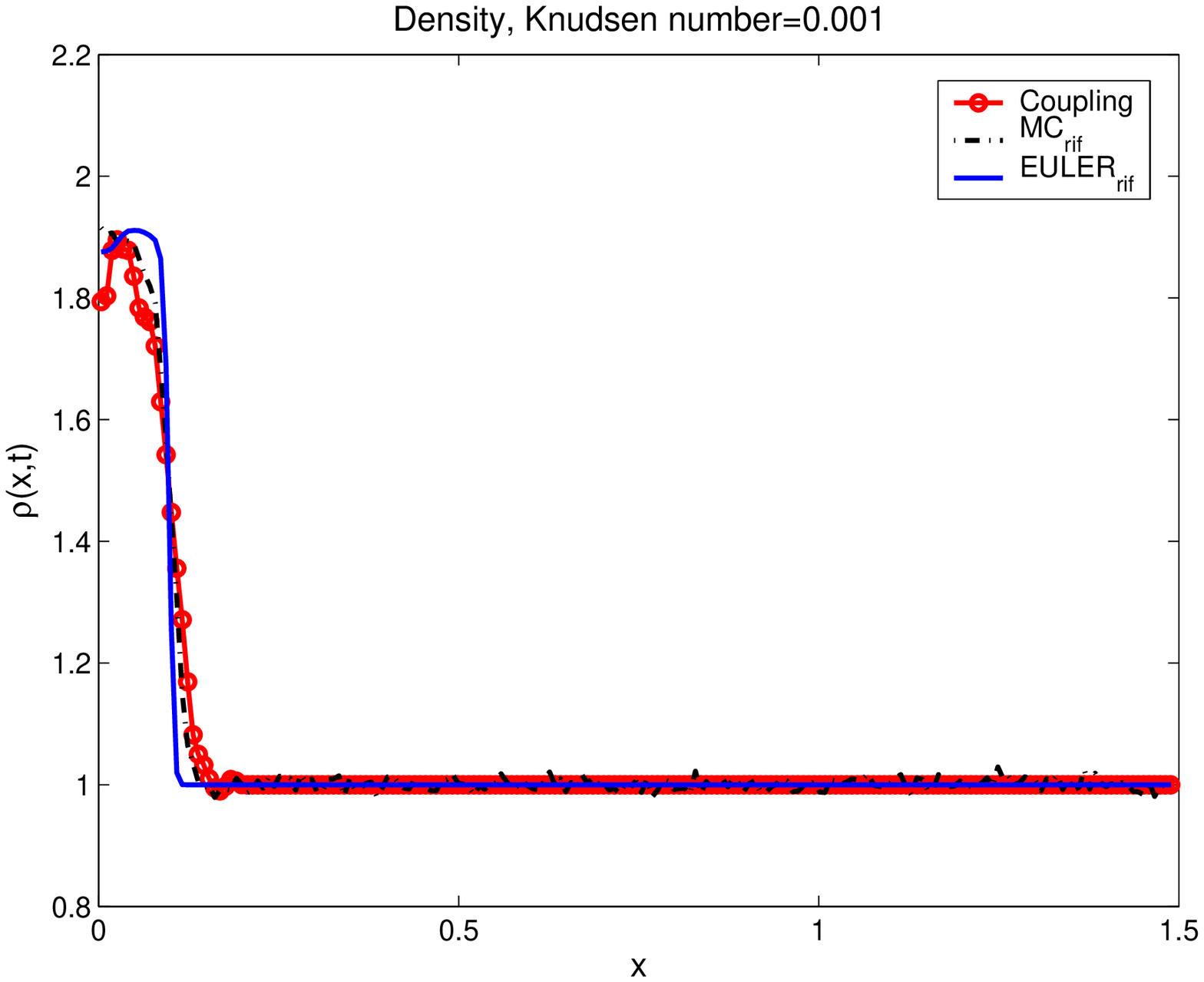}\\
\includegraphics[scale=0.39]{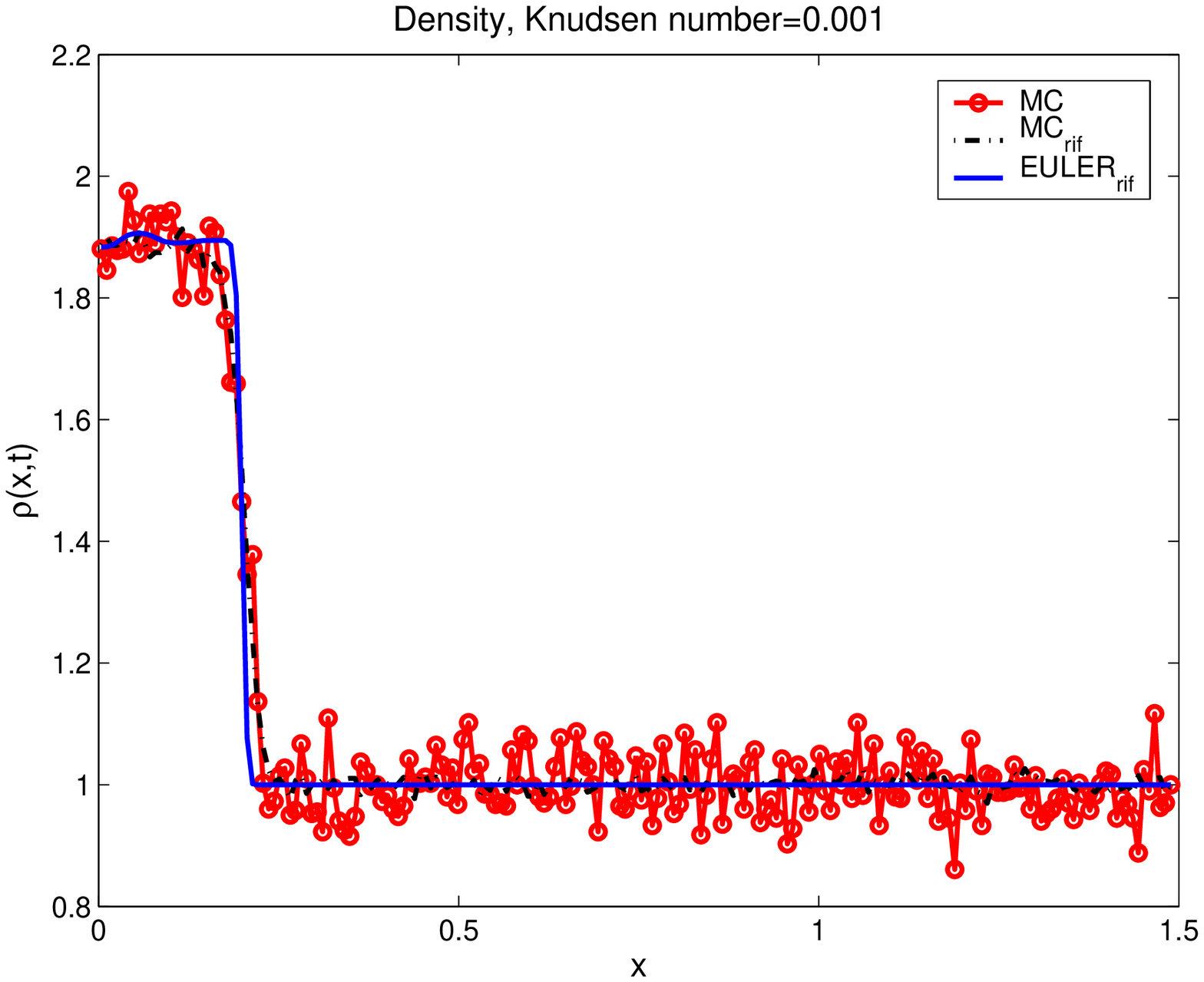}
\includegraphics[scale=0.39]{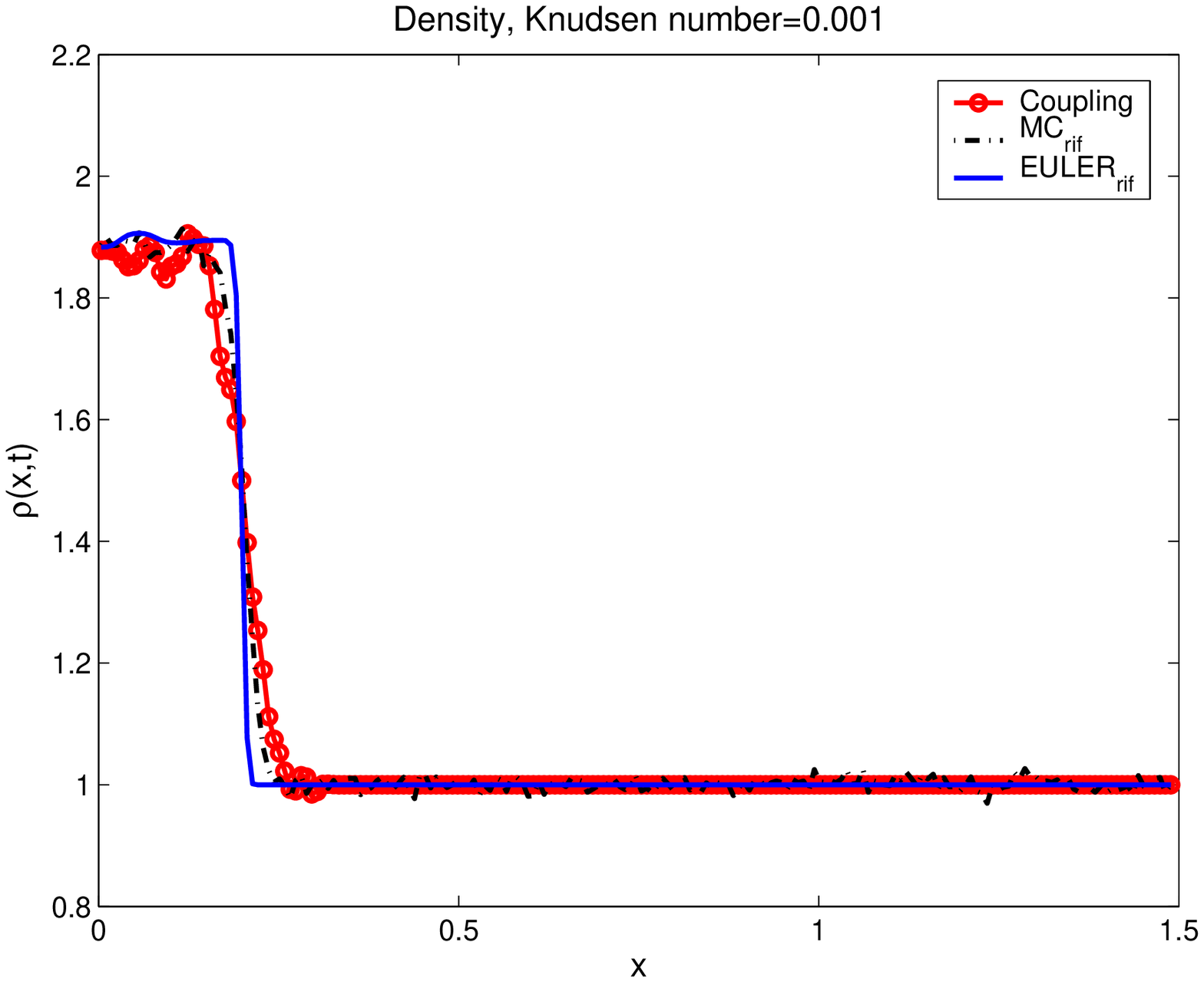}\\
\includegraphics[scale=0.39]{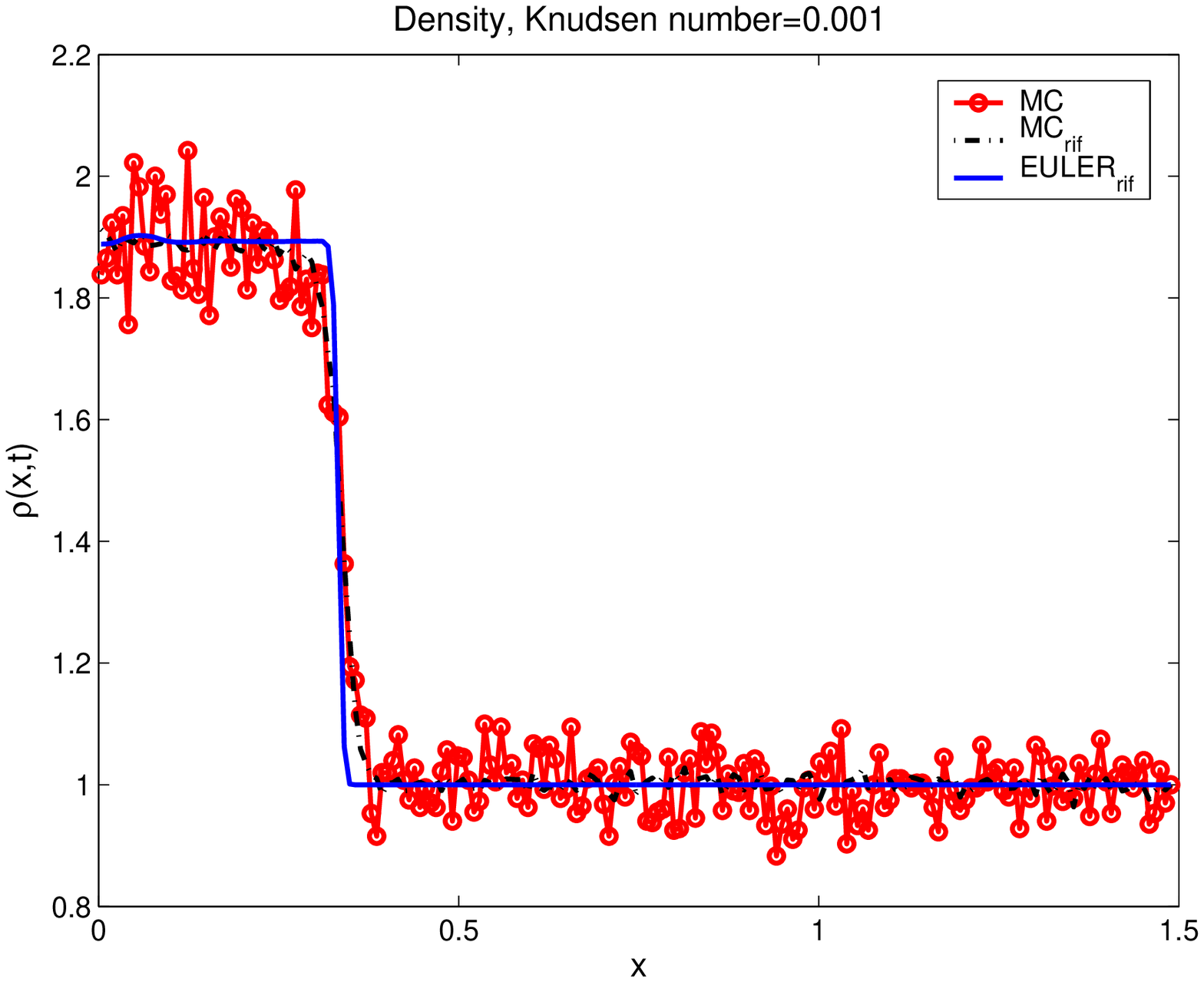}
\includegraphics[scale=0.39]{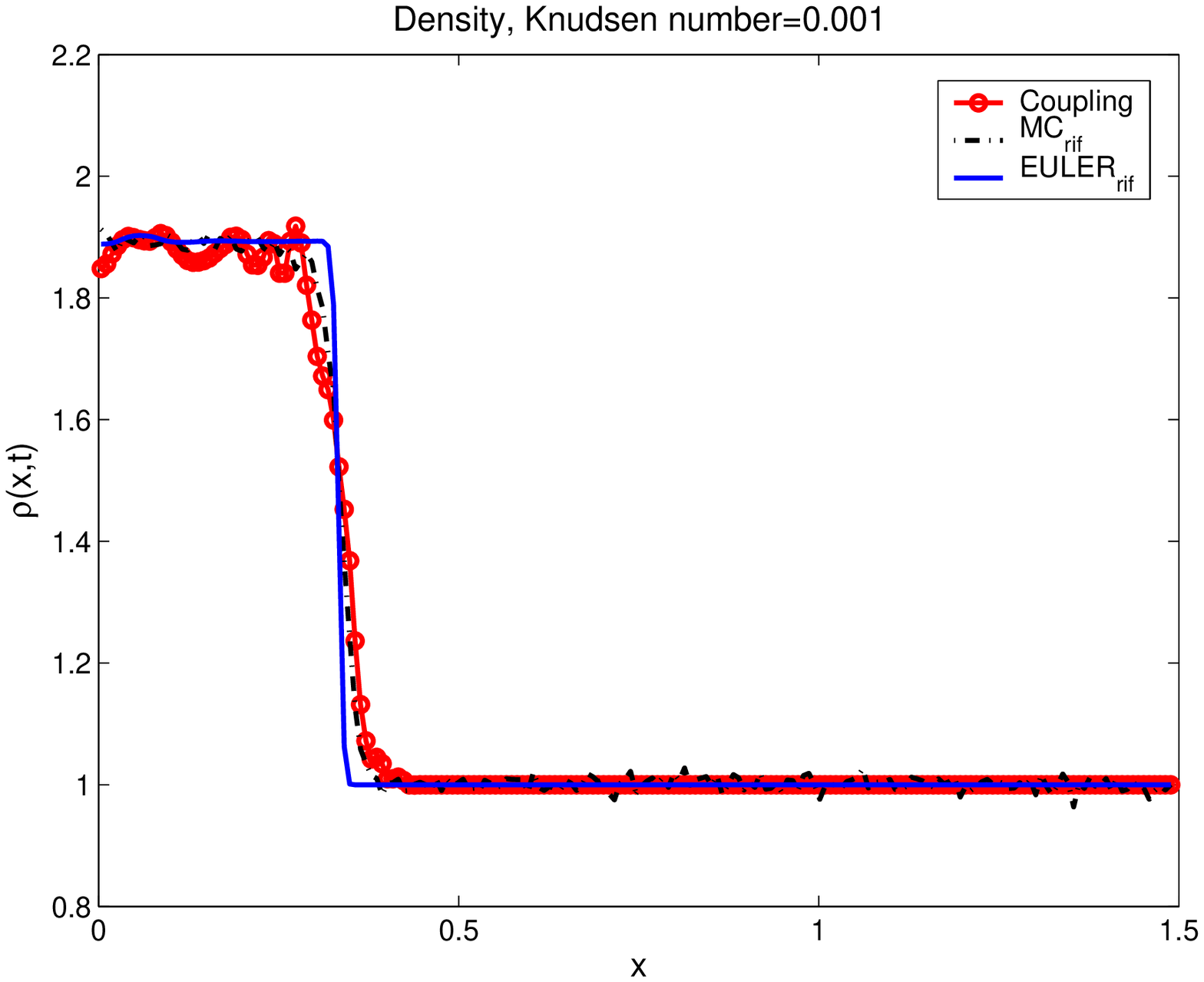}
\caption{Unsteady Shock Test: Solution at $t=0.05$ (top), $t=0.10$
(middle) and $t=0.15$ (bottom) for the density. MC method (left),
Coupling DSMC-Fluid method (right). Knudsen number
$\varepsilon=10^{-3}$. Reference solution (dotted line), Euler solution (continuous line),
DSMC-Fluid or DSMC (circles plus continuous line).}\label{ST20}
\end{center}
\end{figure}

In figure \ref{ST10} the density is reported for
$\varepsilon=10^{-2}$ again for DSMC on the left and DSMC/Fluid
coupling on the right. Finally we report the same results in the case of
$\varepsilon=10^{-3}$ in figure \ref{ST20}.

During the simulation on the left boundary, the transition function
$h$ increases from zero to one, which means that the solution is
computed with the DSMC scheme, while in the rest of the domain the
solution is still computed with the fluid scheme ($h=0$). When the
shock starts to move towards the right the kinetic region increases.
We observe that, in the case of $\varepsilon=10^{-3}$, when the
shock is far from the left boundary the gas returns to
thermodynamical equilibrium and thus automatically $h$ becomes equal
to $0$. This is an interesting results because it shows that the
method is able to recover an equilibrium regime after a time span in
which it was in non equilibrium. While this results is not
surprising when a deterministic scheme is used to solve the kinetic
model, this is not obvious in the case of a DSMC solution of the
kinetic equation because of the fluctuations which prevent the
good evaluation of the macroscopic quantities needed to calculate
the breakdown of the fluid model.

\subsection{Sod tests}
\label{subsec_sod}

\begin{figure}
\begin{center}
\includegraphics[scale=0.27]{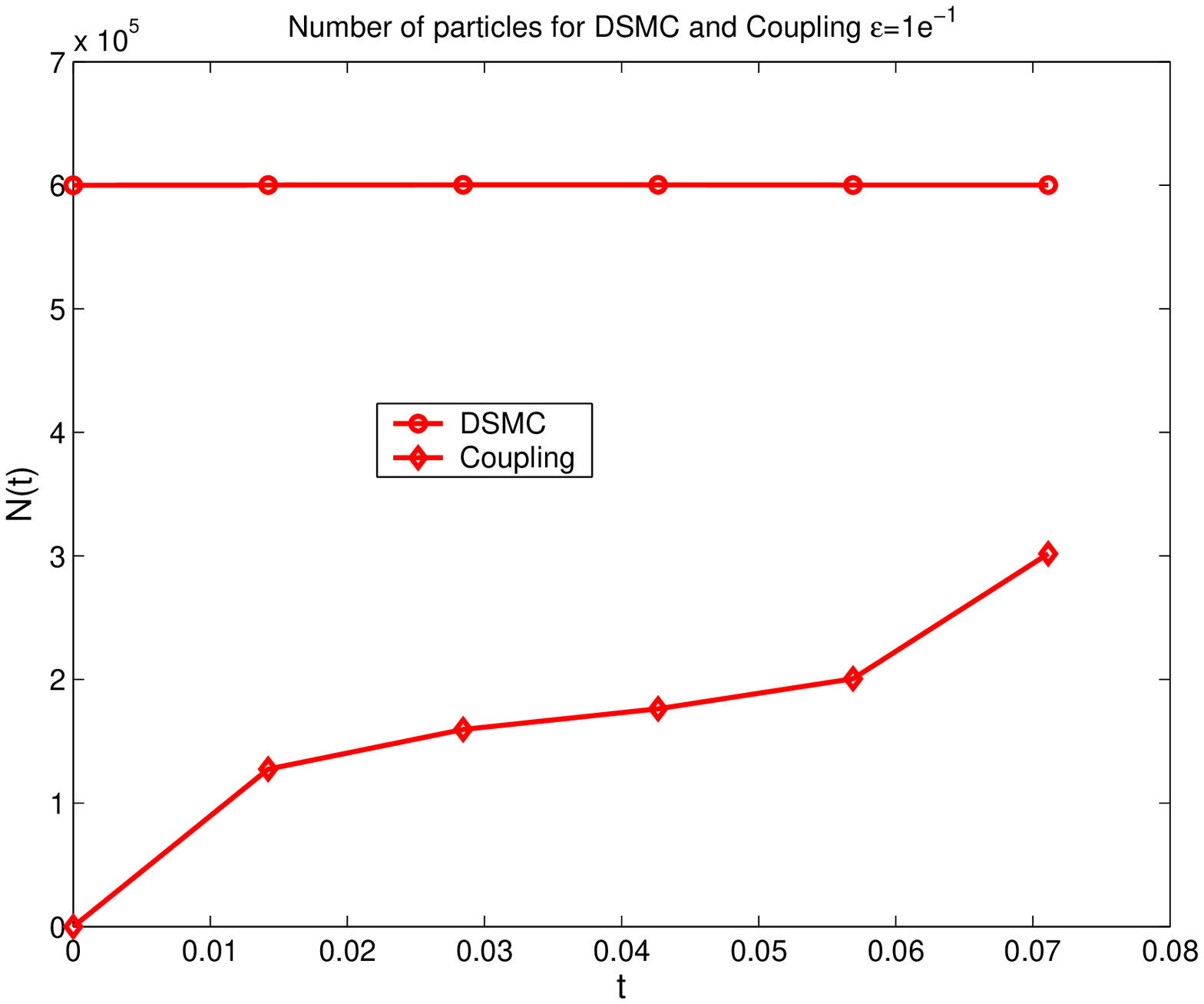}
\includegraphics[scale=0.27]{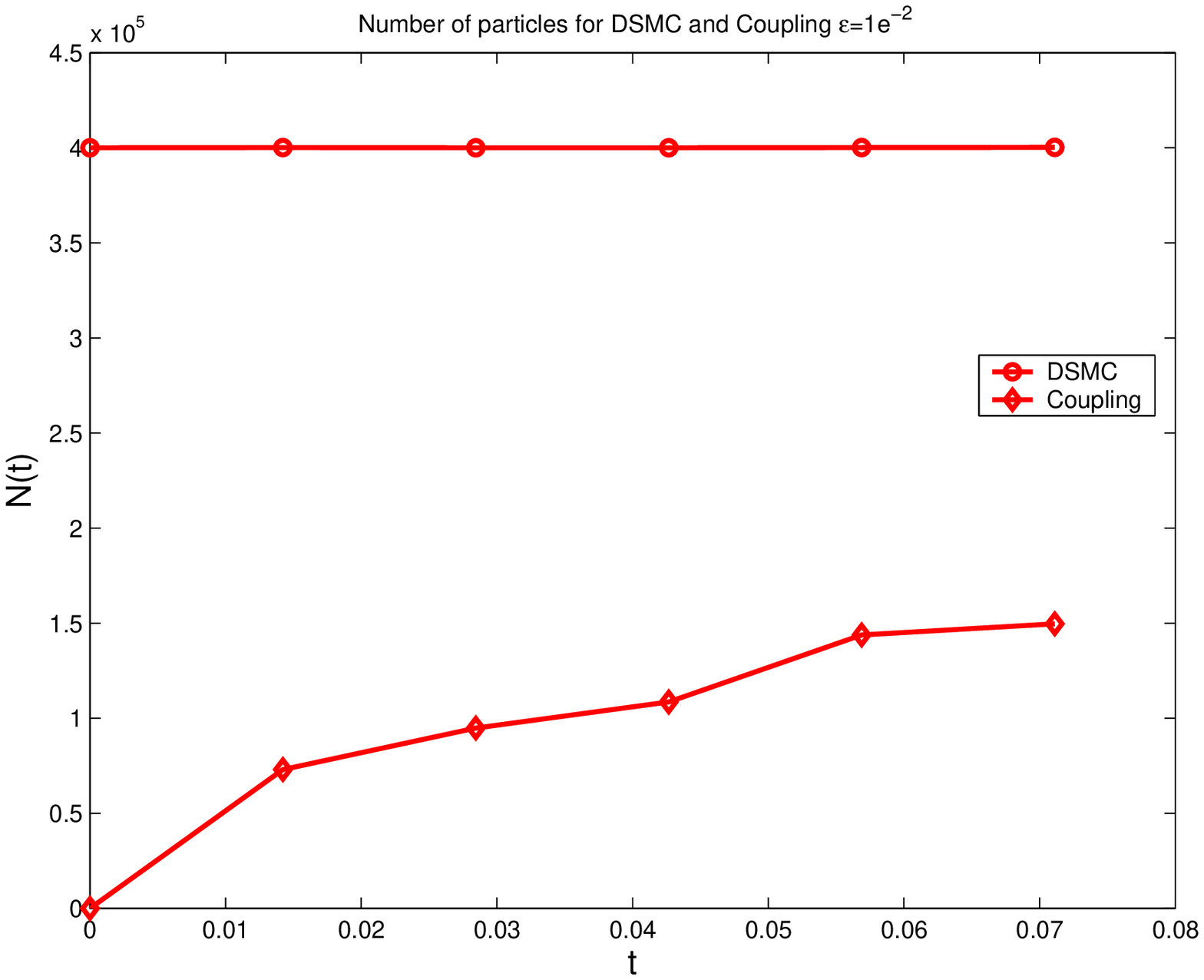}
\includegraphics[scale=0.27]{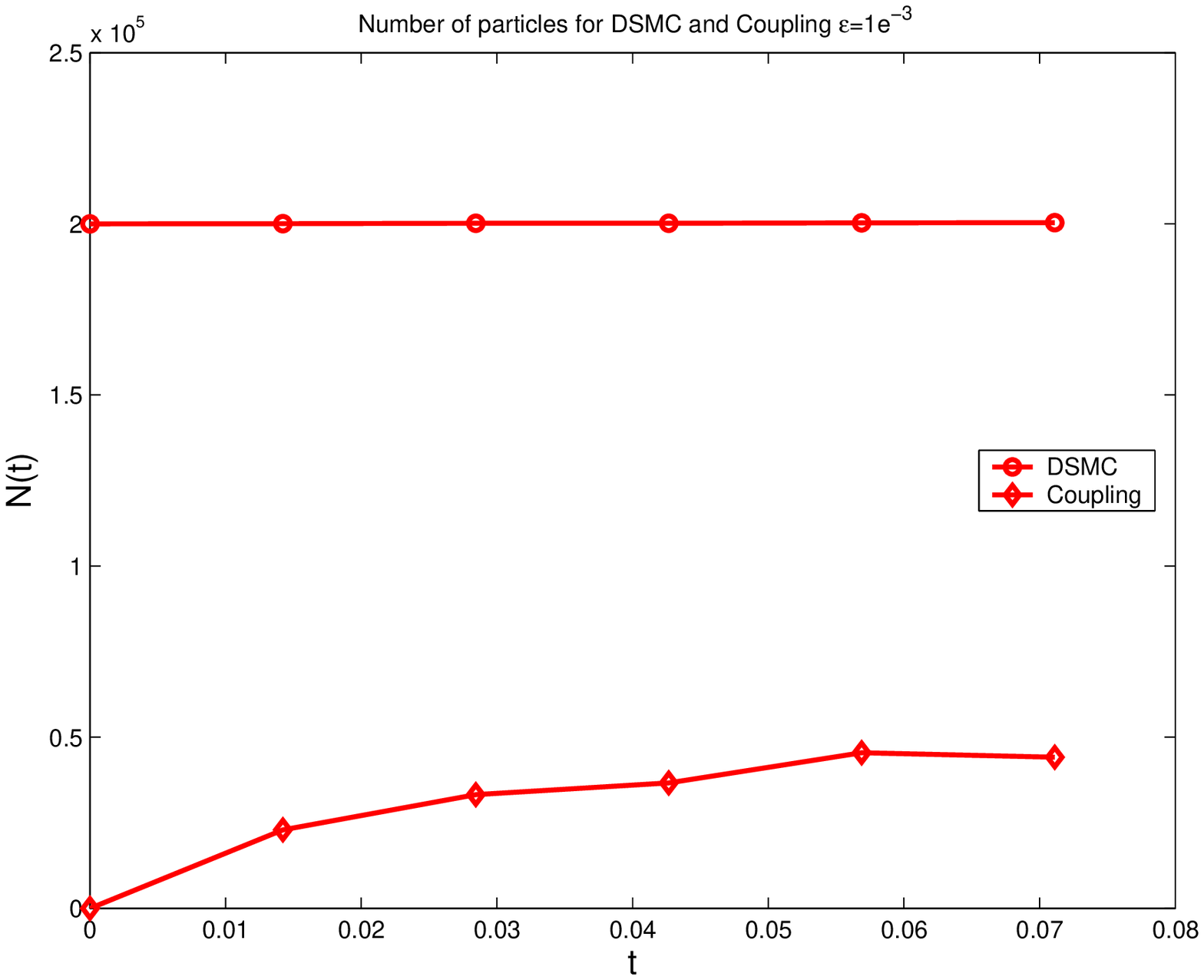}
\caption{Sod Test: Number of particle in time for the Monte Carlo
scheme and the DSMC-Fluid Coupling. Knudsen Number
$\varepsilon=10^{-1}$ (left), $\varepsilon=10^{-2}$ (middle),
$\varepsilon=10^{-3}$ (right).} \label{T24}
\end{center}
\end{figure}

\begin{figure}
\begin{center}
\includegraphics[scale=0.27]{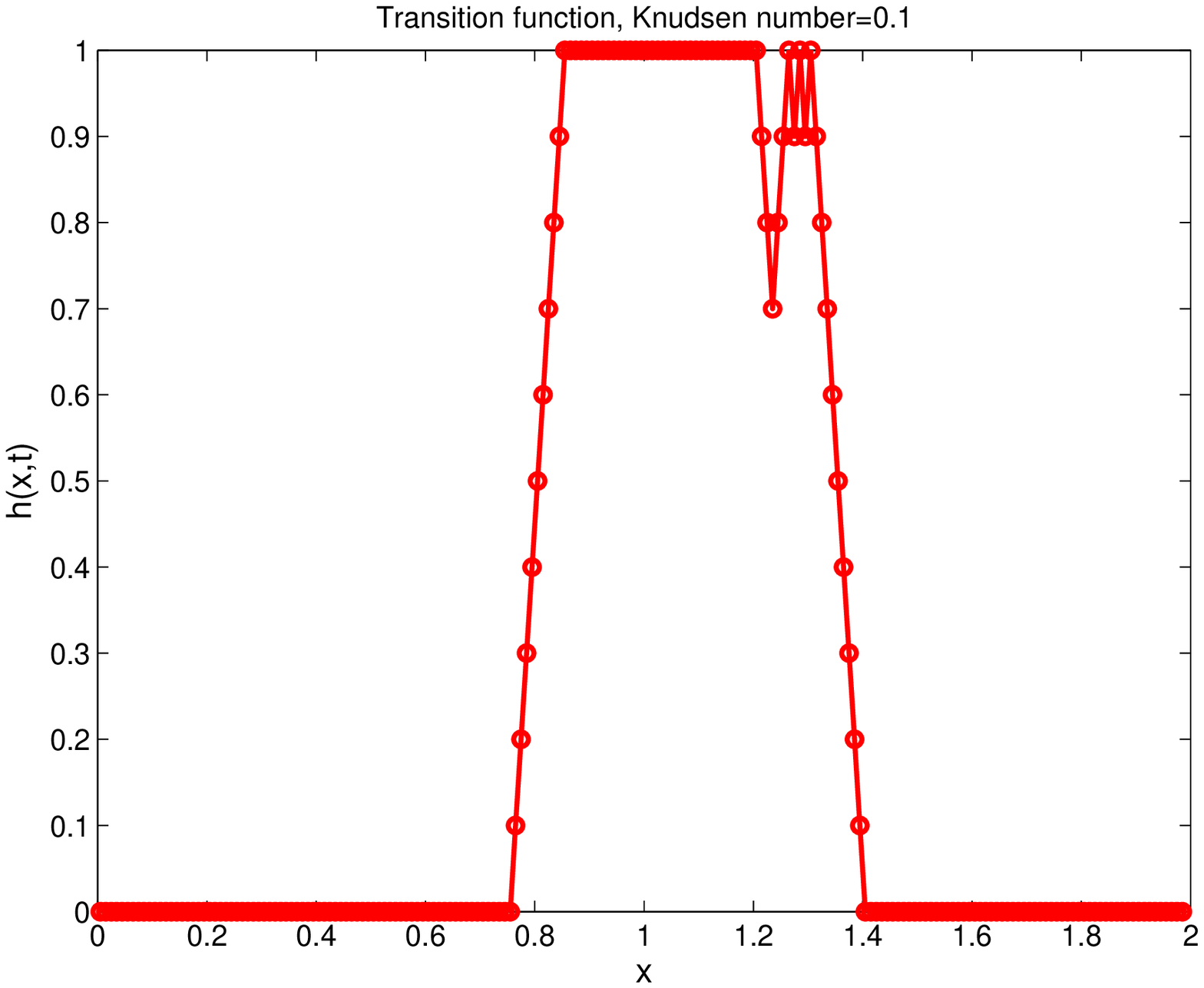}
\includegraphics[scale=0.27]{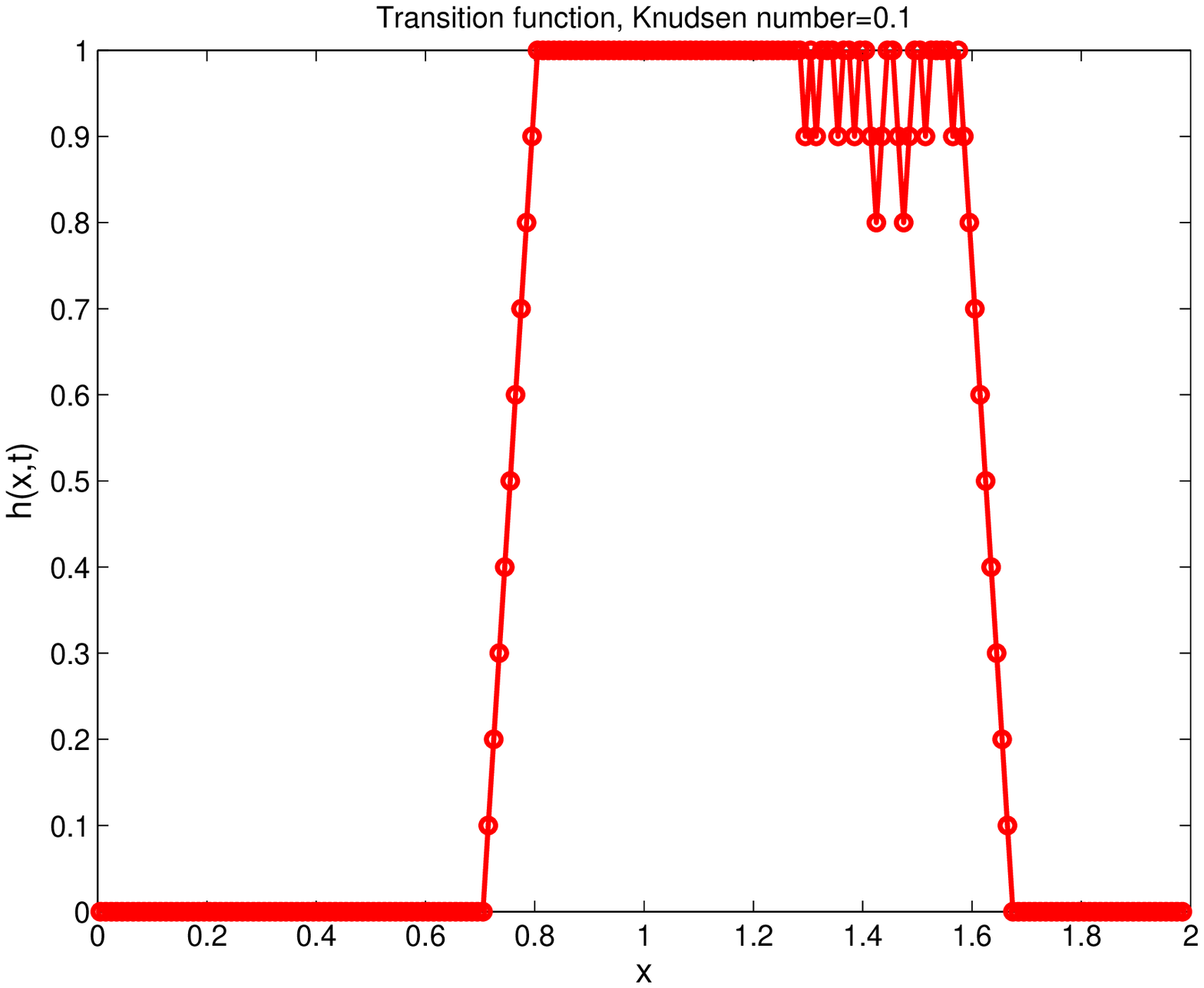}
\includegraphics[scale=0.27]{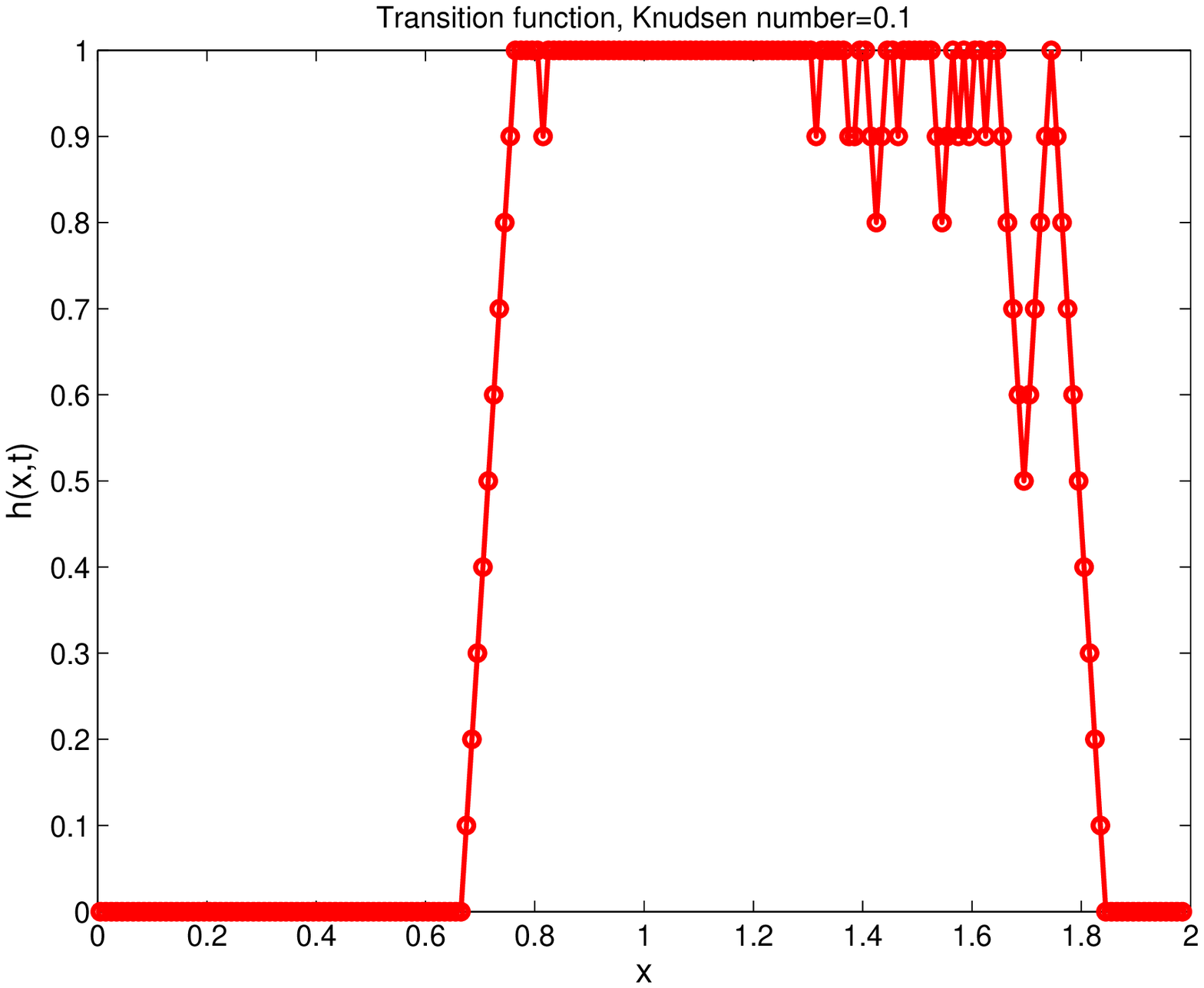}
\caption{Sod Test: Transition function at $t=0.3$ (left), $t=0.6$
(middle) and $t=0.8$ (right). Knudsen number $\varepsilon=10^{-1}$.}
\label{T03}
\end{center}
\end{figure}

\begin{figure}
\begin{center}
\includegraphics[scale=0.27]{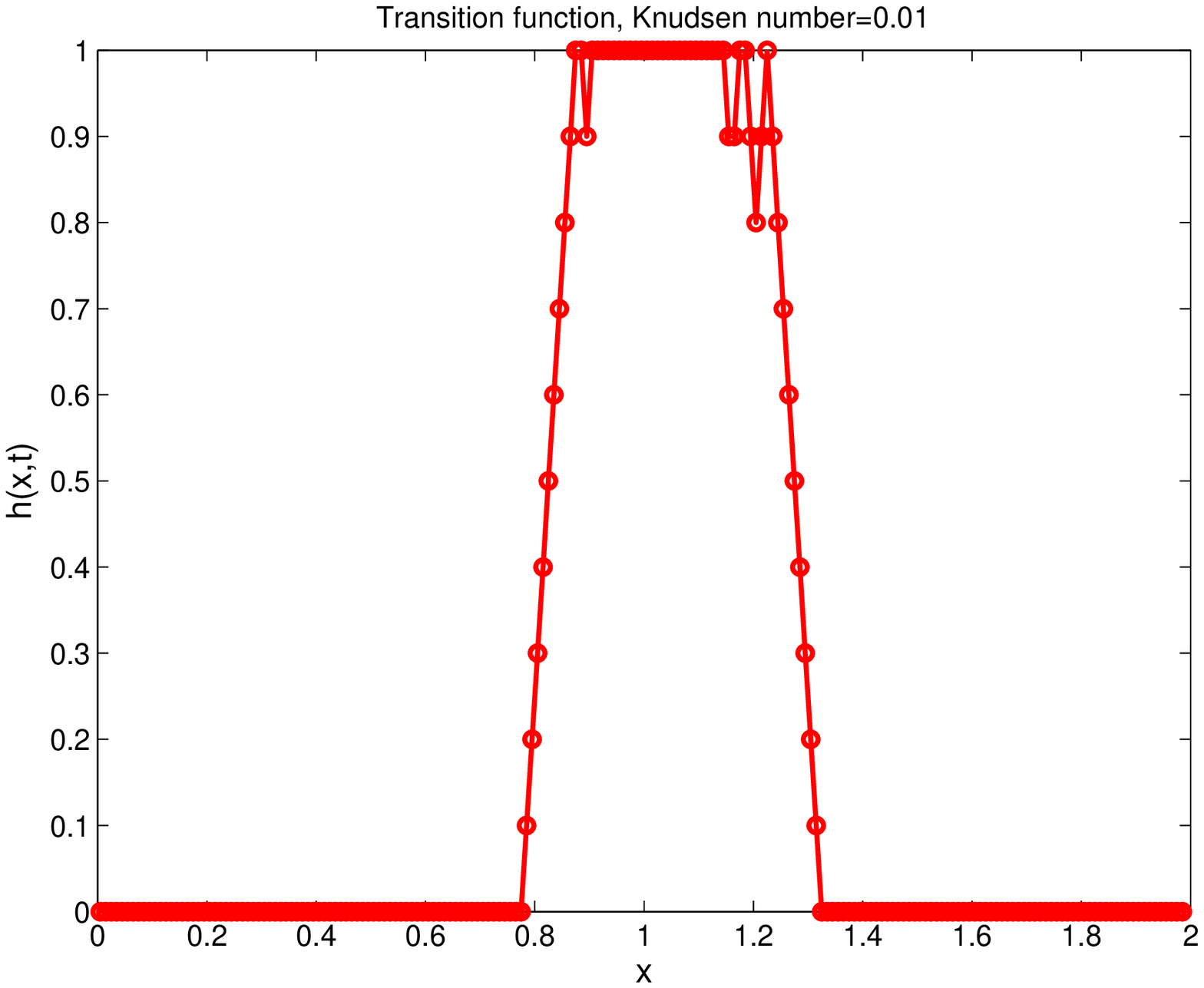}
\includegraphics[scale=0.27]{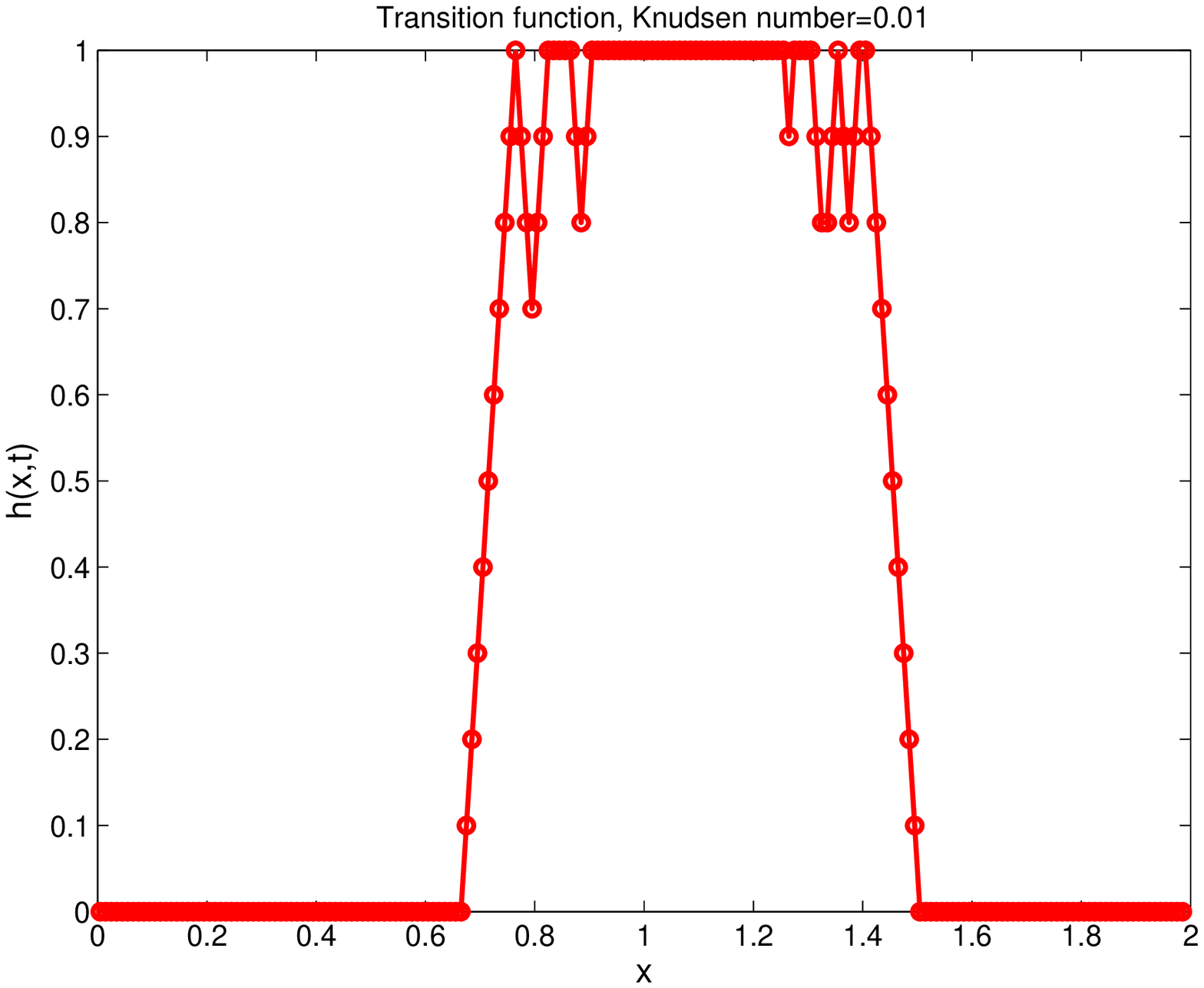}
\includegraphics[scale=0.27]{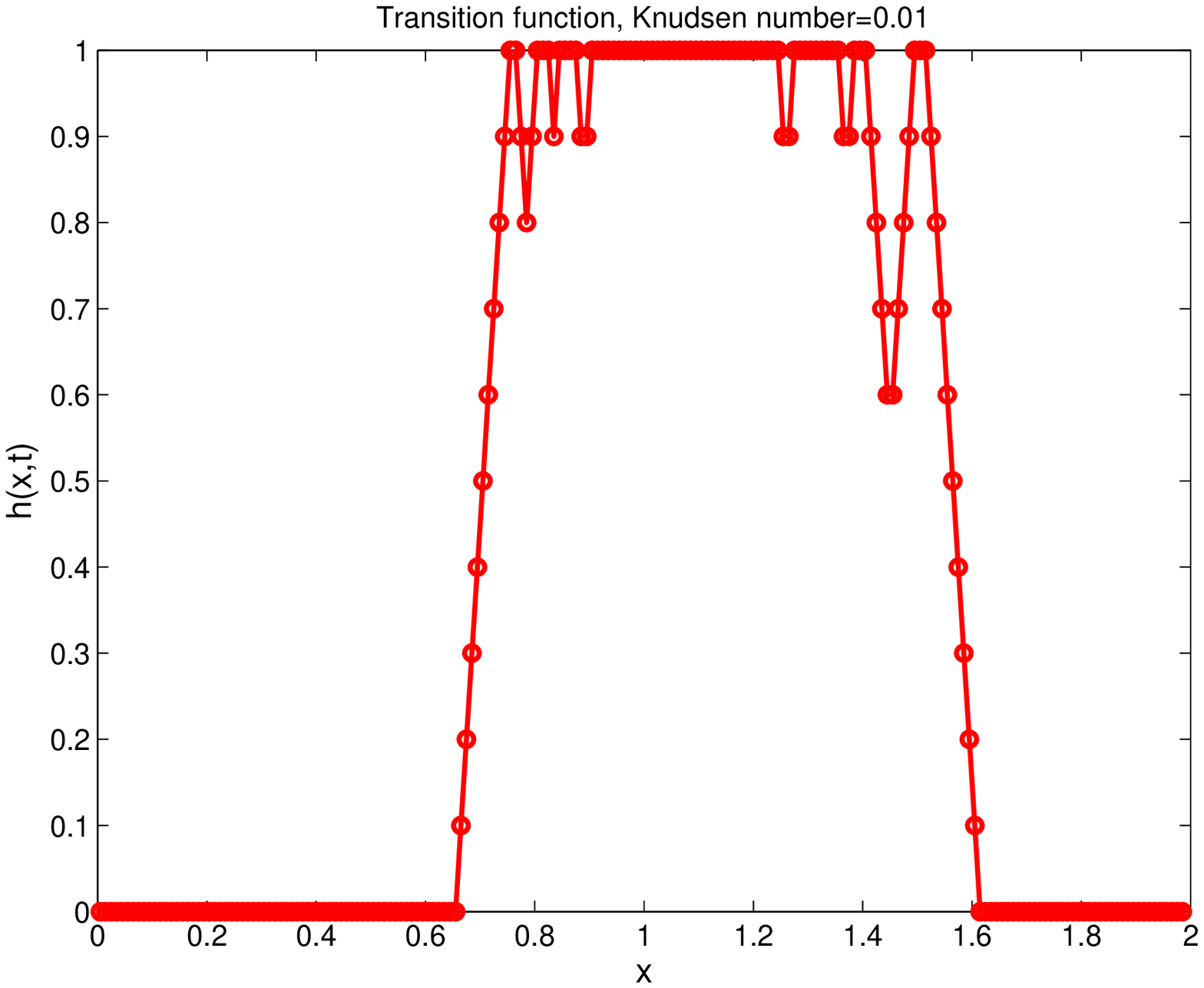}
\caption{Sod Test: Transition function at $t=0.3$ (left), $t=0.6$
(middle) and $t=0.8$ (right). Knudsen number $\varepsilon=10^{-2}$.}
\label{T13}
\end{center}
\end{figure}

In this final series of tests, we consider the classical Sod initial
data in a domain which ranges from $0$ to $2$. The numerical
parameters are the following: $200$ mesh points in space while the
number of particles is chosen different for different values of the
collision frequency. In this test case, we choose a larger number
particles, in comparison to the previous tests, to show that the
method is able to reproduce the correct profiles of the solutions.
The number of particles needed to keep the statistical noise
sufficiently low is a function of the perturbation $g_K$. In fact,
smaller is the perturbation, smaller is the statistical error in the
moment guided method. Observe anyway that, the classical DSMC with
the same number of particles still exhibits some important
fluctuations. The computation is stopped at the final time $t=0.8$.
Finally the simulations are initialized in the thermodynamic
equilibrium case, i.e. $h=0$ everywhere.

We take the following initial conditions: mass density
$\varrho_{L}=1$, mean velocity $u_{L}=0$ and temperature $T_{L}=5$
if $0 \leq x \leq 1 $, while $\varrho_{R}=0.125$, $u_{R}=0$,
$T_{R}=4$ if $1 \leq x \leq 2$. We repeat the same test changing the
collision frequency. They are respectively such that
$\varepsilon=10^{-1}$ in the first case, $\varepsilon=10^{-2}$ in
the second case and $\varepsilon=10^{-3}$ in the last case. The
thickness of the transition regions is fixed for every test
depending on the Knudsen: ten cells for large Knudsen ($10^{-1}$ and
$10^{-2}$) and five cells for small Knudsen ($10^{-3}$).

\begin{figure}
\begin{center}
\includegraphics[scale=0.27]{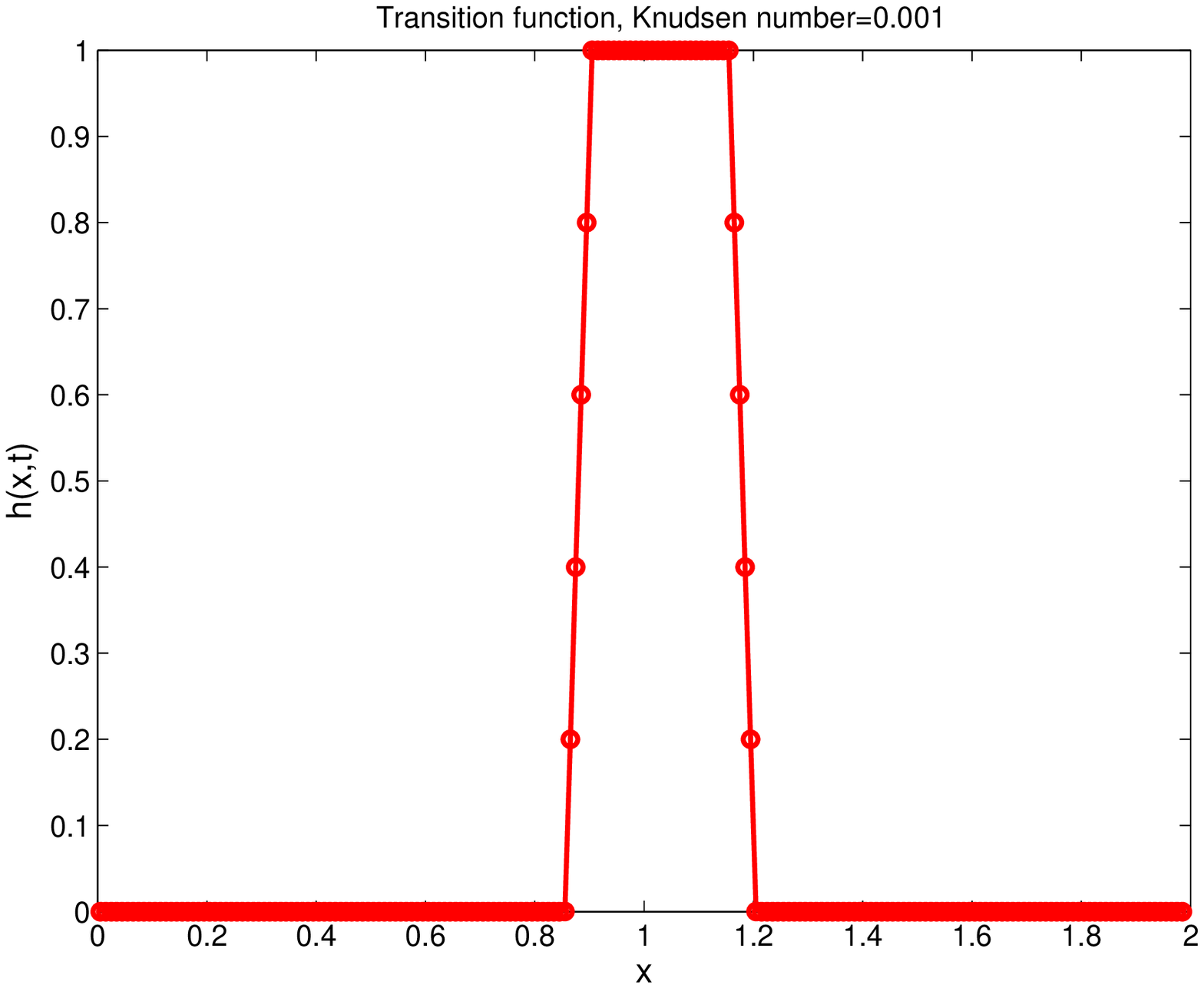}
\includegraphics[scale=0.27]{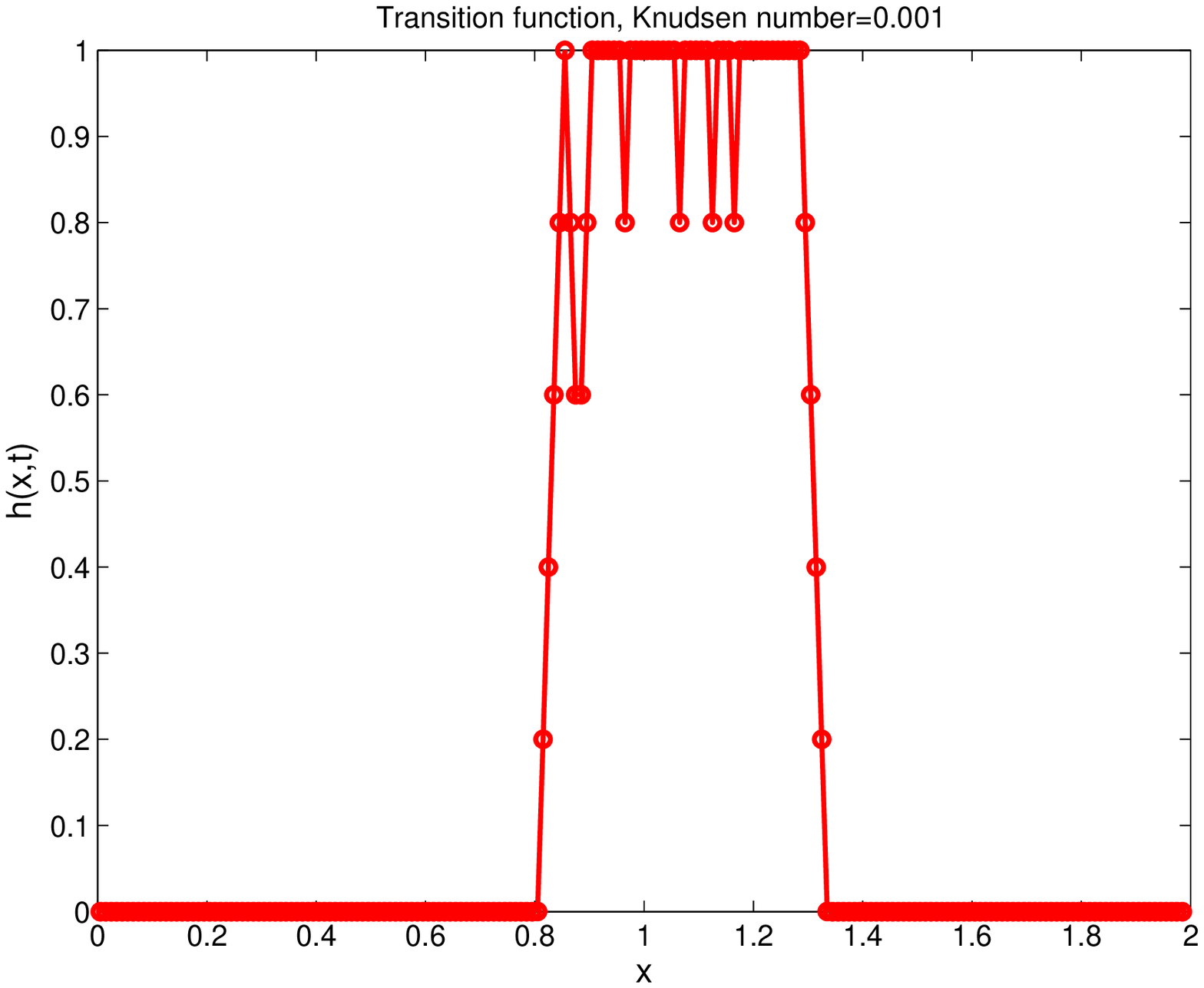}
\includegraphics[scale=0.27]{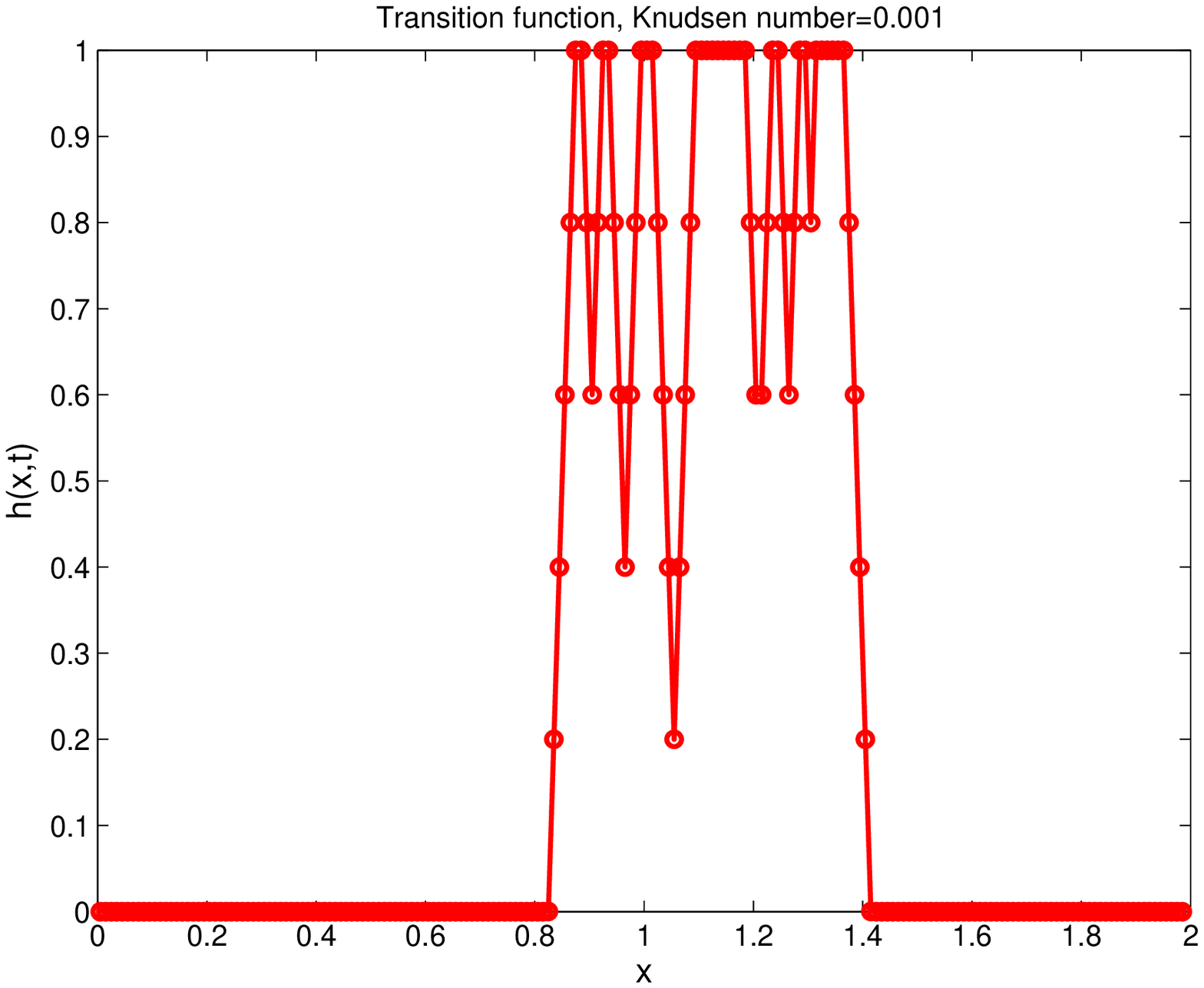}
\caption{Sod Test: Transition function at $t=0.3$ (left), $t=0.6$
(middle) and $t=0.8$ (right). Knudsen number $\varepsilon=10^{-3}$.}
\label{T23}
\end{center}
\end{figure}

In figure \ref{T24} we have reported the number of particles in time
for the DSMC method and the coupling method for different values of
the frequency $\varepsilon^{-1}$ ($\varepsilon=10^{-1}$ left,
$\varepsilon=10^{-2}$ middle, $\varepsilon=10^{-3}$ right). The
initial number of particles (i.e. the number of particles used to
represent the entire solution) is both for the DSMC and the coupling
$6 \ 10^{5}$ for $\varepsilon=10^{-1}$, $4 \ 10^{5}$ for
$\varepsilon=10^{-2}$ and $2 \ 10^{5}$ for $\varepsilon=10^{-3}$.
These numbers are a measure of the computational cost of the method
as explained in the previous paragraph. In figure \ref{T03} the
transition function is depicted for three different times in the
case of $\varepsilon=10^{-1}$. The same function is depicted in
figure \ref{T13} for $\varepsilon=10^{-2}$ and in figure \ref{T23}
for $\varepsilon=10^{-3}$. These figures shows the dynamics of the
kinetic and fluid regions in time and the capability of the method
to follow not only the shock but also regions in which the departure
from the equilibrium is smaller and which moreover moves in the
domain. We remark also that the thickness of the kinetic region
decreases when the collision frequency increases as expected. We
observe that the transition function oscillates in time. This is due
to the statistical fluctuations of the Monte Carlo scheme. If a more
precise evaluation of the zones far from equilibrium is needed, it
is sufficient to increase the number of particles. However, we point
out that even with an oscillatory behavior of the transition
function the profiles of the solution are very well captured.

\begin{figure}
\begin{center}
\includegraphics[scale=0.39]{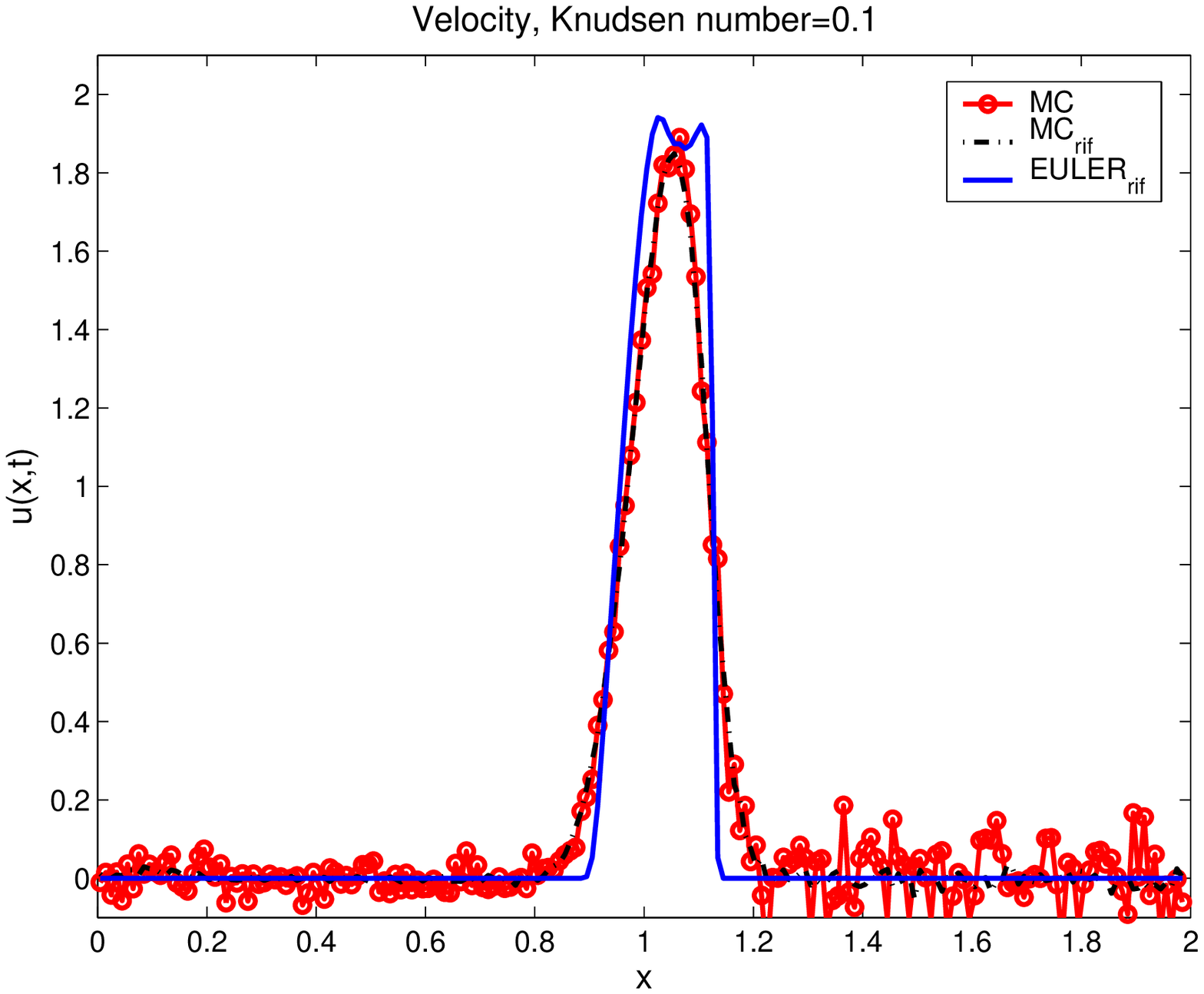}
\includegraphics[scale=0.39]{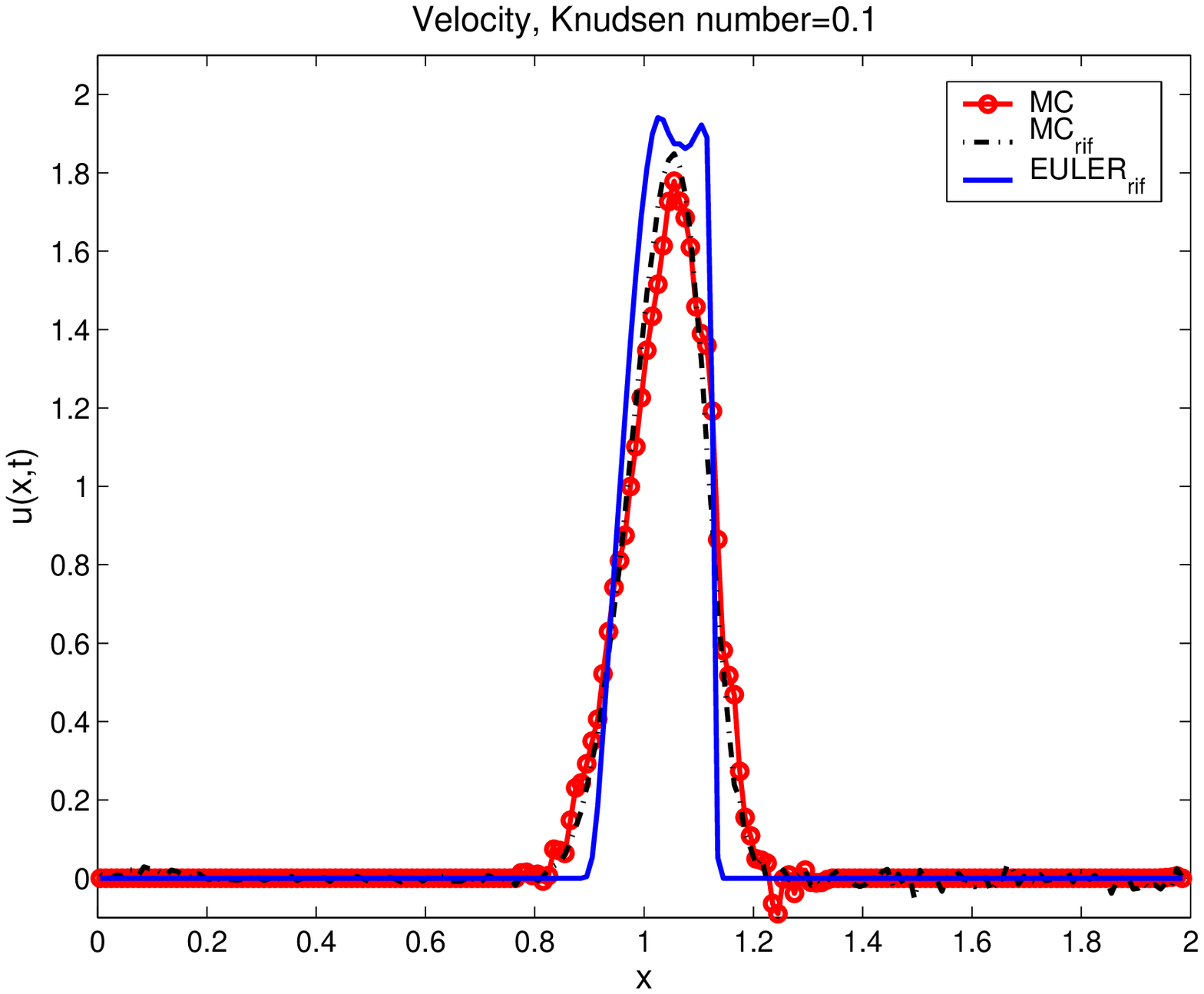}\\
\includegraphics[scale=0.39]{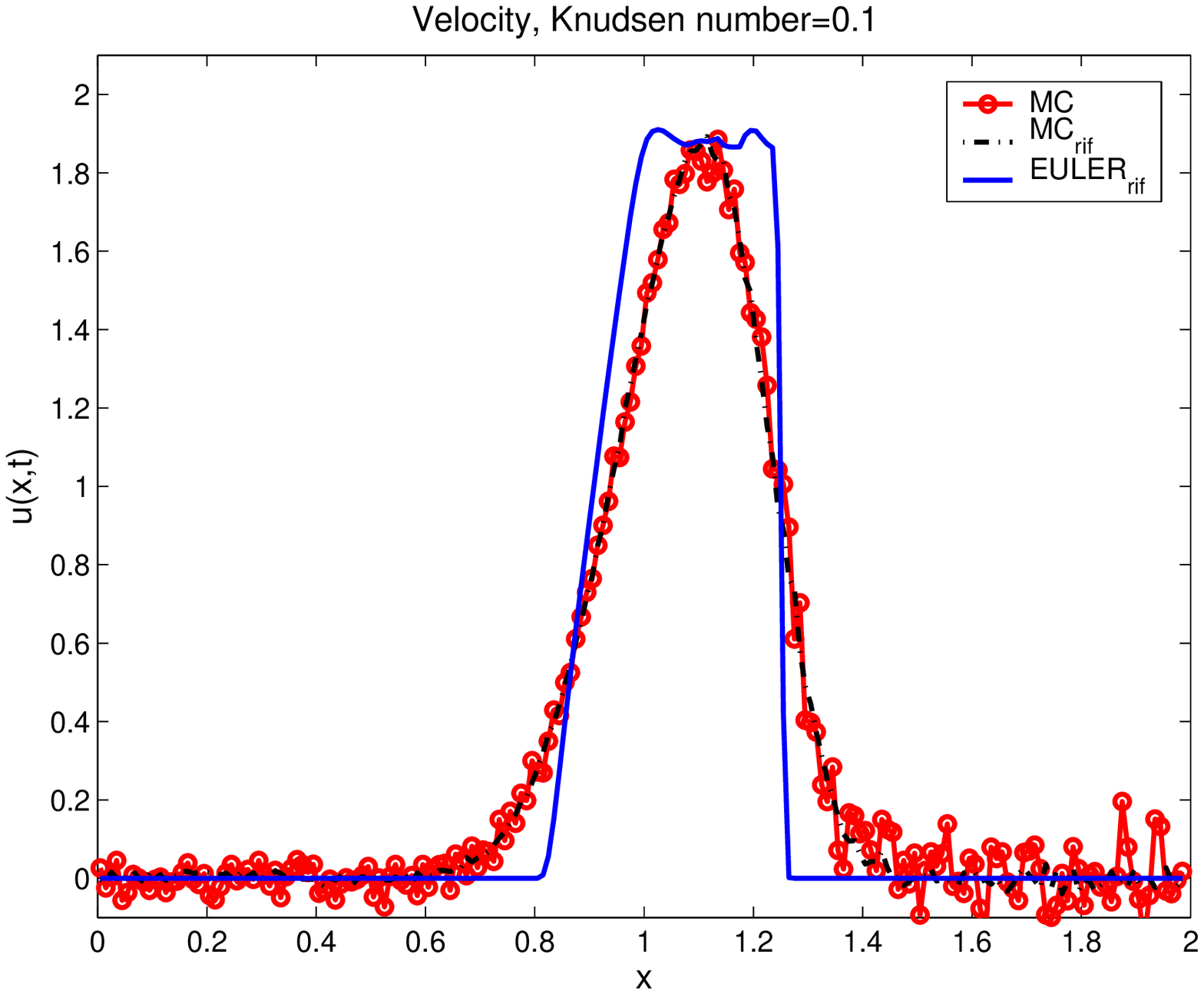}
\includegraphics[scale=0.39]{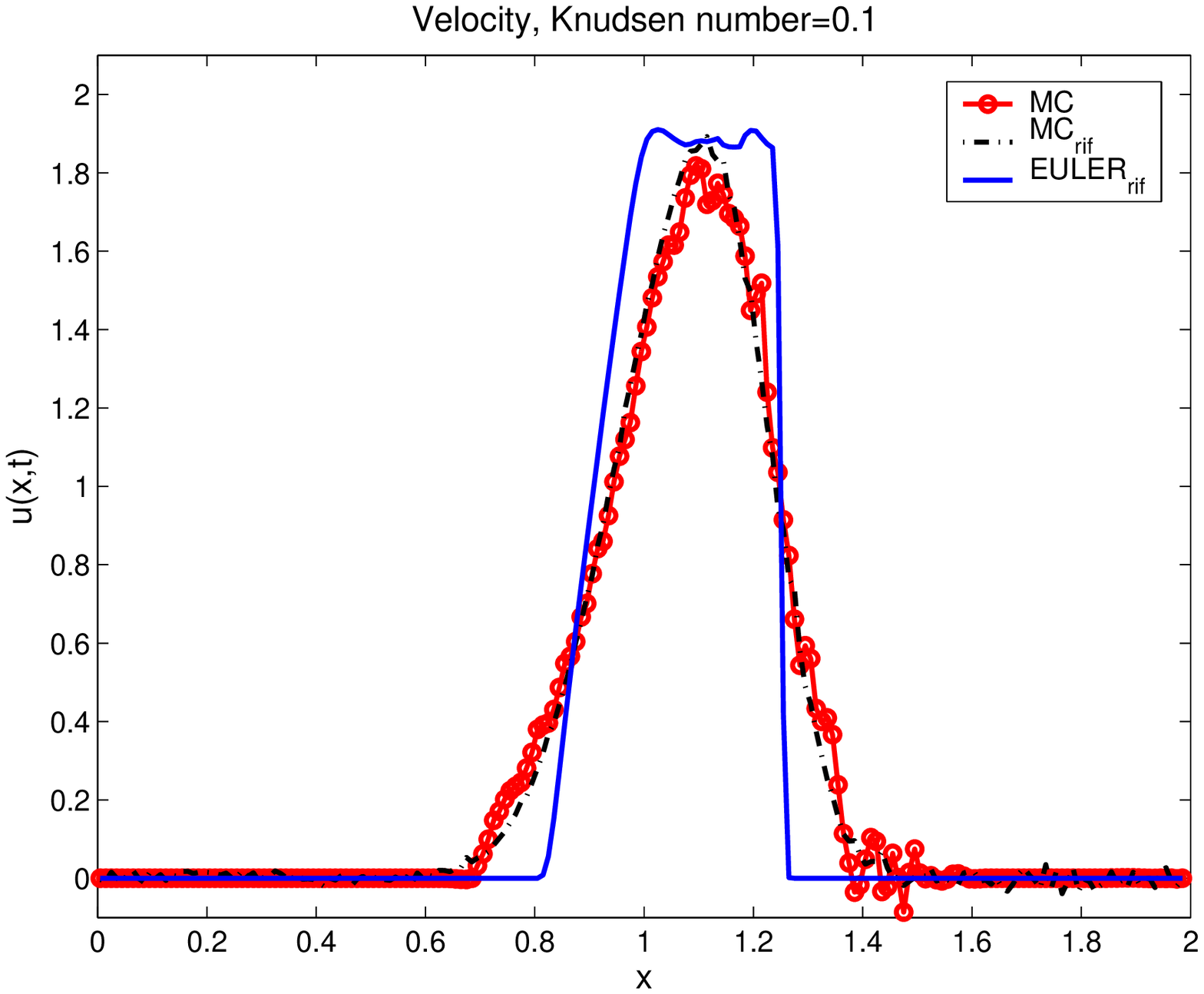}\\
\includegraphics[scale=0.39]{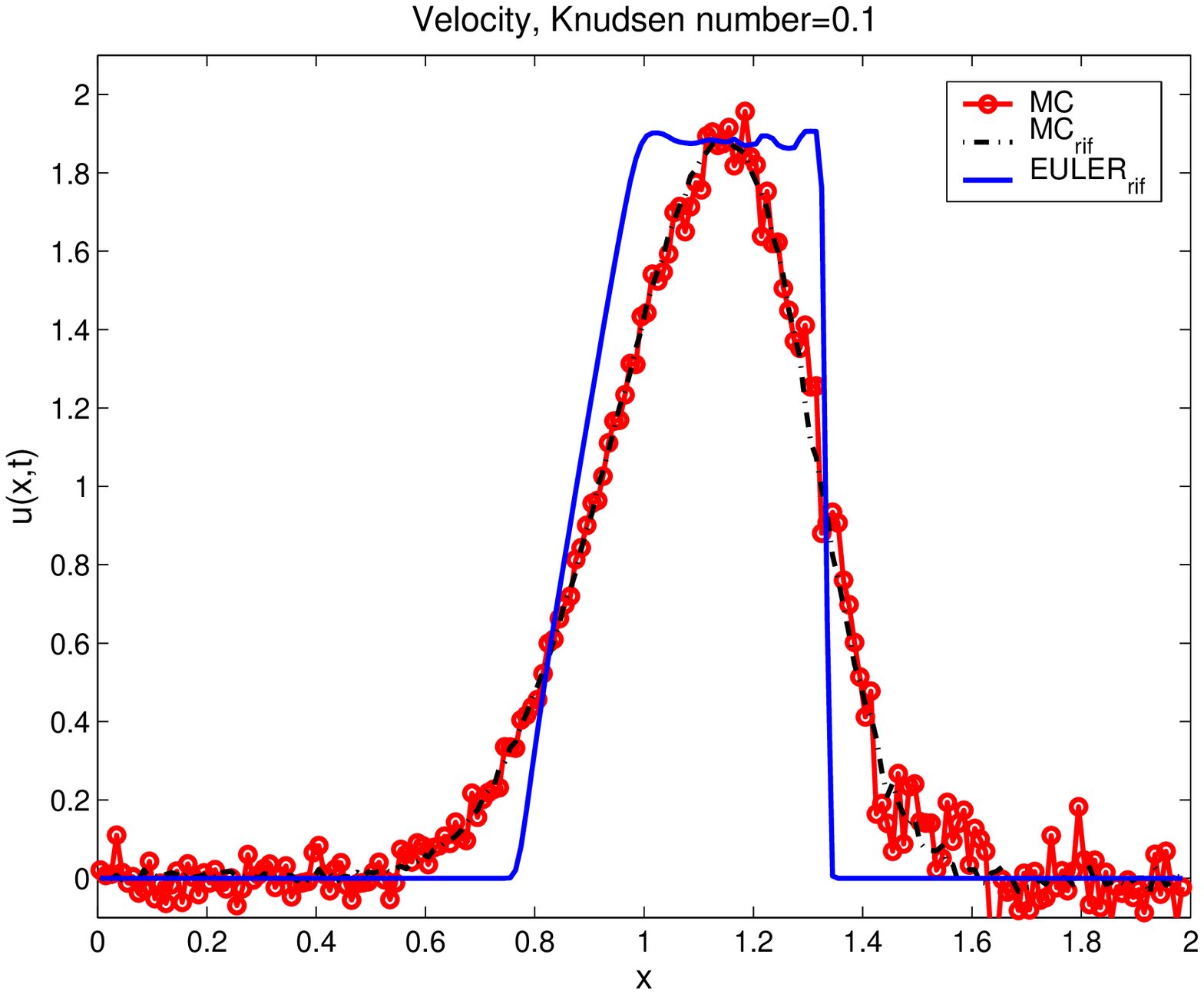}
\includegraphics[scale=0.39]{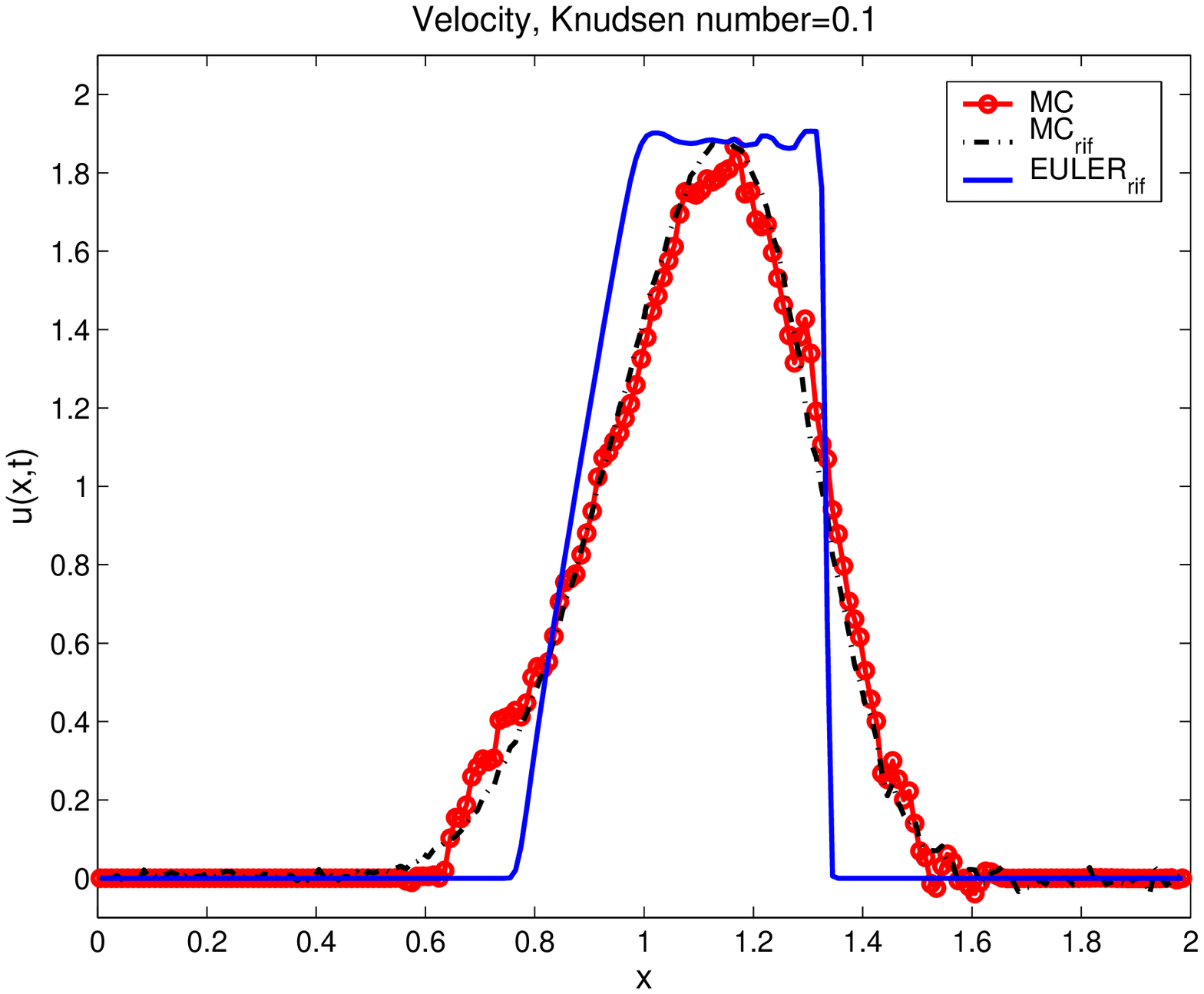}
\caption{Sod Test: Solution at $t=0.3$ (top), $t=0.6$ (middle) and
$t=0.8$ (bottom) for the velocity. MC method (left), Coupling
DSMC-Fluid method (right). Knudsen number $\varepsilon=10^{-1}$.
Reference solution (dotted line), Euler solution (continuous line),
DSMC-Fluid or DSMC (circles plus continuous line).}\label{T01}
\end{center}
\end{figure}

In figure \ref{T01} we have reported the velocity profile on the
left for the DSMC method and the on the right for the DSMC/Fluid
coupling. In this figure the collision frequency is taken equal to
$\varepsilon=10^{-1}$. From top to bottom, time increases from
$t=0.3 $ (top) to $t=0.8 $ (bottom) with $t=0.6$ in the middle. In
each of the plots regarding the macroscopic variables we reported
the solution computed with our coupling method or with the DSMC
method, a reference solution computed with a DMSC method where the
number of particles is taken very high to reduce the statistical
fluctuations and the solution of the compressible Euler equations.

\begin{figure}
\begin{center}
\includegraphics[scale=0.39]{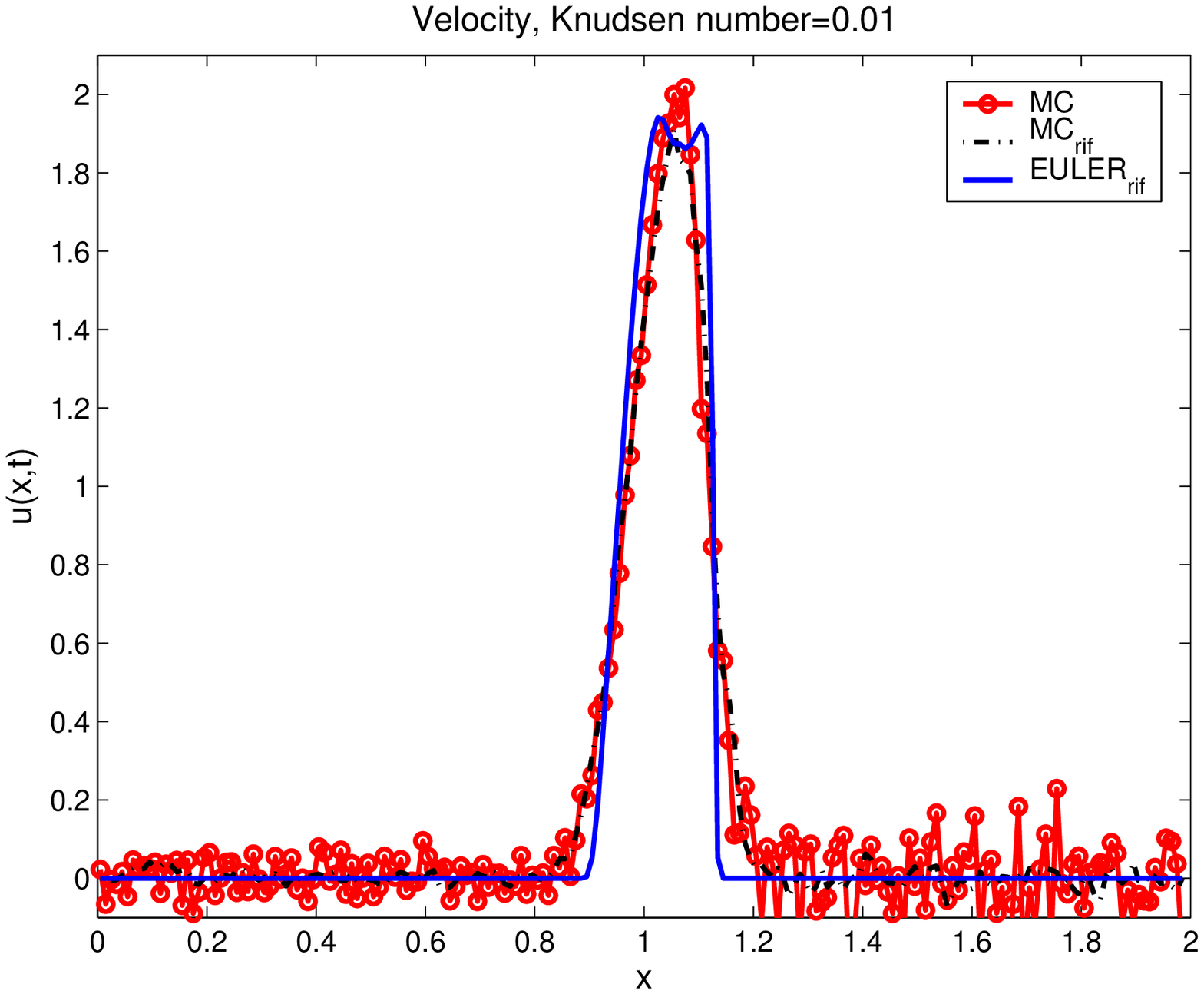}
\includegraphics[scale=0.39]{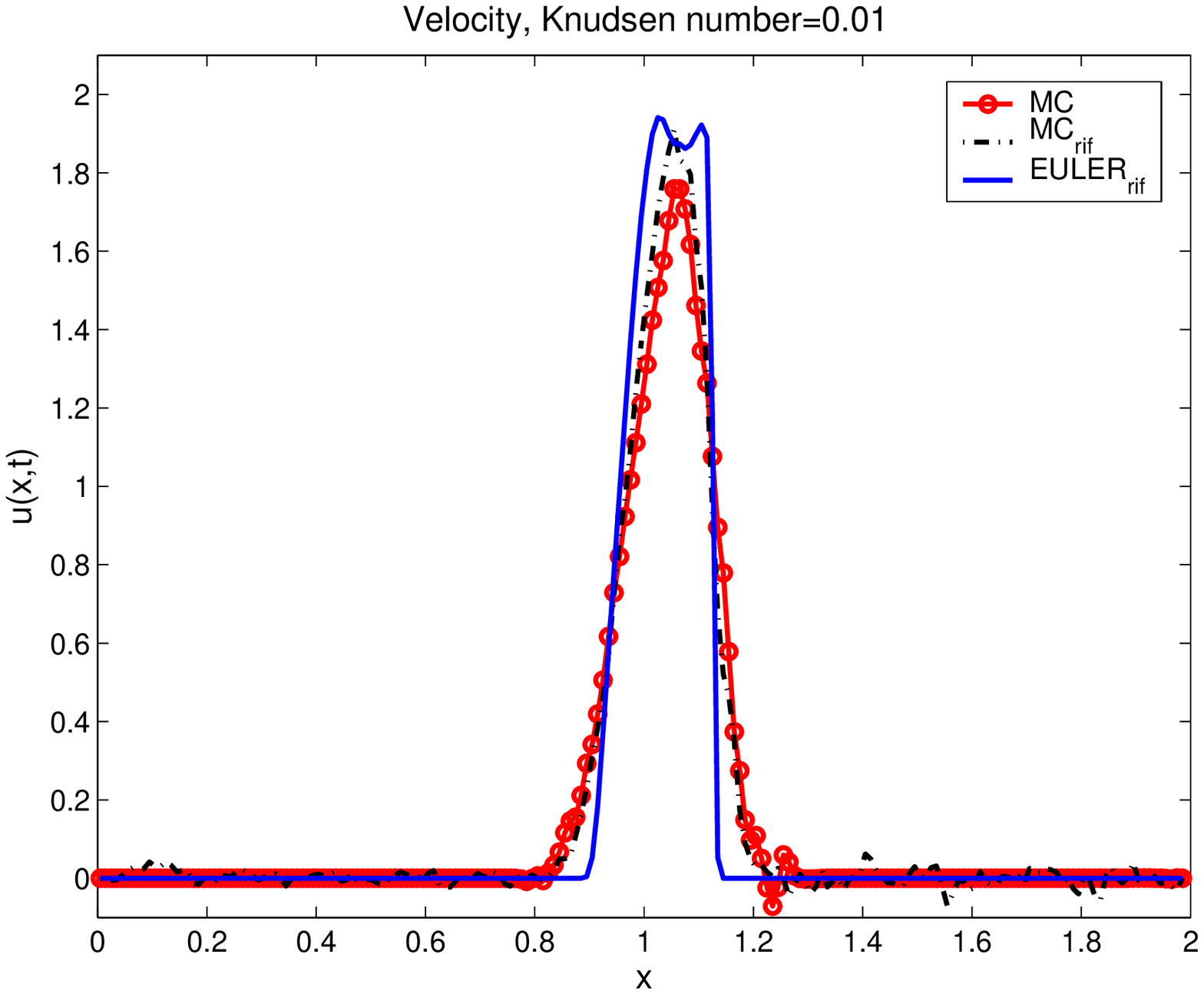}\\
\includegraphics[scale=0.39]{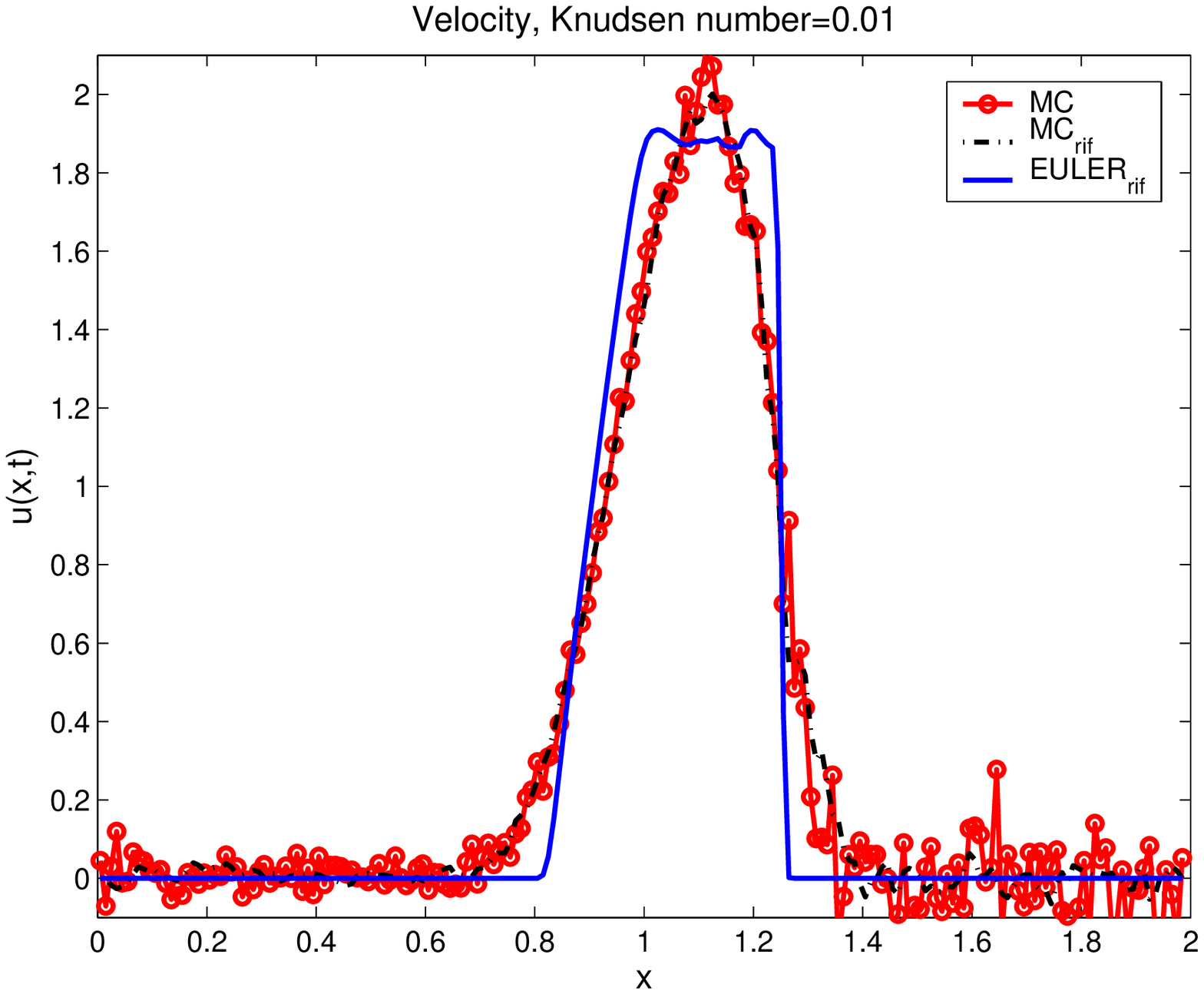}
\includegraphics[scale=0.39]{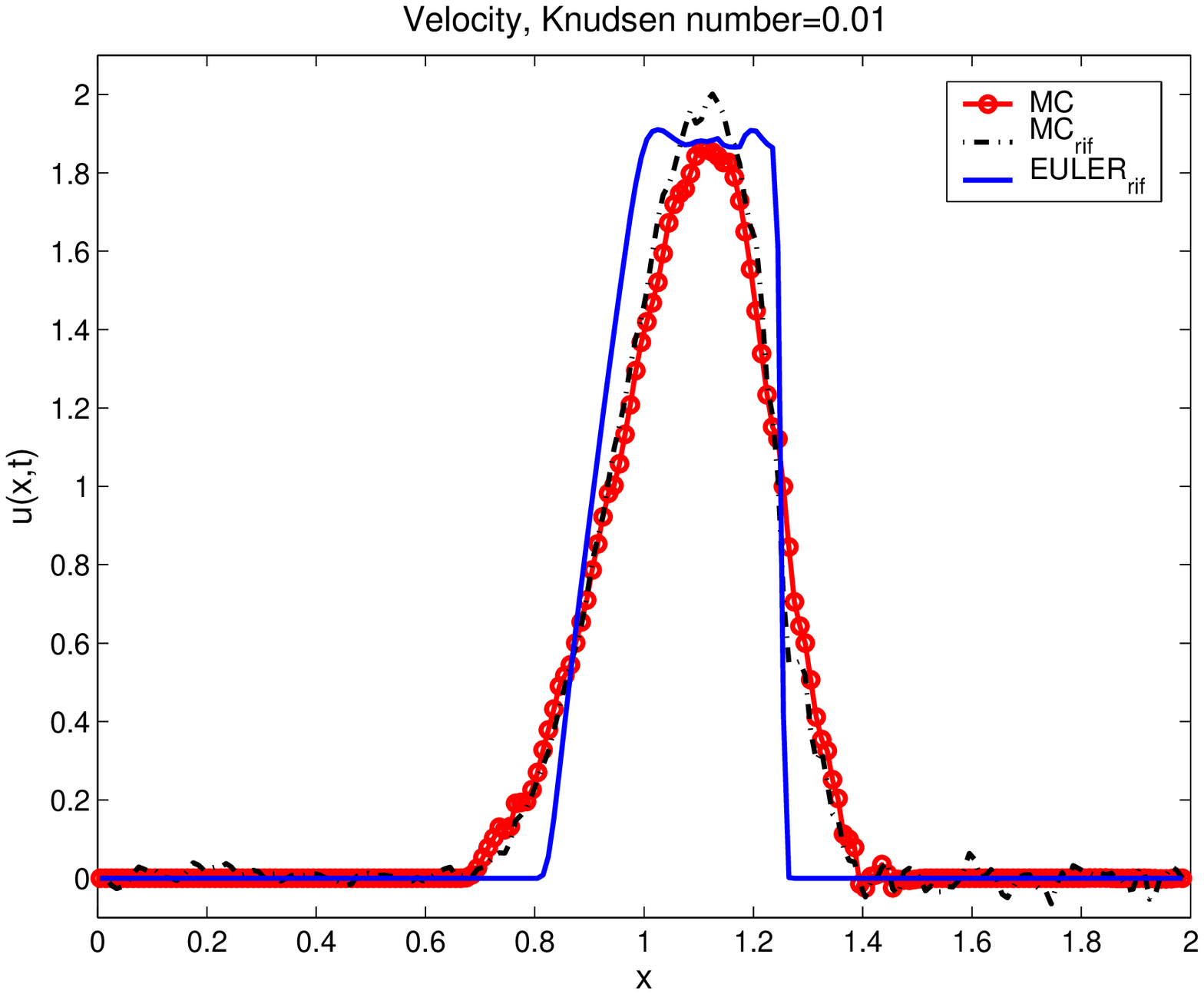}\\
\includegraphics[scale=0.39]{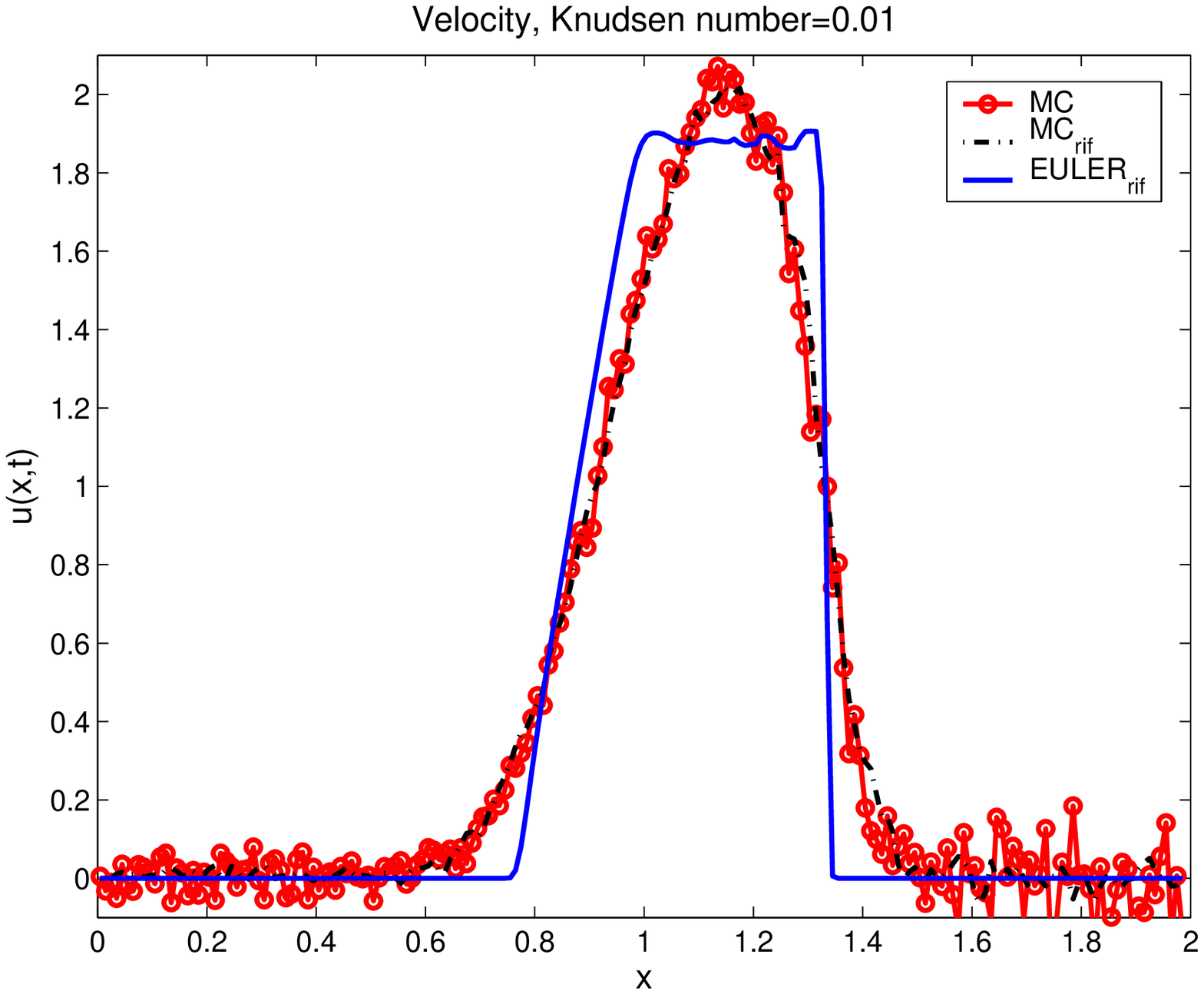}
\includegraphics[scale=0.39]{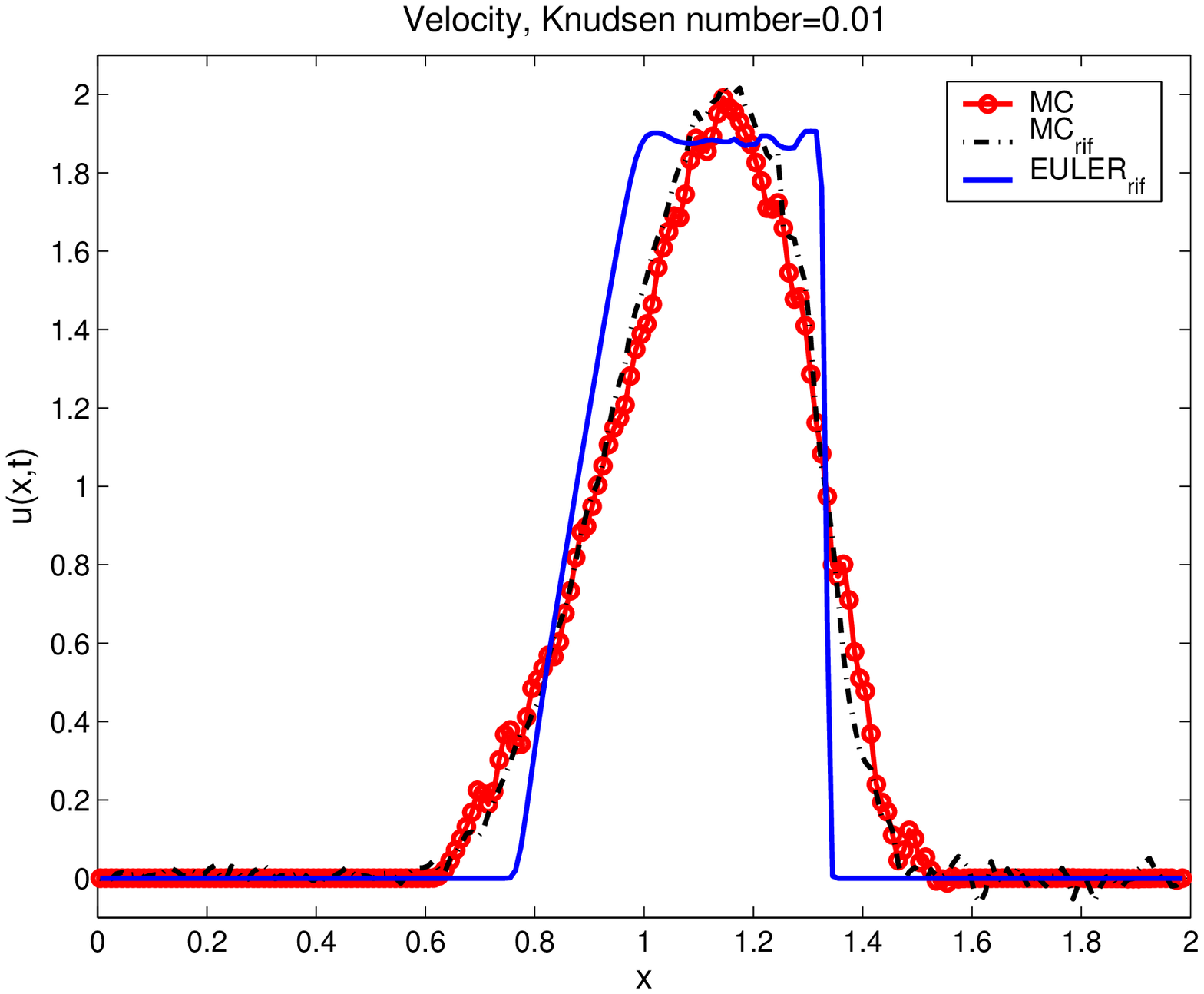}
\caption{Sod Test: Solution at $t=0.3$ (top), $t=0.6$ (middle) and
$t=0.8$ (bottom) for the velocity. MC method (left), Coupling
DSMC-Fluid method (right). Knudsen number $\varepsilon=10^{-2}$.
Reference solution (dotted line), Euler solution (continuous line),
DSMC-Fluid or DSMC (circles plus continuous line).}\label{T11}
\end{center}
\end{figure}

\begin{figure}
\begin{center}
\includegraphics[scale=0.39]{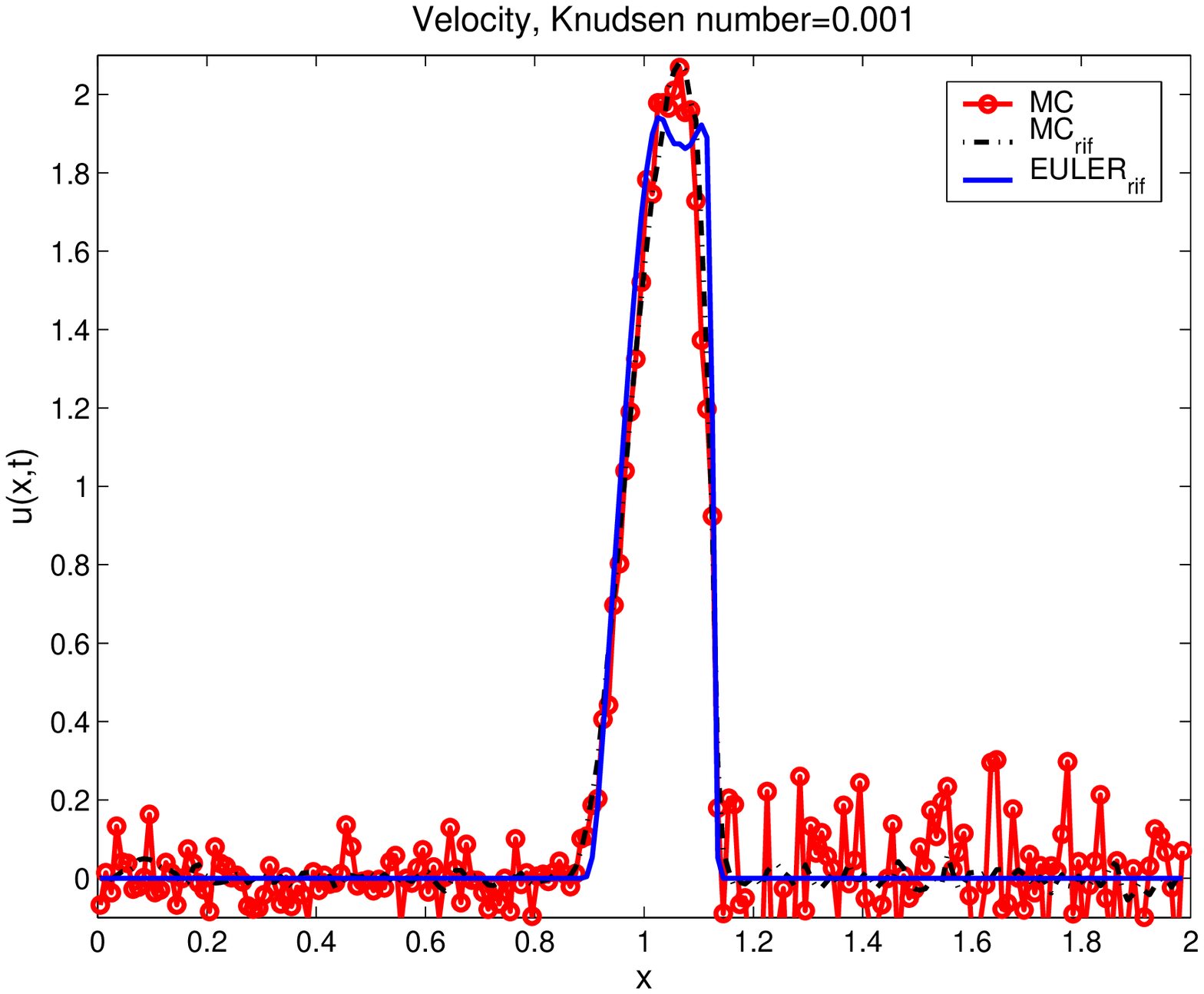}
\includegraphics[scale=0.39]{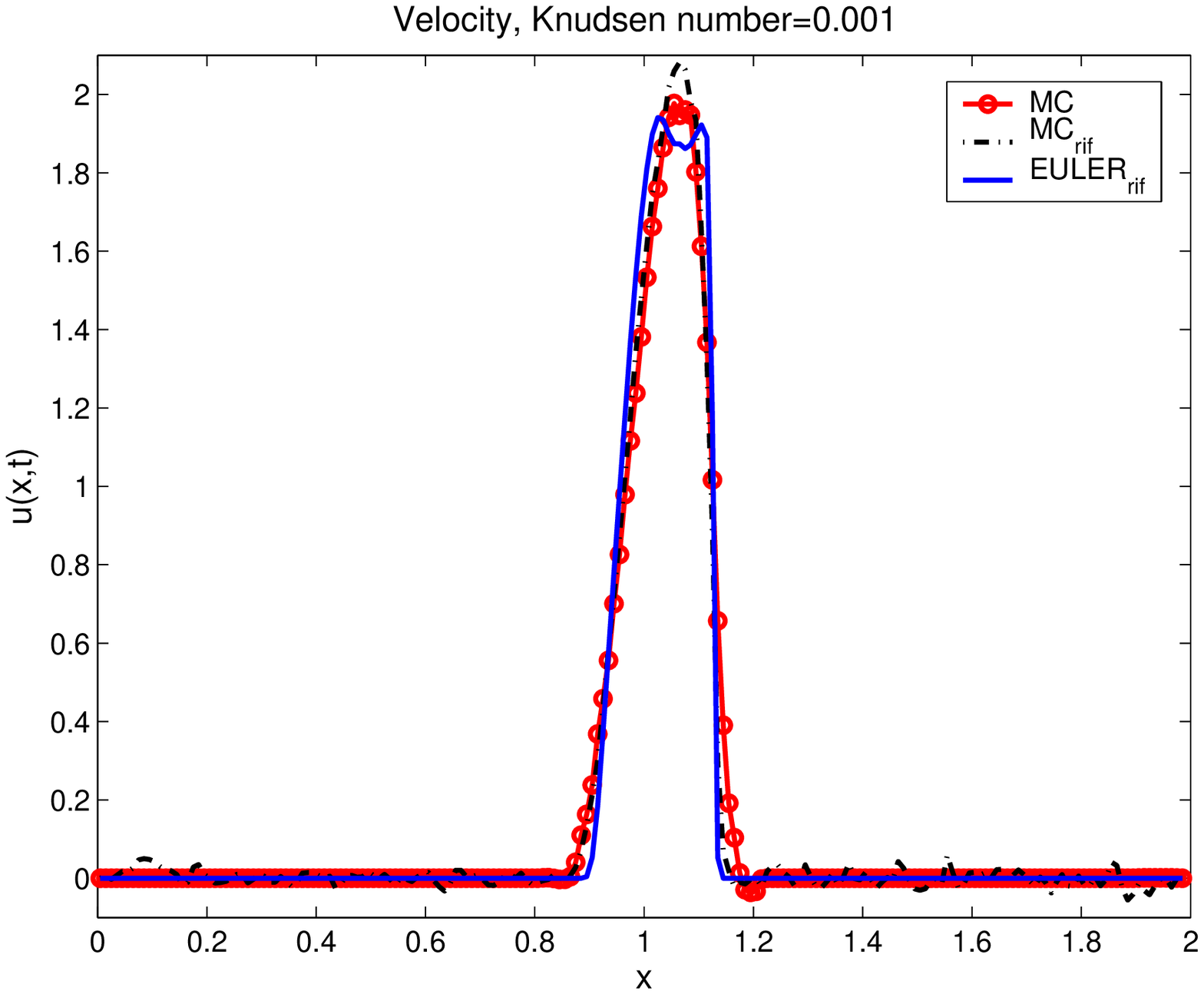}\\
\includegraphics[scale=0.39]{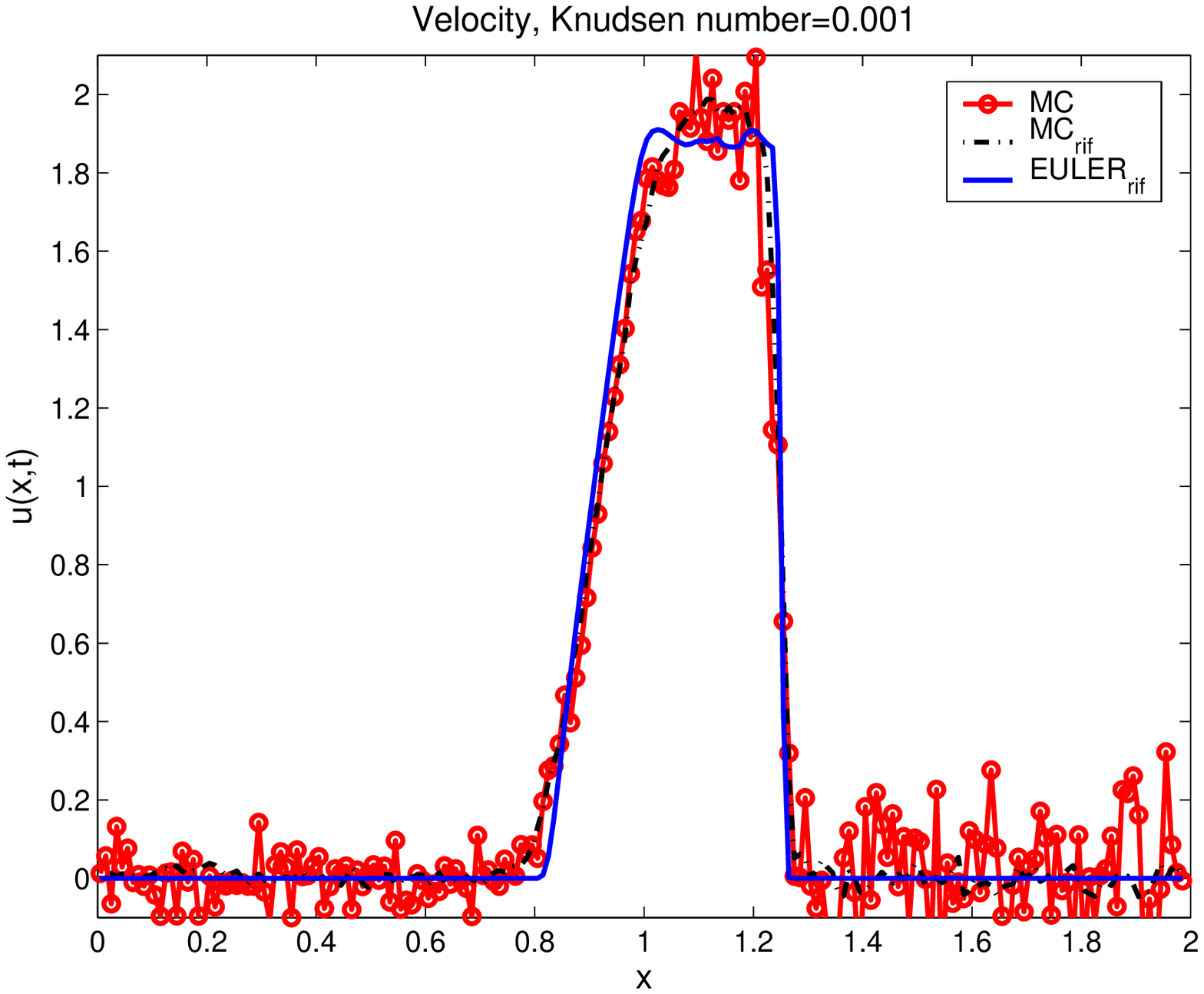}
\includegraphics[scale=0.39]{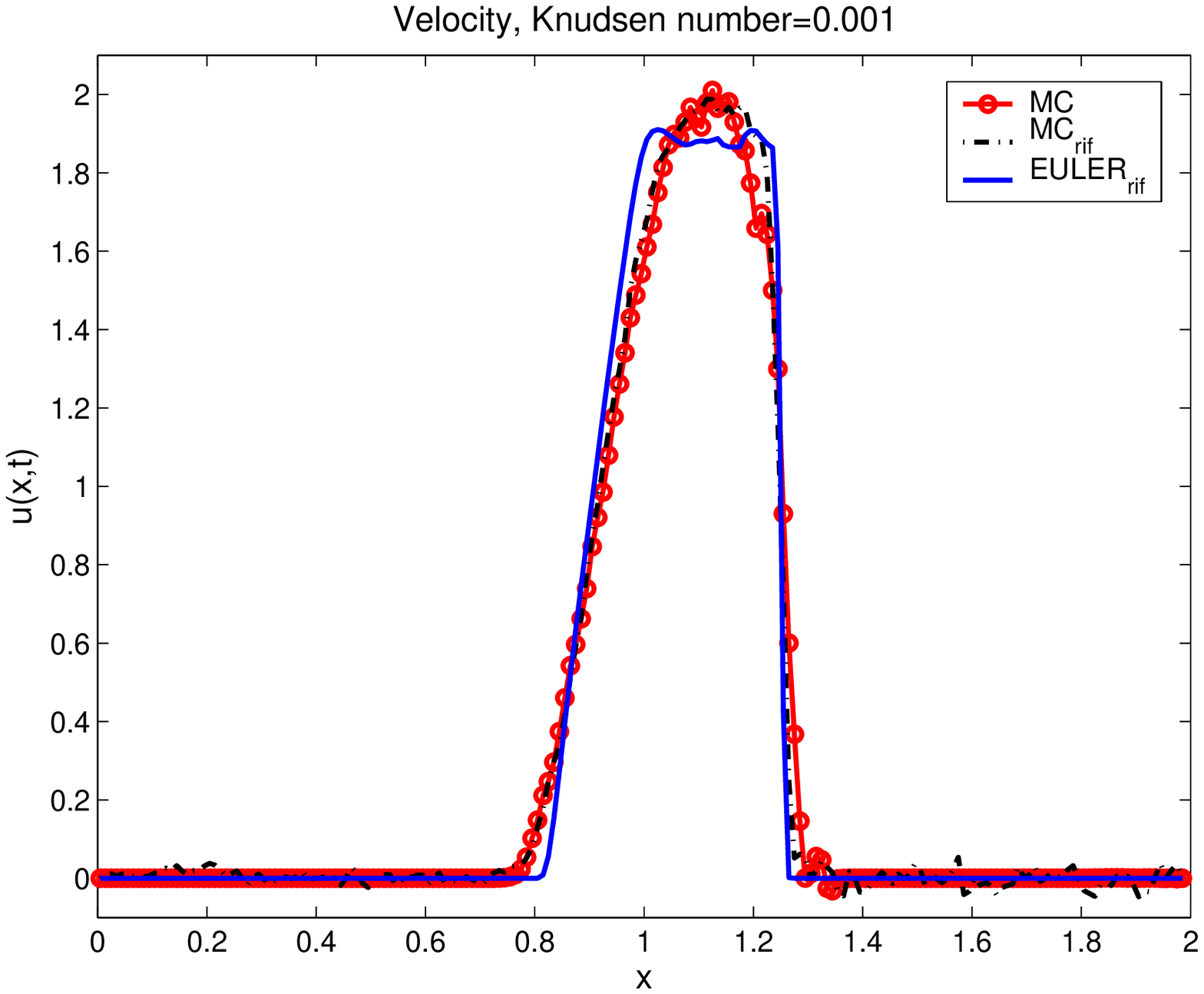}\\
\includegraphics[scale=0.39]{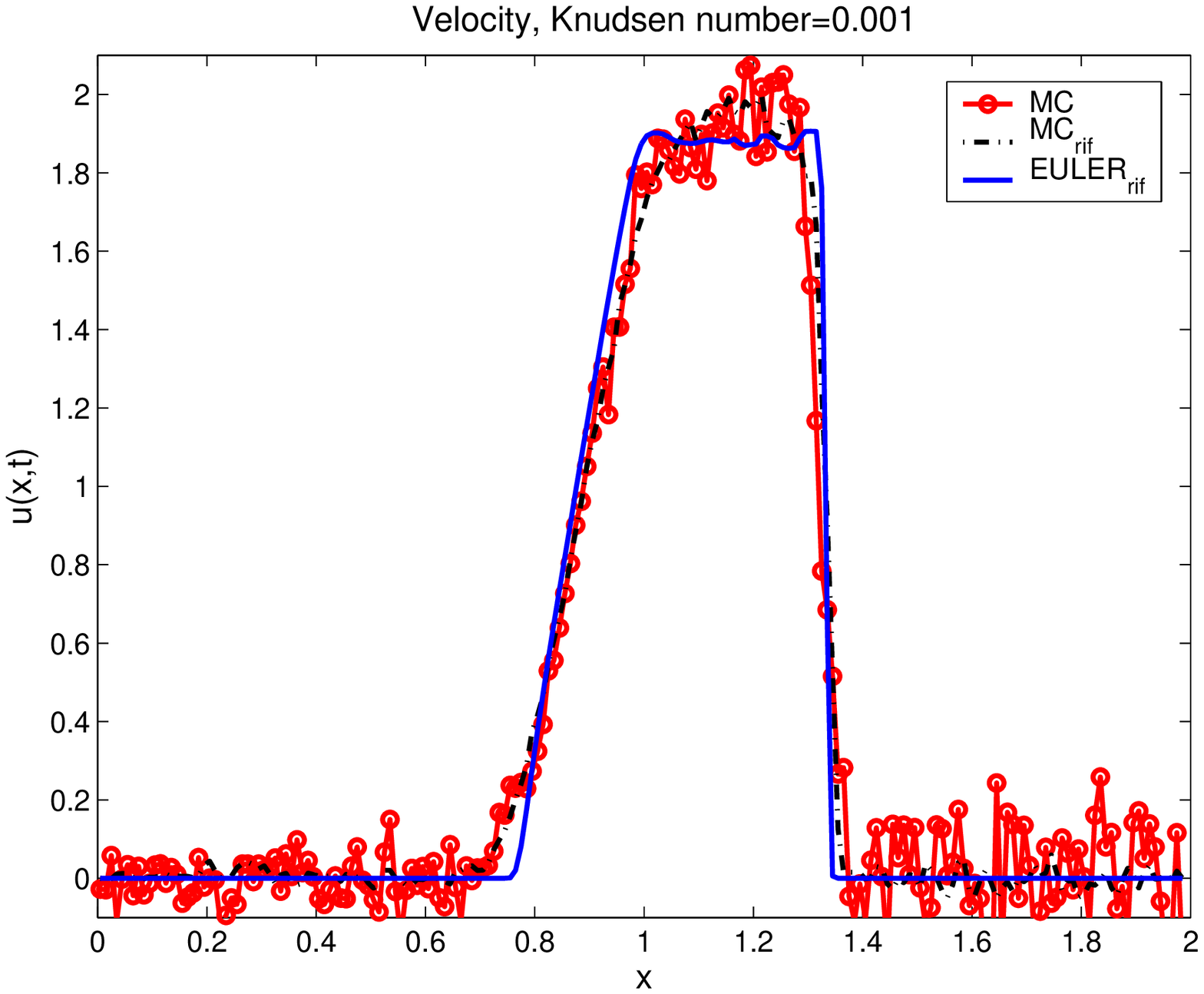}
\includegraphics[scale=0.39]{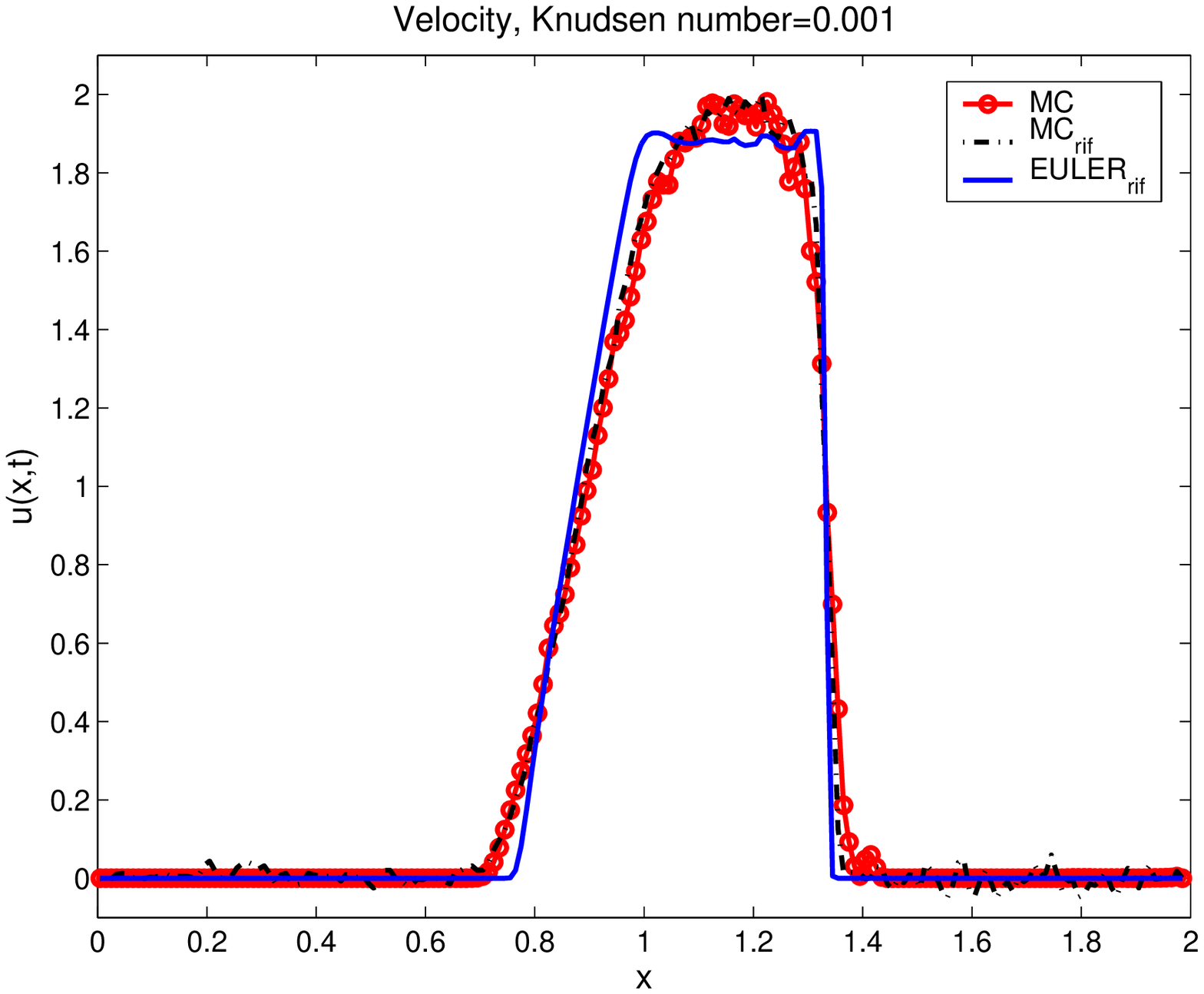}
\caption{Sod Test: Solution at $t=0.3$ (top), $t=0.6$ (middle) and
$t=0.8$ (bottom) for the velocity. MC method (left), Coupling
DSMC-Fluid method (right). Knudsen number $\varepsilon=10^{-3}$.
Reference solution (dotted line), Euler solution (continuous line),
DSMC-Fluid or DSMC (circles plus continuous line).}\label{T21}
\end{center}
\end{figure}

In figures \ref{T11} the velocity profile is reported for
$\varepsilon=10^{-2}$ again for DSMC on the left and DSMC/Fluid
coupling on the right. Finally we report the results in the case of
$\varepsilon=10^{-3}$ in figures \ref{T21} using the same
visualisation criterion of the previous cases.

In spite of the complexity of the solution (in this case we
have more than one zone which is potentially in non equilibrium), the
algorithm shows a very good behavior. In particular, in the last
case, the non equilibrium region is very tiny. Therefore, we largely reduce the
computational cost and we do not lose accuracy in the description of
the problem.

\section{Conclusion}
\label{sec_conclu}

In this paper, we have presented a novel numerical algorithm to
couple a DSMC method for the solution of the Boltzmann equation and
a finite volume like method for the compressible Euler equations.
The method is based on the introduction of a buffer zone which
realizes a smooth transition between the kinetic and fluid regions
which was first introduced by Degond, Jin and Mieussens in
\cite{degond2}. Succesively, this method was extended to the case of
moving interfaces in \cite{dimarco3} and in \cite{dimarco4}.

Here, we have extended the idea of buffer zones and dynamic coupling
to the case in which Monte Carlo methods are used to solve the
kinetic equation. This extension permits to treat more detailed
physical models than in the past. In particular it permits to
resolve the Boltzmann collision integral instead of simplified
relaxation models. In order to better adapt the coupling methodology
to the case of DSMC methods, the coupling terms have been rearranged
in a different way respect to the past which make the coupling more
suitable for Monte Carlo approximations. Moreover, to reduce the
fluctuations caused by Monte Carlo methods which produce spurious
oscillations in the fluid regions, we used a new technique which
permits to reduce the variance of particle methods \cite{dimarco1}.
This technique, called moment guided, is based on matching the
moments of the kinetic solution with the moments of a suitable set
of macroscopic equations which contains reduced statistical error at
each time step of the computation. In addition, the use of this
method permits a more precise estimation of the breakdowns of the
fluid model. This is true even when few particles, in comparison
with the typical number of particle employed in unsteady
computations, are used.

The last part of the present work is centered on the analysis of
several numerical tests. The results clearly demonstrate the
capability of the method to couple Monte Carlo method with finite
volume methods even for unsteady problems and with a relatively
small number of particles. The resulting method is able to
automatically create, cancel and move as many kinetic, fluid or
buffer regions as necessary. It is also able to capture the correct
solutions, to reduce the computational cost with respect to full
DSMC simulations and to avoid the propagation of fluctuations in the
fluid regions.

We finally observe that the computational speed-up will
significantly increase for two or three dimensional simulations,
which we intend to carry out in the future.



\end{document}